\documentclass[12pt]{article}
\usepackage[english]{babel}
\usepackage[utf8]{inputenc}
\usepackage[T1]{fontenc}
\usepackage{parskip}
\usepackage{comment}
\usepackage{cite}

\providecommand{\MR}[1]{}

\usepackage[title,titletoc]{appendix}
\usepackage{amsmath, amssymb, amsthm}
\usepackage{mathtools}
\usepackage{xcolor}
\definecolor{myblue}{RGB}{30,70,180}
\definecolor{mygreen}{RGB}{20,120,60}
\definecolor{myred}{RGB}{180,40,40}
\definecolor{mygray}{gray}{0.35}
\definecolor{myrev}{RGB}{170,20,160}

\usepackage[colorlinks=true,
            linkcolor=myblue,
            citecolor=mygreen,
            urlcolor=myred]{hyperref}
\usepackage[a4paper,margin=2.8cm]{geometry}
\newtheorem{theorem}{Theorem}[section]
\newtheorem{proposition}[theorem]{Proposition}
\newtheorem{lemma}[theorem]{Lemma}
\newtheorem{corollary}[theorem]{Corollary}

\theoremstyle{definition}
\newtheorem{definition}[theorem]{Definition}

\theoremstyle{remark}
\newtheorem{remark}[theorem]{Observation}

\newcommand{\R}{\mathbb{R}}
\newcommand{\N}{\mathbb{N}}

\newcommand{\C}{\mathbb{C}}

\newcommand{\norm}[1]{\left\lVert #1 \right\rVert}

\newcommand{\fourier}[1]{\widehat{#1}}
\newcommand{\embed}{\hookrightarrow}

\title{Orbital Stability and Instability of Slow and Fast Solitary Waves of gILW}
\author{
Diego Yupanqui\thanks{Instituto de Matemática Pura e Aplicada - IMPA, email: diego.huaman@impa.br}
}
\date{}
\begin{document}
\maketitle
\bigskip
\begin{abstract}
    We study the solitary waves of the generalized Intermediate Long Wave equation (gILW), for a power-type nonlinearity of order $p\geq1$. Using concentration-compactness and symmetric rearrangement, we prove existence of suitable (even, strictly positive and radially decreasing) solitary waves of gILW for every $p\geq1$, and show that, as the depth parameter tends to infinity, they converge in the energy space and in $L^q$, $q\in[2,\infty]$, to the solitary waves of the generalized Benjamin-Ono equation (gBO), while, as the depth parameter tends to zero, after a suitable rescaling, they converge in $L^q$, $q\in[2,\infty]$, to the solitary waves of the generalized Korteweg-de Vries equation (gKdV). Our main results concern the orbital stability of slow and fast solitary waves of gILW, that is, of those for which the product $\delta c$ of the depth parameter and the velocity is small or large, respectively. For slow waves we obtain stability for $p=1,2,3$ and instability for $p\geq5$ -- the same threshold as for gKdV -- while for fast waves we obtain stability for $p=1$ and instability for $p\geq3$ -- the same threshold as for gBO. In particular, for $p=3$ the stability of the solitary waves changes from one regime to the other. These results extend to $p>1$ the classical result of Albert, Bona and Henry for the Intermediate Long Wave equation ($p=1$). The proof proceeds by perturbation theory, transferring spectral information from gKdV and gBO to gILW through the convergence of the corresponding solitary-wave operators, and in particular avoids any use of an explicit formula for the solitary wave profile, none being available for $p>1$.
\end{abstract}

\section{Introduction}\label{Section: Introduction}
This paper studies the generalized Intermediate Long Wave (gILW) equation
\begin{equation}\label{Eq: GILW}
    \partial_tu + \frac{1}{\delta}\partial_xu + \mathcal{T}_\delta\partial_{x}^2u + u^pu_x = 0, \quad (x,t)\in\R\times\R.
\end{equation}
Here, $u=u(x,t)$ is a real-valued function, $1\leq p\in \N$, $\delta>0$ is a parameter, called the depth parameter, and
\begin{equation*}
    \mathcal{T}_\delta f(x):= -\frac{1}{2\delta}\text{p.v.}\int \coth\left( \frac{\pi(x-y)}{2\delta} \right)f(y)dy.
\end{equation*}
For $p=1$, equation \eqref{Eq: GILW} becomes the well-known Intermediate Long Wave equation (ILW), which models the propagation of long, weakly nonlinear internal gravity waves at the interface of a two-layer stratified fluid whose lower layer has finite depth, proportional to $\delta$ (see \cite{Bona2008, Chen1979, Craig2005, Joseph1977, Joseph1978, Kubota1978, Satsuma1979} for formal and rigorous derivations).

It is well-known that ILW interpolates between the two classical models of long, weakly nonlinear water waves: the Korteweg-de Vries equation (KdV)
\begin{equation*}
    \partial_tu + uu_x + \partial_x^3u = 0,
\end{equation*}
which arises formally as the shallow-water limit $\delta\to0$, and the Benjamin-Ono equation (BO)
\begin{equation*}
    \partial_tu + uu_x - \mathcal{H}\partial_x^2u = 0,
\end{equation*}
$\mathcal{H}$ denoting the Hilbert transform, which arises as the deep-water limit $\delta\to\infty$. This double limiting behavior was made rigorous by Abdelouhab, Bona, Felland and Saut \cite{Abdelouhab1989}, who proved that, for the same initial data, solutions of ILW converge to the corresponding solutions of BO as $\delta\to\infty$, and, after the rescaling $v(x,t):=3\delta^{-1}u(x,3t/\delta)$, to the corresponding solutions of KdV as $\delta\to0$. The need for this rescaling in the shallow-water regime, absent in the deep-water case, is an early instance of a distinction between the two limits that recurs throughout this paper. For more details on ILW, BO and their relation, we refer the reader to the survey article \cite{Saut2019}; while for recent results on the shallow-water limit of solutions, to \cite{Chapouto2025,Chapouto2026}. 

Among all solutions of gILW and of its two limiting models, this paper is concerned exclusively with solitary waves: solutions that travel without change of shape, at a constant speed, and that vanish at spatial infinity. Writing such a solution as $u(x,t)=\phi(x-ct)$ for a wave speed $c>0$ and a profile $\phi$, solitary waves are known to exist for each of the three classical models above, and, in each case, an explicit formula for $\phi$ is available. For KdV, the profile is given, up to translation, by 
\begin{equation*}
    \phi_c(x)=3c\,\text{sech}^2\left(\frac{\sqrt c}{2}x\right),
\end{equation*}
a classical fact recovered as a special case of the existence theory for gKdV recalled below \cite{Berestycki1983}. For BO, Benjamin \cite{Benjamin1967} found the explicit, algebraically decaying profile 
\begin{equation*}
    \phi_c(x)=\frac{4c}{1+c^2x^2}.
\end{equation*}
For ILW, Joseph \cite{Joseph1977} found the profile
\begin{equation*}
    \phi_{\delta,c}(x) = \frac{b}{\cosh^2(ax)+\frac{b^2}{16a^2}\sinh^2(ax)},
\end{equation*}
where $a,b>0$ are determined in terms of $c,\delta$. This profile decays exponentially, in contrast with the algebraic decay of the BO profile. It is straightforward to check, directly from these three explicit formulas, that the ILW solitary wave $\phi_{\delta,c}$ converges uniformly to the BO solitary wave as $\delta\to\infty$ and, after the same shallow-water rescaling used for the equation, to the KdV solitary wave as $\delta\to0$.

A basic question about any solitary wave is whether it is orbitally stable: informally, whether a solution starting close to the solitary wave stays close to some translate of it for all later times. Stability reflects a balance between dispersion, which tends to spread out a localized profile, and nonlinearity, which tends to concentrate it. To understand how far this balance can be pushed, that is, how strong the nonlinearity may be before stability is lost, one is led to the generalized models gKdV, gBO and gILW, in which the exponent $p\geq1$ is an arbitrary integer.

For gKdV, existence of a solitary wave is classical, and the profile is again explicit and unique up to translation, for every $p\geq1$ \cite{Berestycki1983}; it is orbitally stable for $1\leq p<4$, as proved by Albert, Bona and Henry \cite{Albert1987} and by Bona, Souganidis and Strauss \cite{Bona1987}, and unstable for $p>4$ \cite{Bona1987}, while in the critical case $p=4$ instability was proved by Martel and Merle \cite{Martel2001}. For gBO, existence of suitable solitary waves for every $p\geq1$ was proved by Angulo \cite{Angulo2009}, although an explicit formula is available only for $p=1$; uniqueness up to translation is likewise known, to the best of our knowledge, only for $p=1$, and was proved by Amick and Toland \cite{Amick1991}; another proof was later given by Albert \cite{Albert1995}. For every $p\geq1$, Frank and Lenzmann \cite{Frank2013} proved a weaker form of uniqueness, among variationally characterized solitary waves, together with their nondegeneracy: the kernel of the linearized operator is spanned by the derivative of the wave. This nondegeneracy is the spectral information on gBO on which our analysis of the deep water regime relies. The gBO solitary wave is orbitally stable for $1\leq p<2$ and unstable for $p\geq2$. Stability was proved by Weinstein \cite{Weinstein1987} for $p=1$, and for $1<p<2$ under the assumption that the kernel of the linearized operator is spanned by the derivative of the wave, which is precisely the nondegeneracy property later established by Frank and Lenzmann. Instability for $p>2$ was proved by Bona, Souganidis and Strauss \cite{Bona1987} and by Souganidis and Strauss \cite{Souganidis1990}, while, in the critical case $p=2$, which corresponds to the modified Benjamin-Ono equation, Martel and Pilod \cite{Martel2017} constructed a minimal mass blow-up solution, which implies the instability of the solitary wave. See also the monograph of Angulo \cite{Angulo2009}. For gILW, existence of a solitary wave is known for every $p\geq1$: this follows from variational methods and concentration-compactness, as established by Weinstein \cite{Weinstein1987}, Arnesen \cite{Arnesen2016} and Albert \cite{Albert1999}, although the latter approach, as originally carried out, only covers the range $1\leq p<2$. Moreover, Weinstein \cite{Weinstein1987} showed that, if the semigroup $e^{-M_\delta t}$, $M_\delta$ being the dispersive operator, is positivity preserving, then the solitary waves can be taken nonnegative, and that, if in addition its convolution kernel is radially decreasing (or a suitable uniqueness property holds), they can be taken even and radially decreasing. For gILW, the first property follows from the positivity of the convolution kernels of the operators $(\mu+M_\delta)^{-1}$, $\mu>0$, established by Albert, Bona and Henry \cite{Albert1987} by a residue calculation (see also Example 4 in \cite{Weinstein1987} and Lemma 1 in \cite{Albert1995}), so that the solitary waves of gILW obtained in \cite{Weinstein1987} are positive; however, to the best of our knowledge, the radial monotonicity of the kernel of $e^{-M_\delta t}$ has not been verified, and therefore the existence of even and radially decreasing solitary waves of gILW was not known for $p>1$. Prior to the present work, however, to the best of our knowledge, no orbital stability or instability result was available for gILW when $p>1$.

For ILW itself ($p=1$), by contrast, the question of orbital stability is completely settled, and its resolution proceeded in two steps. The first, due to Albert, Bona and Henry \cite{Albert1987}, established an abstract sufficient condition for orbital stability -- essentially, a spectral condition on the linearized operator about the solitary wave -- and verified it for ILW in the regimes where ILW is close to KdV or to BO, using only the uniform convergence of the ILW solitary wave to the KdV and BO ones, recalled above, together with the perturbation theory of closed operators, thereby avoiding any detailed description of the spectrum of the ILW solitary-wave operator itself. The second, and last, step was taken by Albert and Bona \cite{Albert1991}, who proved, using the explicit formula for the ILW solitary wave, that the required spectral conditions hold for every wave speed $c>0$ and depth parameter $\delta>0$. The key point of their proof is a positivity property: the Fourier transform of the ILW solitary wave is a Pólya frequency function of order 2 (a PF(2) function), which, through the theory of totally positive operators, yields the simplicity of the zero eigenvalue of the linearized operator, while the simplicity of its negative eigenvalue follows from the positivity of the solitary wave itself; the absence of further negative eigenvalues was then deduced from the perturbative result of \cite{Albert1987} by a continuity argument in $\delta$. Shortly afterwards, Albert \cite{Albert1992} showed that the PF(2) property alone is strong enough to control both the number and the simplicity of the non-positive eigenvalues of the linearized operator, which led to a simplified and unified proof of the stability of the KdV, BO and ILW solitary waves.

This second step cannot be replicated for the generalized equation gILW: no explicit formula is available for its solitary waves when $p>1$, so the positivity arguments of \cite{Albert1991, Albert1992} can not be directly applied. The present paper instead takes up the strategy of the first step: it generalizes the application of perturbation theory used by Albert, Bona and Henry, and proves orbital stability or instability, depending on $p$, of the solitary waves of gILW with $\delta c\ll1$ and with $\delta c\gg1$ -- for a fixed depth, the slow and the fast waves, respectively -- by perturbing spectral information already available for gKdV and for gBO, respectively.

The novelties of this paper can be summarized as follows.
\begin{itemize}
    \item Existence of suitable (even, strictly positive and radially decreasing) solitary waves of gILW is established for every $p\geq1$, $c>0$ and $\delta>0$, and for a broader family of nonlinearities which also allows non-integer values of $p$. The proof is direct, based on concentration-compactness and symmetric rearrangements, and requires no information on the semigroup $e^{-M_\delta t}$: the only property of $M_\delta$ that is used is that the convolution kernels of the operators $(\mu+M_\delta)^{-1}$, $\mu>0$, are positive and radially decreasing, which follows from an explicit series representation of these kernels due to Albert (Lemma 1 in \cite{Albert1995}).
    \item Convergence of the solitary waves of gILW to those of gBO as $\delta\to\infty$ is proved, in the energy space and in $L^q$ for $q\in[2,\infty]$; and, after the rescaling $(3/\delta)^{1/p}\phi_{\delta,\delta c/3}$, convergence to the solitary waves of gKdV as $\delta\to0$ is proved in $L^q$ for $q\in[2,\infty]$.
    \item Orbital stability and instability of the slow and fast solitary waves of gILW are established: if $\delta c\ll1$, the solitary waves are stable for $p=1,2,3$ and unstable for $p\geq5$, while if $\delta c\gg1$, they are stable for $p=1$ and unstable for $p\geq3$. For $p=1$ this recovers the result of \cite{Albert1987}; for $p>1$ the results are new. In particular, for $p=3$ stability is lost when passing from slow to fast waves. The cases $p=4$ for slow waves and $p=2$ for fast waves remain open. The approach is based on perturbation theory and on the convergence of solitary waves in the shallow and deep water regimes.
\end{itemize}

We close this introduction by fixing some notation used throughout the paper. We write $A\lesssim B$ if there is a constant $C>0$, independent of the relevant parameters, such that $A\leq CB$; if the constant is allowed to depend on some parameter $a$, we write $A\lesssim_a B$. For $1\leq q\leq\infty$, $|\cdot|_q$ denotes the $L^q(\R)$ norm, and $\norm{\cdot}_s$ the $H^s(\R)$ norm, $s\in\R$; the Fourier transform of $f$ is denoted $\fourier f$, with the normalization $\fourier f(\xi):=\int_\R f(x)e^{-2\pi ix\xi}dx$. As already seen above, the solitary waves considered in this paper are indexed by up to three parameters: the depth parameter $\delta$, the wave speed $c$, and, where relevant, a further parameter $\lambda$ arising from their variational characterization; whenever all three appear together as subscripts, they are always written in the order $\delta,c,\lambda$. Objects associated with the rescaled solitary waves of gILW used in the shallow water regime, introduced in Section \ref{Section: Convergence of solitary waves}, are distinguished by the superscript $(s)$.

The rest of the paper is organized as follows. Section \ref{Section: Preliminaries} collects background material: the general class of Long Wave models to which gILW, gKdV and gBO belong, the abstract criterion for orbital stability and instability, basic properties of gILW and of its scaling symmetry, the known results for gKdV and gBO, including those of Frank and Lenzmann, and the perturbation theory of closed operators. Section \ref{Section: A lower estimate for m delta s} isolates a technical lower estimate for the rescaled symbol $(3/\delta)m_\delta$, used repeatedly afterwards. Section \ref{Section: Existence of suitable solitary waves for GILW} proves existence of suitable solitary waves for gILW, for every $p\geq1$. Section \ref{Section: Convergence of solitary waves} proves their convergence to the solitary waves of gBO as $\delta\to\infty$ and, after rescaling, to those of gKdV as $\delta\to0$. Section \ref{Section: Orbital instability of solitary waves in shallow water} contains the main results of the paper: using perturbation theory to transfer spectral information from gKdV and gBO, we establish orbital stability and instability of the slow and fast solitary waves of gILW. Two appendices collect, for completeness, a proof of Helly's selection theorem and the basic facts about symmetric decreasing rearrangement used throughout.

\section{Preliminaries}\label{Section: Preliminaries}

This section collects the definitions and background results used throughout the paper. Section \ref{Section: General Long Wave type models} introduces the general class of Long Wave models to which gILW, gKdV and gBO belong, together with the associated solitary wave equation, solitary wave operator, and the notion of orbital stability. Section \ref{Section: Abstract stability and instability results} states the abstract criterion for orbital stability and instability that will be applied, in Section \ref{Section: Orbital instability of solitary waves in shallow water}, to the solitary waves of gILW. Section \ref{Section: The gILW equation and its solitary waves} collects some basic properties of the gILW equation, its symbol and its scaling symmetries. Section \ref{Section: Known results for gKdV and gBO} recalls the known results for the solitary waves of gKdV and gBO that will be used later, including the uniqueness and non-degeneracy results of Frank and Lenzmann for gBO. Finally, Section \ref{Section: Perturbation theory tools} recalls the perturbation theory of closed operators, the main tool used in Section \ref{Section: Orbital instability of solitary waves in shallow water} to transfer spectral information from gKdV and gBO to gILW.

\subsection{General Long Wave type models}\label{Section: General Long Wave type models}

Consider the general class of Long Wave models
\begin{equation}\label{Eq: General Long Wave model}
    u_t + u^pu_x - Mu_x = 0, \quad (x,t)\in\R\times\R,
\end{equation}
where $u=u(x,t)$ is a real-valued function, $p\in\N$, $p\geq1$, and $M$ is a dispersive operator defined as the Fourier multiplier
\begin{equation*}
    \fourier{Mf}(\xi) := m(\xi)\fourier{f}(\xi),
\end{equation*}
for a real-valued symbol $m$ satisfying
\begin{equation*}
    0 \leq m(\xi) \lesssim (1+|\xi|)^{\sigma}, \quad \forall\xi\in\R, \quad \text{and} \quad |\xi|^{\sigma} \lesssim m(\xi), \quad\forall\xi\gg1,
\end{equation*}
for some $\sigma\geq1$ (see Section 2.2 in \cite{Angulo2009}).

The gILW equation \eqref{Eq: GILW}, as well as the gKdV and gBO equations, are particular cases of the general Long Wave model \eqref{Eq: General Long Wave model}: gILW corresponds to the symbol $m_\delta$ introduced in Section \ref{Section: The gILW equation and its solitary waves}, gKdV corresponds to $m(\xi):=4\pi^2\xi^2$ (that is, $M=-\partial_x^2$), and gBO corresponds to $m(\xi):=2\pi|\xi|$ (that is, $M=\mathcal{H}\partial_x$, with $\mathcal{H}$ the Hilbert transform).

\begin{theorem}[\textbf{Local well-posedness}]\label{Theorem: Local well-posedness for Long Wave models}
    Let $s>3/2$. Then, the initial value problem associated with equation \eqref{Eq: General Long Wave model} is locally well-posed in $H^s$.
\end{theorem}
\begin{proof}
    See Theorem 4.2 in \cite{Angulo2009} and the references therein.
\end{proof}

A solution $u(x,t)$ of \eqref{Eq: General Long Wave model} is called a solitary wave if it takes the form $u(x,t)=\phi(x-ct)$ for some wave speed $c$, where the wave profile $\phi$ satisfies $\phi(x)\to0$ as $|x|\to\infty$. Replacing this ansatz into the equation and integrating, we obtain the solitary wave equation
\begin{equation}\label{Eq: General solitary wave equation}
    M\phi_c + c\phi_c - \frac{1}{p+1}\phi_c^{p+1} = 0.
\end{equation}
Note that if $\phi_c\in H^{\sigma/2}$ then a bootstrap argument shows that $\phi_c\in H^{\infty}$. In particular, $\phi_c$ satisfies the boundary condition $\phi_c(x)\to0$ as $|x|\to\infty$. For this reason, we define the solitary waves of \eqref{Eq: General Long Wave model} as solutions of \eqref{Eq: General solitary wave equation} in the energy space.

Given a solitary wave $\phi_c$ solution of \eqref{Eq: General solitary wave equation}, its associated solitary wave operator is the linear operator
\begin{equation*}
    \mathcal{L} := M + c - \phi_c^p,
\end{equation*}
obtained by linearizing \eqref{Eq: General solitary wave equation} around $\phi_c$.

\begin{theorem}\label{Theorem: Basic information about spectrum}
    The operator $\mathcal{L}$ is a closed, unbounded, self-adjoint operator on $L^2$ whose spectrum consists of the essential part $[c,\infty)$ together with a finite number of discrete eigenvalues, with finite-dimensional eigenspaces, in the interval $(-\infty,c)$. In particular, $0$ is an eigenvalue of $\mathcal{L}$, with associated eigenfunction $\phi_c'$.
\end{theorem}
\begin{proof}
    See Theorem B.43 in \cite{Angulo2009}. The fact that $0$ is an eigenvalue of $\mathcal{L}$ with eigenfunction $\phi_c'$ follows by differentiating \eqref{Eq: General solitary wave equation} with respect to $x$.
\end{proof}

Associated with equation \eqref{Eq: General Long Wave model}, for any $c>0$ define the functional
\begin{equation*}
    J_c(\psi) := \frac{1}{2}\int_\R \psi M\psi + \frac{c}{2}\int_\R \psi^2, \quad \psi\in H^{\sigma/2},
\end{equation*}
which will be used later for the variational characterization of suitable solitary waves (see Section \ref{Section: Existence of suitable solitary waves for GILW}).

\begin{definition}[\textbf{Suitable solitary wave}]\label{Def: Suitable solitary wave}
    A solitary wave $\phi$ of \eqref{Eq: General Long Wave model} is called suitable if it is even, strictly positive, and radially decreasing.
\end{definition}

We can now state the notion of orbital stability studied throughout this paper.

\begin{definition}[\textbf{Orbital stability}]\label{Def: Orbital stability}
    Let $\phi_c$ be a solitary wave of \eqref{Eq: General Long Wave model}. Following \cite{Albert1987}, we say that $\phi_c$ is orbitally stable if:
    \begin{enumerate}
        \item[(i)] there exists $s_0>3/2$ and a neighborhood $U$ of $\phi_c$ in $H^{s_0}$ such that, for every $\varphi\in U$, there is a unique solution $u\in C([0,\infty);H^{s_0})$ of \eqref{Eq: General Long Wave model} with $u(x,0)=\varphi(x)$; and
        \item[(ii)] for every $\epsilon>0$, there exists $\eta>0$ such that, for every $\varphi\in U$ satisfying
        \begin{equation*}
            \inf_{\gamma\in\R}\|\varphi-\phi_c(\cdot-\gamma)\|_{\sigma/2}<\eta,
        \end{equation*}
        the corresponding solution $u$ given by (i) satisfies
        \begin{equation*}
            \inf_{\gamma\in\R}\|u(\cdot,t)-\phi_c(\cdot-\gamma)\|_{\sigma/2}<\epsilon, \quad \text{for all } t>0.
        \end{equation*}
    \end{enumerate}
    Condition (i) is included because, unlike the local theory given by Theorem \ref{Theorem: Local well-posedness for Long Wave models}, no global well-posedness result in the energy space $H^{\sigma/2}$ is available, in general, for the Long Wave model \eqref{Eq: General Long Wave model}; not even for gILW is such a result known. Condition (i) circumvents this gap by requiring, as part of the very definition of orbital stability, the existence of unique global-in-time solutions for initial data in some $H^{s_0}$-neighborhood $U$ of $\phi_c$; orbital stability, in the sense of (ii), is then only asserted for that restricted class of initial data.
\end{definition}

\subsection{Abstract stability and instability results}\label{Section: Abstract stability and instability results}

We now recall a result that gives sufficient conditions for the orbital stability/instability of solitary waves of the general Long Wave model \eqref{Eq: General Long Wave model}. This result is taken in the form stated by Albert \cite{Albert1992} and combines the contributions of several authors: its stability part goes back to the sufficient conditions of Bona, Souganidis and Strauss \cite{Bona1987} and of Weinstein \cite{Weinstein1987}, and relies on the total positivity theory of Albert and Bona \cite{Albert1991}, while its instability part is due to Souganidis and Strauss \cite{Souganidis1990}.

Let $l(\xi):=\xi m(\xi)$. Consider the following conditions on the symbol $m$:
\begin{itemize}
    \item[(C1)] $m\in C^1(\R)\cap C^3(\R\backslash\{0\})$,
    \item[(C2)] $l''(\xi)$ is bounded away from $0$ as $|\xi|\to\infty$ and $l''$ has only a finite number of zeros, all of which are non-degenerate,
    \item[(C3)] there are some $\kappa>0$, $r\geq0$, such that $|\xi|^{r}\lesssim |l''(\xi)|$ for any $|\xi|<\kappa$.
\end{itemize}

\begin{theorem}[\textbf{Stability and instability criterion}]\label{Theorem: Stability instability criterion Albert}
    Let $\phi$ be a solitary wave solution of \eqref{Eq: General Long Wave model} satisfying condition (i) in Definition \ref{Def: Orbital stability}, and suppose that its solitary wave operator $\mathcal{L}$ has the following properties:
    \begin{enumerate}
        \item $\mathcal{L}$ has a simple, negative eigenvalue $\mu$,
        \item $\mathcal{L}$ has no negative eigenvalue other than $\mu$,
        \item the eigenvalue $0$ of $\mathcal{L}$ is simple.
    \end{enumerate}
    Choose $\psi\in L^2$ such that $\mathcal{L}\psi = \phi$, and define $I:= \langle\psi,\phi\rangle$. If $I<0$, then $\phi$ is orbitally stable. If $I>0$ and the symbol $m$ satisfies conditions (C1), (C2) and (C3), then $\phi$ is unstable.
\end{theorem}
\begin{proof}
    The stability part is Theorem 3.1 in \cite{Albert1992}; for the details of its proof see Theorem 2 in \cite{Albert1991}. The instability part, under conditions (C1), (C2) and (C3) on the symbol, is due to Souganidis and Strauss; see Section 4 in \cite{Souganidis1990}, \cite{Bona1987}, and Remark (ii) after Theorem 3.1 in \cite{Albert1992}. Note that, by condition 3 and since $\langle\phi,\phi'\rangle=0$, the equation $\mathcal{L}\psi=\phi$ has $L^2$-solutions, and $I$ does not depend on the choice of $\psi$.
\end{proof}

\subsection{The gILW equation and its solitary waves}\label{Section: The gILW equation and its solitary waves}

Consider the gILW equation \eqref{Eq: GILW}. It is a particular case of \eqref{Eq: General Long Wave model} with $M=M_\delta$, the Fourier multiplier with symbol
\begin{equation*}
    m_\delta(\xi):=2\pi\xi\coth(2\pi\delta\xi)-\frac{1}{\delta}.
\end{equation*}
Following Section \ref{Section: General Long Wave type models}, we denote by $\phi_{\delta,c}$ a solitary wave of gILW with velocity $c$, making explicit the dependence on the depth parameter $\delta$. It satisfies the solitary wave equation
\begin{equation}\label{Eq: GILW solitary wave equation}
    M_\delta\phi_{\delta,c} + c\phi_{\delta,c} - \frac{1}{p+1}\phi_{\delta,c}^{p+1} = 0,
\end{equation}
and its solitary wave operator is $\mathcal{L}_{\delta,c}:=M_\delta+c-\phi_{\delta,c}^p$.

\begin{lemma}\label{Lemma: m delta and |xi|}
    Fix an arbitrary $\delta>0$. Then, for any $\xi\in\R$,
    \begin{equation*}
       2\pi|\xi|-\frac{2}{\delta} \leq m_\delta(\xi)\leq 2\pi|\xi|.
    \end{equation*}
    In particular, for any $c>0$,
    \begin{equation*}
        1+|\xi| \leq \max\left\{\frac{1}{2\pi},\frac{1}{c}+\frac{1}{\delta\pi c}\right\} (c+m_\delta(\xi))
    \end{equation*}
    and
    \begin{equation*}
        c+m_\delta(\xi) \leq \max\{c,2\pi\}(1+|\xi|).
    \end{equation*}
\end{lemma}
\begin{proof}
    See Lemma 4.1 in \cite{Abdelouhab1989}.
\end{proof}

Note also that $m_\delta\geq0$, since $y\coth y\geq1$ for every $y\in\R$. Following Section \ref{Section: General Long Wave type models}, denote by $J_{\delta,c}$ the functional $J_c$ particularized to $M=M_\delta$, that is,
\begin{equation*}
    J_{\delta,c}(\psi) := \frac{1}{2}\int_\R \psi M_\delta\psi + \frac{c}{2}\int_\R\psi^2.
\end{equation*}

\begin{corollary}\label{Corollary: J delta and norm 1/2}
    For any $\delta>0$,
    \begin{equation*}
        J_{\delta,c}(\psi)\lesssim_c \norm{\psi}_{1/2}^2 \lesssim_{\delta,c} J_{\delta,c}(\psi)
    \end{equation*}
    and the upper estimate is uniform in $\delta$ (the constant is independent of $\delta$) as $\delta\to\infty$. In particular, for any fixed values $\delta,c>0$, the functional $\sqrt{J_{\delta,c}(\cdot)}$ defines a norm in $H^{1/2}$ equivalent to the standard one.
\end{corollary}

In the deep water limit $\delta\to\infty$, the symbol $m_\delta$ converges to $2\pi|\xi|$, the symbol of gBO. In the shallow water limit $\delta\to0$, instead, gILW converges to gKdV only after the rescaling $v(x,t):=(3/\delta)^{1/p}u(x,3t/\delta)$, which multiplies the symbol by $3/\delta$; indeed, $(3/\delta)m_\delta(\xi)\to4\pi^2\xi^2$ as $\delta\to0$ \cite{Abdelouhab1989}. In this paper, this rescaling is only used at the level of the solitary waves, as a technical device; see Section \ref{Section: Convergence of solitary waves}.

We finally describe the scaling symmetry of gILW. For $\theta>0$ we define
\begin{equation*}
    D_\theta\varphi(x):=\theta^{1/p}\varphi(\theta x) \qquad\text{ and }\qquad \tau_\theta f(x):=f(\theta x).
\end{equation*}
A direct computation from the definition of $m_\delta$ shows that $m_{\delta/\theta}(\xi)=\theta\, m_\delta(\xi/\theta)$ for every $\xi\in\R$.

\begin{lemma}[\textbf{Scaling symmetries of gILW}]\label{Lemma: Scaling symmetries}
    Let $\delta,c,\theta>0$. Then:
    \begin{enumerate}
        \item For every $s\in\R$, $D_\theta$ is a bounded linear isomorphism of $H^s$ onto itself, with inverse $D_{1/\theta}$. Moreover, $D_\theta$ preserves evenness, strict positivity and radial monotonicity, and $D_\theta[f(\cdot-\gamma)]=(D_\theta f)(\cdot-\gamma/\theta)$ for every $\gamma\in\R$.
        \item If $u\in C([0,\infty);H^s)$, $s>3/2$, is a solution of gILW with depth parameter $\delta$, then
        \begin{equation*}
            \tilde u(x,t):=\theta^{1/p}u(\theta x,\theta^{2}t), \quad\text{ that is, }\quad \tilde u(\cdot,t)=D_\theta\left[u(\cdot,\theta^{2}t)\right],
        \end{equation*}
        defines a solution $\tilde u\in C([0,\infty);H^s)$ of gILW with depth parameter $\delta/\theta$, with $\tilde u(\cdot,0)=D_\theta[u(\cdot,0)]$.
        \item A function $\phi$ is a solitary wave of gILW with depth parameter $\delta$ and velocity $c$ if and only if $D_\theta\phi$ is a solitary wave of gILW with depth parameter $\delta/\theta$ and velocity $\theta c$. Moreover, $\phi$ is suitable if and only if $D_\theta\phi$ is suitable.
        \item If $\phi_{\delta,c}$ is a solitary wave of gILW with depth parameter $\delta$ and velocity $c$, and $\phi_{\delta/\theta,\theta c}:=D_\theta\phi_{\delta,c}$ is the solitary wave given by item 3, then their solitary wave operators satisfy
        \begin{equation}\label{Eq: scaling of solitary wave operators}
            \mathcal{L}_{\delta/\theta,\theta c}=\theta\left[\tau_\theta\mathcal{L}_{\delta,c}\tau_\theta^{-1}\right].
        \end{equation}
        Consequently, $\text{spec}(\mathcal{L}_{\delta/\theta,\theta c})=\{\theta\mu:\ \mu\in\text{spec}(\mathcal{L}_{\delta,c})\}$, and $\varphi$ is an eigenfunction of $\mathcal{L}_{\delta/\theta,\theta c}$ with eigenvalue $\theta\mu$ if and only if $\tau_\theta^{-1}\varphi$ is an eigenfunction of $\mathcal{L}_{\delta,c}$ with eigenvalue $\mu$, corresponding eigenvalues having the same multiplicity.
        \item The product $\delta c$ is invariant under the transformation $(\delta,c)\mapsto(\delta/\theta,\theta c)$. Conversely, if $\delta',c'>0$ satisfy $\delta'c'=\delta c$, then $(\delta',c')=(\delta/\theta,\theta c)$ with $\theta:=\delta/\delta'$.
        \item If a solitary wave $\phi$ of gILW with depth parameter $\delta$ satisfies condition (i) in Definition \ref{Def: Orbital stability}, with associated regularity $s_0$ and neighborhood $U\subset H^{s_0}$, then the solitary wave $D_\theta\phi$ of gILW with depth parameter $\delta/\theta$ also satisfies condition (i), with the same $s_0$ and the neighborhood $D_\theta(U)$.
    \end{enumerate}
\end{lemma}
\begin{proof}
    The proof follows from straightforward computations.
\end{proof}

\subsection{Known results for gKdV and gBO}\label{Section: Known results for gKdV and gBO}

We collect below the known existence, uniqueness, and stability results for the suitable solitary waves of gKdV and gBO.

For gKdV, existence and uniqueness hold for every admissible $p\geq1$, as a consequence of the one-dimensional theory of ordinary differential equations.

\begin{theorem}[\textbf{Existence and uniqueness for gKdV}]\label{Theorem: Existence and uniqueness for gKdV}
    For any $c>0$, any admissible order of nonlinearity $p\geq1$ and any of the nonlinear terms $g\in\{g_I,g_N,g_L\}$ of Section \ref{Section: Existence of suitable solitary waves for GILW}, there exists a suitable solitary wave $\phi_{0,c}$ of gKdV, and it is unique up to translations. Moreover, $\phi_{0,c}$ is given explicitly by
    \begin{equation*}
        \phi_{0,c}(x) = \left[\frac{c(p+1)(p+2)}{2}\,\text{sech}^2\left(\frac{p\sqrt{c}}{2}x\right)\right]^{1/p}.
    \end{equation*}
\end{theorem}
\begin{proof}
    This is Theorem 5 in \cite{Berestycki1983}, applied to the equation $-\phi''+c\phi=\frac{1}{p+1}\phi^{p+1}$; the explicit formula follows from direct substitution.
\end{proof}

It is well known that the solitary waves of gKdV are orbitally stable for $1\leq p<4$, as proved by Albert, Bona and Henry \cite{Albert1987} (see also \cite{Bona1987}), and unstable for $p>4$ \cite{Bona1987} and for $p=4$ \cite{Martel2001}.

For gBO we have the following results.

\begin{theorem}[\textbf{Existence for gBO}]\label{Theorem: Existence for gBO}
    For any $c>0$, any admissible order of nonlinearity $p\geq1$ and any of the nonlinear terms $g\in\{g_I,g_N,g_L\}$ of Section \ref{Section: Existence of suitable solitary waves for GILW}, there exists a suitable solitary wave $\phi_{\infty,c}$ of gBO.
\end{theorem}
\begin{proof}
    For $g=g_I$ this is Theorem 8.3 in \cite{Angulo2009}. The proof follows the same variational scheme, based on concentration-compactness and rearrangement inequalities, that we use in Section \ref{Section: Existence of suitable solitary waves for GILW} to prove existence of suitable solitary waves for gILW, and that scheme applies verbatim, with $M_\delta$ replaced by $\mathcal{H}\partial_x$, to the nonlinearities $g_N$ and $g_L$.
\end{proof}

For $p=1$, uniqueness (up to translation) of the positive solitary wave of gBO, among all positive solutions in $H^{1/2}$, was proved by Amick and Toland \cite{Amick1991}; another proof was later given by Albert \cite{Albert1995}. For every admissible $p\geq1$, uniqueness of the suitable solitary wave of gBO obtained variationally, as in Theorem \ref{Theorem: Existence for gBO}, follows from the work of Frank and Lenzmann \cite{Frank2013}; see Proposition \ref{Prop: gBO waves are FL ground states} at the end of this subsection.

We now give a complete description of the spectrum of $\mathcal{L}_{0,c}$ near $0$, valid for every $p\geq1$; this is what will allow us, in Section \ref{Section: Orbital instability of solitary waves in shallow water}, to transfer spectral information from gKdV to gILW for every $p$.

\begin{theorem}[\textbf{Spectrum of the gKdV solitary wave operator}]\label{Theorem: Spectrum of gKdV solitary wave operator}
    Let $c>0$ and $p\geq1$, $p\in\N$, and let $\phi_{0,c}$ be the suitable solitary wave of gKdV given by Theorem \ref{Theorem: Existence and uniqueness for gKdV}. Then the associated solitary wave operator $\mathcal{L}_{0,c}=-\partial_x^2+c-\phi_{0,c}^p$ has exactly one negative eigenvalue $\mu_0$, and it is simple; the eigenvalue $0$ is simple, with associated eigenfunction $\partial_x\phi_{0,c}$; and there exists $\beta_0>0$ such that every $\lambda\in\text{spec}(\mathcal{L}_{0,c})\backslash\{\mu_0,0\}$ satisfies $\lambda>\beta_0$. In particular, $\mathcal{L}_{0,c}$ satisfies conditions 1, 2 and 3 in Theorem \ref{Theorem: Stability instability criterion Albert}, and $0$ is an isolated point of $\text{spec}(\mathcal{L}_{0,c})$, so that
    \begin{equation*}
        \gamma_0 := \text{dist}\left(0,\text{spec}(\mathcal{L}_{0,c})\backslash\{0\}\right) \geq \min\{|\mu_0|,\beta_0\} > 0.
    \end{equation*}
\end{theorem}
\begin{proof}
    This is Theorem B.61 in \cite{Angulo2009}, applied to the smooth nonlinearity $f(t)=t^{p+1}/(p+1)$, which satisfies $f(0)=f'(0)=0$ and $f'(\phi_{0,c})=\phi_{0,c}^p$, so that equation (B.36) there becomes the gKdV solitary wave equation and the operator considered there becomes $\mathcal{L}_{0,c}$. The only hypothesis left to verify is that $\partial_x\phi_{0,c}$ has exactly one zero, which follows from the explicit formula in Theorem \ref{Theorem: Existence and uniqueness for gKdV}: $\phi_{0,c}$ is even, strictly positive and strictly decreasing on $(0,\infty)$, so that $\partial_x\phi_{0,c}>0$ on $(-\infty,0)$, $\partial_x\phi_{0,c}(0)=0$ and $\partial_x\phi_{0,c}<0$ on $(0,\infty)$.
\end{proof}

For gBO, a description of the spectrum of $\mathcal{L}_{\infty,c}$ near $0$ analogous to Theorem \ref{Theorem: Spectrum of gKdV solitary wave operator}, valid for every admissible $p\geq1$, follows from the work of Frank and Lenzmann \cite{Frank2013} on ground states of fractional Laplacians. We state their results in the particular case $s=1/2$, $\alpha=p$, which is the one relevant for gBO; note that $(-\Delta)^{1/2}=\mathcal{H}\partial_x$, the Fourier multiplier with symbol $2\pi|\xi|$, and that, in this case, the critical exponent $\alpha_{\max}(1/2)$ in \cite{Frank2013} is infinite, so that every real $p>0$ is allowed. Consider the equation
\begin{equation}\label{Eq: FL ground state equation}
    \mathcal{H}\partial_xQ+Q-Q^{p+1}=0, \quad\text{ in }\R,
\end{equation}
and the functional
\begin{equation*}
    \mathcal{J}(u):=\left(\int_\R|u|^{p+2}\right)^{-1}\left(\int_\R u\,\mathcal{H}\partial_xu\right)^{p/2}\int_\R|u|^2, \quad u\in H^{1/2}\backslash\{0\}.
\end{equation*}
An even and positive solution $Q\in H^{1/2}$ of \eqref{Eq: FL ground state equation} which minimizes $\mathcal{J}$ is called a ground state.

\begin{theorem}[\textbf{Frank-Lenzmann}]\label{Theorem: Frank-Lenzmann}
    Let $p\geq1$ be real. Then:
    \begin{enumerate}
        \item There exists a ground state $Q$ of \eqref{Eq: FL ground state equation}, and it is unique. Moreover, $Q$ is strictly decreasing in $|x|$, and every minimizer of $\mathcal{J}$ is of the form $\beta Q(\gamma(\cdot+y))$ for some $\beta\neq0$, $\gamma>0$ and $y\in\R$.
        \item $Q\in H^2$, and the function $R:=\frac{1}{p}Q+xQ'$ also belongs to $H^2$.
        \item The linearized operator $\mathcal{L}_+:=\mathcal{H}\partial_x+1-(p+1)Q^p$, acting on $L^2$, is non-degenerate: $\ker\mathcal{L}_+=\text{span}\{Q'\}$.
    \end{enumerate}
\end{theorem}
\begin{proof}
    Item 1 follows from Proposition 1.1 (i) and Theorem 2.4 in \cite{Frank2013}; item 2 follows from Proposition 1.1 (iii) in \cite{Frank2013} and from the remark following (4.5) there, where it is shown that $R\in H^{2s+1}$, with $s=1/2$; and item 3 is Theorem 2.3 in \cite{Frank2013}.
\end{proof}

\begin{proposition}[\textbf{Spectrum of the gBO solitary wave operator}]\label{Prop: gBO waves are FL ground states}
    Let $c>0$ and an admissible $p\geq1$, and let $\phi_{\infty,c}=\eta\,\psi$ be a suitable solitary wave of gBO obtained, as in Theorem \ref{Theorem: Existence for gBO}, from a minimizer $\psi$ of the variational problem
    \begin{equation*}
        I_{\infty,c,\lambda} := \inf\left\{J_{\infty,c}(u):\ u\in H^{1/2},\ \int_\R u^{p+2}=\lambda\right\},
    \end{equation*}
    for some $\lambda>0$, where $J_{\infty,c}$ denotes the functional $J_c$ of Section \ref{Section: General Long Wave type models} with $M=\mathcal{H}\partial_x$, and $\eta>0$ is the constant determined by the corresponding Lagrange multiplier (see Section \ref{Section: Existence of suitable solitary waves for GILW}). Then:
    \begin{enumerate}
        \item $\phi_{\infty,c}=A\,\tau_cQ$, where $A:=((p+1)c)^{1/p}$ and $Q$ is the ground state of Theorem \ref{Theorem: Frank-Lenzmann}. In particular, $\phi_{\infty,c}$ depends neither on $\lambda$ nor on the chosen minimizer $\psi$, so that the suitable solitary wave of gBO given by Theorem \ref{Theorem: Existence for gBO} is unique.
        \item The solitary wave operator $\mathcal{L}_{\infty,c}=\mathcal{H}\partial_x+c-\phi_{\infty,c}^p$ satisfies $\mathcal{L}_{\infty,c}=c\,\tau_c\mathcal{L}_+\tau_c^{-1}$. It has exactly one negative eigenvalue $\mu_\infty$, and it is simple; the eigenvalue $0$ is simple, with associated eigenfunction $\partial_x\phi_{\infty,c}$; and there exists $\beta_\infty>0$ such that every $\lambda\in\text{spec}(\mathcal{L}_{\infty,c})\backslash\{\mu_\infty,0\}$ satisfies $\lambda>\beta_\infty$. In particular, $\mathcal{L}_{\infty,c}$ satisfies conditions 1, 2 and 3 in Theorem \ref{Theorem: Stability instability criterion Albert}, and $0$ is an isolated point of $\text{spec}(\mathcal{L}_{\infty,c})$, so that
        \begin{equation*}
            \gamma_\infty := \text{dist}\left(0,\text{spec}(\mathcal{L}_{\infty,c})\backslash\{0\}\right) \geq \min\{|\mu_\infty|,\beta_\infty\} > 0.
        \end{equation*}
    \end{enumerate}
\end{proposition}
\begin{proof}
    \emph{Item 1.} By homogeneity, $\psi$ minimizes the quotient $W(u):=2J_{\infty,c}(u)/|u|_{p+2}^2$ over $H^{1/2}\backslash\{0\}$. Since the term $\int u\,\mathcal{H}\partial_xu$ is invariant under the dilations $u\mapsto u(\mu\,\cdot)$, $\mu>0$, minimizing $W(u(\mu\,\cdot))$ over $\mu>0$ shows that $\inf_{\mu>0}W(u(\mu\,\cdot))=C_{p,c}\,\mathcal{J}(u)^{2/(p+2)}$ for every $u\neq0$, with a constant $C_{p,c}>0$ independent of $u$; hence $\psi$ is also a minimizer of $\mathcal{J}$. By item 1 of Theorem \ref{Theorem: Frank-Lenzmann}, and since $\psi$ is suitable, $\psi=\beta Q(\gamma\,\cdot)$ for some $\beta,\gamma>0$, so that $\phi_{\infty,c}=\eta\beta\,Q(\gamma\,\cdot)$. Inserting this expression into the gBO solitary wave equation $\mathcal{H}\partial_x\phi+c\phi-\frac{1}{p+1}\phi^{p+1}=0$ and using \eqref{Eq: FL ground state equation}, we obtain $\gamma=c$ and $(\eta\beta)^p=(p+1)c$, that is, $\phi_{\infty,c}=A\,\tau_cQ$.

    \emph{Item 2.} The identity $\mathcal{L}_{\infty,c}=c\,\tau_c\mathcal{L}_+\tau_c^{-1}$ follows from item 1 and $\mathcal{H}\partial_x\tau_c=c\,\tau_c\mathcal{H}\partial_x$. Hence, by item 3 of Theorem \ref{Theorem: Frank-Lenzmann}, $\ker\mathcal{L}_{\infty,c}=\text{span}\{\partial_x\phi_{\infty,c}\}$, that is, the eigenvalue $0$ is simple. Next, since $\psi$ minimizes $J_{\infty,c}$ under the constraint $\int u^{p+2}=\lambda$, the second variation gives $\langle\mathcal{L}_{\infty,c}v,v\rangle\geq0$ for every $v$ in the domain of $\mathcal{L}_{\infty,c}$ with $\langle v,\psi^{p+1}\rangle=0$, a subspace of codimension one; by the min-max principle, $\mathcal{L}_{\infty,c}$ has at most one negative eigenvalue, counted with multiplicity. On the other hand, $\mathcal{L}_{\infty,c}\phi_{\infty,c}=-\frac{p}{p+1}\phi_{\infty,c}^{p+1}$, so that $\langle\mathcal{L}_{\infty,c}\phi_{\infty,c},\phi_{\infty,c}\rangle<0$, and therefore $\mathcal{L}_{\infty,c}$ has exactly one negative eigenvalue $\mu_\infty$, which is simple. Finally, by Theorem \ref{Theorem: Basic information about spectrum}, the spectrum of $\mathcal{L}_{\infty,c}$ in $(-\infty,c)$ consists of finitely many eigenvalues, which gives $\beta_\infty>0$.
\end{proof}

\begin{remark}\label{Remark: consequences of FL for the deep water limit}
    Uniqueness of the positive solitary waves of gBO among all positive solutions of the solitary wave equation, which is the statement proved in \cite{Albert1995,Amick1991} for $p=1$, is not covered by \cite{Frank2013} when $p>1$.
\end{remark}

\subsection{Perturbation theory tools}\label{Section: Perturbation theory tools}

The following classical tools from the perturbation theory of closed operators will be used later to study how the spectral properties of a solitary wave operator behave under perturbations. Proofs of the results in this section can be found in \cite{Kato1995}.

Let $L^2\times L^2$ be the Hilbert space with norm
\begin{equation*}
    \|(f,g)\|_{L^2\times L^2}=(\|f\|^2+\|g\|^2)^{1/2}.
\end{equation*}
For any closed operator $T:L^2\to L^2$, with domain $D(T)$, define its graph as
\begin{equation*}
    G(T) := \{(f,g)\in L^2\times L^2: f\in D(T), g=T(f)\}.
\end{equation*}
Let $Cl(L^2)$ be the space of closed operators in $L^2$ with metric $d$ defined as
\begin{equation*}
    d(S,T) := \|P_S-P_T\|_{B(L^2\times L^2)}
\end{equation*}
where $P_S$ and $P_T$ are the orthogonal projections onto the sets $G(S)$ and $G(T)$, respectively, and $\|\cdot\|_{B(L^2\times L^2)}$ denotes the operator norm in the space of bounded operators on $L^2\times L^2$.

\begin{proposition}\label{Prop: Basic properties of gap metric}
    Let $S,T\in Cl(L^2)$ and suppose $A$ is a bounded operator in $L^2$ with operator norm $\|A\|_{2,2}$. Then:
    \begin{itemize}
        \item $d(T+A,T)\leq \|A\|_{2,2}$,
        \item $d(S+A,T+A)\leq 2(1+\|A\|_{2,2}^2)d(S,T)$.
    \end{itemize}
    If, in addition, $S$ and $T$ are invertible, then
    \begin{itemize}
        \item $d(S,T)=d(S^{-1},T^{-1})$.
    \end{itemize}
\end{proposition}

\begin{proposition}\label{Prop: Perturbation theorem}
    Let $T\in Cl(L_2)$ and let $\Omega$ denote an open set of the complex plane whose boundary is a smooth contour $\Gamma$. Suppose that $\text{spec}(T)\cap \Gamma = \emptyset$ and $\text{spec}(T)\cap \Omega$ consists of a finite number of eigenvalues of $T$, each with finite (algebraic) multiplicity. Then, there exists some $\epsilon>0$ such that if $S\in Cl(L_2)$ and $d(S,T)<\epsilon$, then $\text{spec}(S)\cap \overline{\Omega}$ consists of a finite number of eigenvalues of finite multiplicity, the sum of the multiplicities been equal to the sum of the multiplicities of the eigenvalues of $T$ in $\Omega$.

    In particular, suppose $\text{spec}(T)\cap\Omega$ consists of a single, simple eigenvalue $\mu$. If $\{S_n\}$ is a sequence in $Cl(L_2)$ such that $d(S_n,T)\to0$ as $n\to\infty$, then, for $n$ large, $\text{spec}(S_n)\cap\Omega$ consists of a single simple eigenvalue $\mu_n$, and $\mu_n\to\mu$ as $n\to\infty$.
\end{proposition}

\section{A lower estimate for the symbol of gILW}\label{Section: A lower estimate for m delta s}

This short section isolates a single technical estimate, used repeatedly from Section \ref{Section: Convergence of solitary waves} onward: a lower bound for the rescaled symbol $(3/\delta)m_\delta$, where $m_\delta$ is the symbol of gILW, valid for every $\delta>0$ with constants independent of $\delta$. The factor $3/\delta$ is the natural normalization in the shallow water regime, since $(3/\delta)m_\delta(\xi)\to4\pi^2\xi^2$, the symbol of gKdV, as $\delta\to0$. Together with the upper bound of Lemma \ref{Lemma: m delta and |xi|}, this estimate is the main quantitative input behind the uniform-in-$\delta$ estimates for the solitary waves in the shallow water regime established later in the paper.

\begin{lemma}\label{Lemma: Lower estimate for m_s_delta}
    Let $\alpha_1<\alpha_2$ be constants defined as
    \begin{equation*}
        \alpha_1 := \frac{1}{\pi^2}\left(\frac{90}{96}\right)^{1/2} \quad \text{and} \quad \alpha_2:=\frac{3}{\pi}.
    \end{equation*}
    Then, for any $\delta>0$ and $\xi\in\R$:
    \begin{equation*}
        \frac{3}{\delta}m_{\delta}(\xi) \geq \left\{\begin{array}{cc}
            5\pi|\xi|^2, & |\xi|\leq \alpha_1/\delta, \\
            5\pi\alpha_1^2/\delta^2, & \alpha_1/\delta \leq |\xi| \leq \alpha_2/\delta, \\
            5\pi|\xi|/\delta, & \alpha_2/\delta\leq |\xi|.
        \end{array}\right.
    \end{equation*}
\end{lemma}
\begin{proof}
    By evenness of $m_\delta$ (since $\xi\coth(2\pi\delta\xi)$ is an even function of $\xi$), it suffices to consider $\xi>0$; the bound for $\xi<0$ then follows with $|\xi|$ in place of $\xi$. By Lemma 8.2.1 in \cite{Abdelouhab1989},
    \begin{equation*}
        \frac{3}{\delta}m_{\delta}(\xi) = 4\pi^2\xi^2 - \xi^2h(\xi,\delta), \quad\text{ where }\quad h(\xi,\delta) = \sum_{n=1}^{\infty}\frac{96\delta^2\xi^2}{n^2(n^2+4\delta^2\xi^2)},
    \end{equation*}
    and we treat separately the three ranges of $\xi$ appearing in the statement.

    \emph{Step 1: small frequencies, $0<\xi\leq\alpha_1/\delta$.} Dropping the term $4\delta^2\xi^2$ from the denominator,
    \begin{equation*}
        h(\xi,\delta) \leq \sum_{n=1}^{\infty}\frac{96\delta^2\xi^2}{n^4} = 96\delta^2\xi^2\sum_{n=1}^{\infty}\frac{1}{n^4},
    \end{equation*}
    and, using $\delta\xi\leq\alpha_1$ in this range together with $\sum_{n=1}^{\infty}n^{-4}=\pi^4/90$,
    \begin{equation*}
        h(\xi,\delta) \leq 96\alpha_1^2\sum_{n=1}^{\infty}\frac{1}{n^4} = 96\cdot\frac{1}{\pi^4}\left(\frac{90}{96}\right)\cdot\frac{\pi^4}{90} = 1,
    \end{equation*}
    thus $(3/\delta)m_{\delta}(\xi) \geq (4\pi^2-1)\xi^2$.

    \emph{Step 2: intermediate frequencies, $\alpha_1/\delta\leq\xi\leq\alpha_2/\delta$.} Since $y\mapsto y\coth(y)$ is increasing for $y>0$ and $y=2\pi\delta\xi$ is an increasing function of $\xi$, the map $\xi\mapsto m_\delta(\xi)$ is increasing on $(0,\infty)$, and therefore, by Step 1,
    \begin{equation*}
        \frac{3}{\delta}m_{\delta}(\xi) \geq \frac{3}{\delta}m_{\delta}\left(\frac{\alpha_1}{\delta}\right)\geq (4\pi^2-1)\left(\frac{\alpha_1}{\delta}\right)^2.
    \end{equation*}

    \emph{Step 3: large frequencies, $\xi\geq\alpha_2/\delta$.} Since $\coth(x)\geq1$ for any $x>0$,
    \begin{equation*}
        \frac{3}{\delta}m_{\delta}(\xi) = \frac{3}{\delta}\left(2\pi\xi\coth(2\pi\delta\xi)-\frac{1}{\delta}\right) \geq \frac{6\pi\xi}{\delta} - \frac{3}{\delta^2} \geq \frac{6\pi\xi}{\delta} - \frac{3\xi}{\delta\alpha_2} = \frac{5\pi\xi}{\delta}.
    \end{equation*}

    \emph{Step 4: conclusion.} Since $4\pi^2-1\geq5\pi$, Steps 1-3 give the stated bound in each of the three ranges.
\end{proof}

\section{Existence of suitable solitary waves for gILW}\label{Section: Existence of suitable solitary waves for GILW}
In this section, we prove the existence of strictly positive, even, and radially decreasing solitary waves for gILW. In fact, we establish a stronger result by considering a broader family of nonlinearities than those initially defining this equation. Our approach follows the one used for the gBO model (see Theorem $8.3$ in \cite{Angulo2009}) and relies on concentration-compactness and symmetric rearrangement: we first construct minimizers for a variational problem $N_{\delta,c,\lambda}$ associated with the nonlinearity $g_N$, ruling out vanishing and dichotomy in Lions concentration-compactness principle; symmetric decreasing rearrangement is then used to show that the resulting minimizer can be taken to be suitable, which simultaneously transfers the existence result to the two remaining nonlinearities $g_I$ and $g_L$.

Consider the following broader class of gILW equations
\begin{equation}\label{Eq: Broader class of GILW}
    \partial_tu + \frac{1}{\delta}\partial_xu + \mathcal{T}_\delta\partial_{x}^2u + \partial_x(g(u)) = 0, \quad (x,t)\in\R\times\R.
\end{equation}
Here, the nonlinear term $g:\R\to\R$ can be any of the following three types, which we label $g_I$, $g_N$ and $g_L$ for later reference:
\begin{equation*}
    g(t) = \left\{\begin{array}{cc}
        g_I(t):=t^{p+1}/(p+1), & 1\leq p\in\N, \\
        g_N(t):=t|t|^{p}/(p+1), & 1\leq p \in \R,\\
        g_L(t):=|t|^{p+1}/(p+1), & 1\leq p \in \R,
    \end{array}\right.
\end{equation*}
and the other parameters remain as before.

Solitary waves of \eqref{Eq: Broader class of GILW} are defined as $H^{1/2}$ solutions of
\begin{equation}\label{Eq: Broader GILW solitary wave equation}
    M_\delta\phi_{\delta,c} + c\phi_{\delta,c} - g(\phi_{\delta,c}) = 0.
\end{equation}
It is well-known that solutions of \eqref{Eq: Broader GILW solitary wave equation} can be obtained from solutions of a variational problem \cite{Arnesen2016, Weinstein1987}. In fact, depending on the nonlinear term $g$, we can use one of the following variational problems:
\begin{align*}
    I_{\delta,c,\lambda} &= \inf\{J_{\delta,c}(\psi): \psi\in H^{1/2}, H_I(\psi)=\lambda\}, \\
    N_{\delta,c,\lambda} &= \inf\{J_{\delta,c}(\psi): \psi\in H^{1/2}, H_N(\psi)=\lambda\}, \\
    L_{\delta,c,\lambda} &= \inf\{J_{\delta,c}(\psi): \psi\in H^{1/2}, H_L(\psi)=\lambda\}.
\end{align*}
Here, the objective functional $J_{\delta,c}$ is the one introduced in Section \ref{Section: The gILW equation and its solitary waves} (see also the general definition in Section \ref{Section: General Long Wave type models}), particularized to $M=M_\delta$. The restriction functionals depend on the choice of $g$ as follows:
\begin{align*}
    H_{I}(\psi) = &\int G_I(\psi)=\int \psi^{p+2}, \quad H_{N}(\psi)=\int G_N(\psi)=\int|\psi|^{p+2}, \\
    \quad &\text{ and } \quad H_{L}(\psi) = \int G_L(\psi)=\int \psi|\psi|^{p+1},
\end{align*}
where $G_I(t)=t^{p+2}$, $G_N(t)=|t|^{p+2}$ and $G_L(t)=t|t|^{p+1}$, so that each $G_k$, $k\in\{I,N,L\}$, satisfies $G_k'= (p+1)(p+2)g_k$ for $g_I, g_N$ and $g_L$ the three different kinds of nonlinearities.

The set of minimizers for the variational problems above are denoted as $\mathcal{I}_{\delta,c,\lambda}$, $\mathcal{N}_{\delta,c,\lambda}$, and $\mathcal{L}_{\delta,c,\lambda}$, respectively.

From now on, our goal is to prove the existence of suitable solutions for our variational problems. For technical reasons (see Observation \ref{Obs: why to choose N delta c lambda}), we begin to work with $N_{\delta,c,\lambda}$. 

The main tool in the proof of existence is the concentration-compactness principle of Lions (Lemma 1.1 in \cite{Lions1984a}). It states that any sequence of nonnegative functions with a fixed integral has a subsequence which falls into one of three alternatives: \emph{compactness}, in which, up to translations, the mass of the sequence remains concentrated in a bounded region; \emph{vanishing}, in which the mass spreads out, so that the amount contained in any interval of fixed length tends to zero uniformly; and \emph{dichotomy}, in which the mass splits into two pieces, each carrying a fixed fraction of the total, whose supports drift infinitely far apart. Existence of minimizers then follows once vanishing and dichotomy are ruled out.

The connection between the concentration-compactness principle and the existence of solutions for $N_{\delta,c,\lambda}$ is given through the so called concentration function $\rho$. Let $\psi_n$ be a minimizing sequence and define, for each $n$, the function $\rho_n:[0,\infty)\to[0,\lambda]$ as follows:
\begin{equation*}
    \rho_n(r)=\sup_{y\in\R}\int_{|x-y|<r}|\psi_n(x)|^{p+2}dx.
\end{equation*}
Since $\rho_n$ is a sequence of uniformly bounded and non-decreasing functions, Helly's selection theorem (Theorem \ref{Theorem: Helly selection theorem} in Appendix \ref{Appendix: Helly selection theorem}) guarantees the existence of a point-wise convergent sub-sequence, whose limit is also non-decreasing. Denoting this Helly's convergent sub-sequence also as $\rho_n$, we define the concentration function $\rho$, for any $r\in [0,\infty)$, as
\begin{equation*}
    \rho(r):=\lim_{n\to\infty}\rho_n(r)
\end{equation*}
and the concentration parameter $\alpha$ as
\begin{equation*}
    \alpha := \lim_{r\to\infty}\rho(r).
\end{equation*}

From now on, we work with the Helly-convergent subsequence, but retain the notation $\rho_n$ and $\psi_n$.

\begin{remark}\label{Obs: why to choose N delta c lambda}
    Note that in order to use Helly's selection theorem, we need the sequence $\rho_n$ to be uniformly bounded. In the variational problem $N_{\delta,c,\lambda}$, the positivity of $G_N$ immediately gives the bounds $0\leq\rho_n\leq\lambda$. This is the reason why we choose to work first with $N_{\delta,c,\lambda}$.
\end{remark}

Concentration-compactness principle applied to the sequence $|\psi_n|^{p+2}$ can be stated in terms of the concentration parameter. In fact, the compactness, vanishing, or dichotomy case occur if and only if $\alpha=\lambda$, $\alpha=0$, or $0<\alpha<\lambda$, respectively. As a consequence, since the compactness case leads to the existence of a non-trivial accumulation function, to prove the existence of minimizers for $N_{\delta,c,\lambda}$ it is enough to rule out the cases $0\leq\alpha<\lambda$.

As mentioned before, the application of the concentration-compactness principle to the variational problem $N_{\delta,c,\lambda}$ closely follows the gBO case presented in \cite{Angulo2009}. Differences arise only in arguments that rely on scaling or on the specific properties of $\mathcal{H}$. For the sake of completeness, we present here the whole procedure but skip some proofs that remain valid without any change.

The first step is to show that $N_{\delta,c,\lambda}$ is well-defined as an infimum.
\begin{lemma}\label{Lemma: Any minimizing sequence is bounded}
    We have $0<N_{\delta,c,\lambda}<\infty$. Moreover, any minimizing sequence for $N_{\delta,c,\lambda}$ is bounded in $H^{1/2}$.
\end{lemma}
\begin{proof}
    Let $\psi$ be any admissible function for $N_{\delta,c,\lambda}$. Then, by Sobolev embedding theorem and Corollary \ref{Corollary: J delta and norm 1/2},
    \begin{equation*}
        0 < \lambda^{2/(p+2)} = |\psi|_{p+2}^2 \lesssim \norm{\psi}_{1/2}^2 \lesssim_{\delta,c} J_{\delta,c}(\psi) < \infty.
    \end{equation*}
    In particular, 
    \begin{equation*}
        0 < \lambda^{2/(p+2)} \lesssim_{\delta,c} N_{\delta,c,\lambda} \leq J_{\delta,c}(\psi) < \infty.
    \end{equation*}
    The last part of the statement follows since for fixed $\delta,c>0$ the functional $\sqrt{J_{\delta,c}(\cdot)}$ is a norm on $H^{1/2}$ equivalent to the standard one (Corollary \ref{Corollary: J delta and norm 1/2}), therefore any minimizing sequence is bounded in $H^{1/2}$.
\end{proof}

The second step is to rule out vanishing. We obtain this as a consequence of the following lemma, whose proof does not depend on the particular function $m_\delta$ but on the energy space $H^{1/2}$. Since the energy space is the same for the gBO model, the proof of the lemma does not need any modification from its version stated in \cite{Angulo2009}; therefore, we omit the proof.

\begin{lemma}\label{Lemma: to rule out vanishing}
    Let $B>0$, $\epsilon>0$ be given. If $f\in H^{1/2}$ satisfies $\norm{f}_{1/2}\leq B$ and $|f|_{p+2}\geq\epsilon$, then, there exists some $\eta=\eta(B,\epsilon)>0$ such that
    \begin{equation*}
        \sup_{y\in\R}\int_{|x-y|<2}|f(x)|^{p+2}dx\geq\eta.
    \end{equation*}
\end{lemma}
\begin{proof}
    See Lemma 8.7 in \cite{Angulo2009}.
\end{proof}
Now we are in position to rule out vanishing.

\begin{theorem}[\textbf{Ruling out vanishing}]\label{Theorem: Ruling out vanishing}
    For any minimizing sequence $\psi_n$ for $N_{\delta,c,\lambda}$, we have $\alpha>0$.
\end{theorem}
\begin{proof}
    We know that $|\psi_n|_{p+2}^{p+2}=\lambda$ and $\norm{\psi_n}_{1/2}\leq B$ for some $B=B(\delta,c,\lambda)$. Then, Lemma \ref{Lemma: to rule out vanishing} implies that 
    \begin{equation*}
        \sup_{y\in\R}\int_{|x-y|<2}|\psi_n(x)|^{p+2}dx\geq\eta, \quad \forall n\in\N.
    \end{equation*}
    In particular, $\rho_n(2) \geq \eta > 0$ for all $n\in\N$. Letting $n\to\infty$ we get $\rho(2)=\lim_n\rho_n(2)\geq\eta$. Since each $\rho_n$ is non-decreasing, so is its limit $\rho$; therefore, $\alpha = \lim_{r\to\infty}\rho(r) \geq \rho(2) \geq\eta>0$.
\end{proof}

The third step is to rule out dichotomy. To do this we need some technical lemmas, the first one being a commutator estimate from \cite{Albert1999}.

\begin{lemma}\label{Lemma: Commutaror for Dichotomy}
    There exists a constant $A>0$ such that if $\theta$ is any continuously differentiable function with $\theta$ and $\theta'$ in $L^{\infty}$, and $f$ is any $L^2$ function, then
    \begin{equation*}
        \norm{[M_\delta,\theta]f} \leq A|\theta'|_{\infty}\norm{f},
    \end{equation*}
    where $[M_\delta,\theta]f$ denotes the commutator $M_\delta(\theta f)-\theta M_\delta f$.
\end{lemma}
\begin{proof}
    See Lemma 3.7 in \cite{Albert1999}.
\end{proof}

We introduce some notation that is going to be used in the next two lemmas.

Let $\varphi,\zeta\in C^{\infty}$ be such that $\varphi\equiv 1$ in $[-1,1]$, $\varphi(x)=0$ for $x\not\in[-2,2]$, and $\varphi^2+\zeta^2=1$ on $\R$. For each $w\in\R\backslash\{0\}$ and $r\in\R$, define
\begin{equation*}
    \varphi_{w,r}(x) := \varphi\left(\frac{x-r}{w}\right) \quad \text{ and } \quad \zeta_{w,r}(x) := \zeta\left(\frac{x-r}{w}\right).
\end{equation*}
The following two lemmas use the commutator estimate above to develop a good way to ``cut'' the functions $\psi_n$ into pieces that almost preserve the value of the functional $J_{\delta,c}$ and that are ``close'' to the admissible set of the variational problem $N_{\delta,c,\lambda}$.

\begin{lemma}\label{Lemma: Pre 1 Ruling out dichotomy}
    Let $\epsilon>0$, $B>0$ and $r\in\R$ be given. Assume that $\{\psi_n\}_n\subset H^{1/2}$ is such that $\norm{\psi_n}_{1/2}\leq B$. Then, there exists $w_0 = w_0(B,\epsilon)>0$ such that for all $w\geq w_0$ and $n\in\N$:
    \begin{equation*}
        \Big|J_{\delta,c}(\psi_n)-J_{\delta,c}(g_n)-J_{\delta,c}(h_n)\Big|<\epsilon,
    \end{equation*}
    where $g_n := \varphi_{w,r}\psi_n$ and $h_n := \zeta_{w,r}\psi_n$.
\end{lemma}
\begin{proof}
    \begin{align*}
        2J_{\delta,c}(g_n) &= \int\varphi_{w,r}\psi_n M_\delta(\varphi_{w,r}\psi_n) + c\int \varphi_{w,r}^2\psi_n^2 \\
        &= \int\varphi_{w,r}^2\psi_n M_\delta\psi_n + \int\varphi_{w,r}\psi_n [M_\delta,\varphi_{w,r}]\psi_n + c\int \varphi_{w,r}^2\psi_n^2 \\
        &= \int\varphi_{w,r}^2\{\psi_n M_\delta\psi_n + c\psi_n^2\} +\int\varphi_{w,r}\psi_n [M_\delta,\varphi_{w,r}]\psi_n.
    \end{align*}
    Similarly,
    \begin{equation*}
        2J_{\delta,c}(h_n) =\int\zeta_{w,r}^2\{\psi_n M_\delta\psi_n + c\psi_n^2\} +\int\zeta_{w,r}\psi_n [M_\delta,\zeta_{w,r}]\psi_n
    \end{equation*}
    and
    \begin{equation*}
        2J_{\delta,c}(\psi_n) = \int(\varphi_{w,r}^2+\zeta_{w,r}^2)\{\psi_n M_\delta\psi_n + c\psi_n^2 \}.
    \end{equation*}
    Therefore, using Lemma \ref{Lemma: Commutaror for Dichotomy},
    \begin{align*}
        \Big|J_{\delta,c}(\psi_n)-J_{\delta,c}(g_n)-J_{\delta,c}(h_n)\Big| &\lesssim \int|\varphi_{w,r}| |\psi_n [M_\delta,\varphi_{w,r}]\psi_n| + \int|\zeta_{w,r}| |\psi_n [M_\delta,\zeta_{w,r}]\psi_n| \\
        &\leq \norm{\psi_n}\left\{\norm{[M_\delta,\varphi_{w,r}]\psi_n} + \norm{[M_\delta,\zeta_{w,r}]\psi_n}\right\}\\
        &\lesssim_B |\varphi_{w,r}'|_{\infty}\norm{\psi_n} + |\zeta_{w,r}'|_{\infty}\norm{\psi_n} \\
        &\lesssim_B \frac{1}{w}.
    \end{align*}
\end{proof}

\begin{lemma}\label{Lemma: Pre 2 Ruling out dichotomy}
    For any $\epsilon>0$, there exists some $n_0\in\N$ and sequences $\{g_{n_0},g_{n_0+1},\cdots\}$ and $\{h_{n_0},h_{n_0+1},\cdots\}$ in $H^{1/2}$ such that for all $n\geq n_0$:
    \begin{enumerate}
        \item $|H_N(g_{n})-\alpha|<\epsilon$,
        \item $|H_N(h_n)-(\lambda-\alpha)|<\epsilon$,
        \item $J_{\delta,c}(\psi_n) \geq J_{\delta,c}(g_n) + J_{\delta,c}(h_n)-\epsilon$.
    \end{enumerate}
\end{lemma}
\begin{proof}
    From the definition of the concentration parameter $\alpha$ we know there is some $w_1\in\N$ such that for $w\geq w_1$:
    \begin{equation*}
        \alpha-\frac{\epsilon}{2} < \rho(w)\leq \rho(2w)\leq \alpha.
    \end{equation*}
    Let $w_0$ be the one from Lemma \ref{Lemma: Pre 1 Ruling out dichotomy} and fix $w=\max\{w_0,w_1\}$. Pick some $n_0\in\N$ such that for all $n\geq n_0$: $|\rho_n(w)-\rho(w)|<\epsilon/2$ and $|\rho_n(2w)-\rho(2w)|<\epsilon$. Now we have
    \begin{equation*}
        \alpha-\epsilon<\rho_n(w)\leq \rho_n(2w)<\alpha+\epsilon, \quad \forall n\geq n_0,
    \end{equation*}
    and therefore, by definition of $\rho_n$, we can find some $y_n$ such that
    \begin{equation*}
       \alpha-\epsilon < \int_{|x-y_n|<w}|\psi_n|^{p+2}dx \leq \int_{|x-y_n|<2w}|\psi_n|^{p+2}dx \leq \sup_{y\in\R}\int_{|x-y|<2w}|\psi_n|^{p+2}dx < \alpha+\epsilon.
    \end{equation*}
    For each $n\geq n_0$ define $g_n:= \varphi_{w,y_n}\psi_n$ and $h_n:= \zeta_{w,y_n}\psi_n$ and note they are in $H^{1/2}$ for $w\neq0$. Now, for any $n\geq n_0$, the inequalities
    \begin{align*}
        H_N(g_n) &= \int_\R |g_n|^{p+2} = \int_\R |\varphi_{w,y_n}\psi_n|^{p+2} \\
        &= \int_{y_n-2w}^{y_n+2w} \left|\varphi\left(\frac{x-y_n}{w}\right)\right|^{p+2}|\psi_n(x)|^{p+2}dx\\
        &\leq \int_{y_n-2w}^{y_n+2w} |\psi_n(x)|^{p+2}dx\\
        &<\alpha+\epsilon
    \end{align*}
    and
    \begin{align*}
        H_N(g_n) &= \int_{y_n-2w}^{y_n+2w} \left|\varphi\left(\frac{x-y_n}{w}\right)\right|^{p+2}|\psi_n(x)|^{p+2}dx\\
        &\geq \int_{y_n-w}^{y_n+w} \left|\varphi\left(\frac{x-y_n}{w}\right)\right|^{p+2}|\psi_n(x)|^{p+2}dx\\
        &= \int_{y_n-w}^{y_n+w} |\psi_n(x)|^{p+2}dx\\
        &> \alpha-\epsilon
    \end{align*}
    show item 1. The same argument gives item 2 from the inequalities
    \begin{align*}
        H_N(h_n) &= \int_\R |h_n|^{p+2} = \int_\R |\zeta_{w,y_n}\psi_n|^{p+2} \\
        &= \int_{|x-y_n|>w} \left|\zeta\left(\frac{x-y_n}{w}\right)\right|^{p+2}|\psi_n(x)|^{p+2}dx\\
        &\leq \int_{|x-y_n|>w}|\psi_n(x)|^{p+2}dx\\
        &= \int_{|x-y_n|>w}|\psi_n(x)|^{p+2}dx + \int_{|x-y_n|<w}|\psi_n(x)|^{p+2}dx - \int_{|x-y_n|<w}|\psi_n(x)|^{p+2}dx\\
        &= \lambda - \int_{|x-y_n|<w}|\psi_n(x)|^{p+2}dx\\
        &<\lambda + \epsilon - \alpha
    \end{align*}
    and
    \begin{align*}
        H_N(h_n) &= \int_{|x-y_n|>w} \left|\zeta\left(\frac{x-y_n}{w}\right)\right|^{p+2}|\psi_n(x)|^{p+2}dx\\
        &\geq \int_{|x-y_n|>2w} \left|\zeta\left(\frac{x-y_n}{w}\right)\right|^{p+2}|\psi_n(x)|^{p+2}dx\\
        &= \int_{|x-y_n|>2w} |\psi_n(x)|^{p+2}dx\\
        &= \int_{|x-y_n|>2w} |\psi_n(x)|^{p+2}dx + \int_{|x-y_n|<2w} |\psi_n(x)|^{p+2}dx - \int_{|x-y_n|<2w} |\psi_n(x)|^{p+2}dx\\ 
        &> \lambda -\epsilon-\alpha.
    \end{align*}
    Finally, taking $r=y_n$ in Lemma \ref{Lemma: Pre 1 Ruling out dichotomy} we get item 3.
\end{proof}

The last lemma needed to rule out dichotomy is the sub-additivity property of the value $N_{\delta,c,\lambda}$ with respect to $\lambda$.

\begin{lemma}\label{Lemma: Sub Aditivity to rule out dichotomy}
    For all $\theta\in(0,\lambda)$,
    \begin{equation*}
        N_{\delta,c,\lambda} < N_{\delta,c,\lambda-\theta} + N_{\delta,c,\theta}.
    \end{equation*}
\end{lemma}
\begin{proof}
    For $t>0$ and $\psi\in H^{1/2}$ we have $H_N(t\psi)=t^{p+2}H_N(\psi)$ and $J_{\delta,c}(t\psi)=t^2J_{\delta,c}(\psi)$. Hence, for any $\mu>0$, the map $\psi\mapsto\mu^{1/(p+2)}\psi$ is a bijection between the admissible sets of $N_{\delta,c,\lambda}$ and $N_{\delta,c,\mu\lambda}$, which multiplies $J_{\delta,c}$ by $\mu^{2/(p+2)}$; therefore $N_{\delta,c,\mu\lambda}=\mu^{2/(p+2)}N_{\delta,c,\lambda}$. In particular, $N_{\delta,c,\lambda}=\lambda^sN_{\delta,c,1}$ for every $\lambda>0$, where $s:=2/(p+2)\in(0,1)$, and $N_{\delta,c,1}>0$ by Lemma \ref{Lemma: Any minimizing sequence is bounded}. Since $(a+b)^s<a^s+b^s$ for all $a,b>0$, we obtain, for every $\theta\in(0,\lambda)$,
    \begin{equation*}
        N_{\delta,c,\lambda-\theta} + N_{\delta,c,\theta} = N_{\delta,c,1}\left((\lambda-\theta)^s+\theta^s\right) > N_{\delta,c,1}\,\lambda^s = N_{\delta,c,\lambda}.
    \end{equation*}
\end{proof}
Now we are in position to rule out dichotomy.
\begin{theorem}[\textbf{Ruling out dichotomy}]
    For every minimizing sequence $\psi_n$ for the variational problem $N_{\delta,c,\lambda}$, we have $\alpha=\lambda$.
\end{theorem}
\begin{proof}
    Suppose that $0<\alpha<\lambda$. We are going to show that this assumption implies the inequality
    \begin{equation*}
        N_{\delta,c,\lambda} \geq N_{\delta,c,\lambda-\alpha} + N_{\delta,c,\alpha}.
    \end{equation*}
    Given $\epsilon>0$, apply Lemma \ref{Lemma: Pre 2 Ruling out dichotomy} to get $n_0\in\N$ and the sequences $\{g_{n_0},g_{n_0+1},\cdots\}$ and $\{h_{n_0},h_{n_0+1},\cdots\}$ in $H^{1/2}$satisfying items 1,2, and 3 of the lemma. For all $n\geq n_0$, define
    \begin{equation*}
        \tilde{g}_n:=\frac{\alpha^{1/(p+2)}}{|g_n|_{p+2}}g_n \quad \text{and} \quad \tilde{h}_n:=\frac{(\lambda-\alpha)^{1/(p+2)}}{|h_n|_{p+2}}h_n.
    \end{equation*}
    Then, $H_N(\tilde{g}_n)=\alpha$ and $H_N(\tilde{h}_n)=\lambda-\alpha$, so $J_{\delta,c}(\tilde{g}_n)\geq N_{\delta,c,\alpha}$ and $J_{\delta,c}(\tilde{h}_n)\geq N_{\delta,c,\lambda-\alpha}$.
    
    Note that
    \begin{equation*}
        J_{\delta,c}(\tilde g_n) = \frac{\alpha^{2/(p+2)}}{|g_n|_{p+2}^2}J_{\delta,c}(g_n) \quad \text{and} \quad J_{\delta,c}(\tilde h_n) = \frac{(\lambda-\alpha)^{2/(p+2)}}{|h_n|_{p+2}^2}J_{\delta,c}(h_n),
    \end{equation*}
    thus,
    \begin{equation*}
        J_{\delta,c}(g_n)\geq\frac{|g_n|_{p+2}^2}{\alpha^{2/(p+2)}}N_{\delta,c,\alpha} \quad \text{and} \quad J_{\delta,c}(h_n)\geq\frac{|h_n|_{p+2}^2}{(\lambda-\alpha)^{2/(p+2)}}N_{\delta,c,\lambda-\alpha}.
    \end{equation*}
    Now, using item 3 in Lemma \ref{Lemma: Pre 2 Ruling out dichotomy} we get
    \begin{equation*}
        J_{\delta,c}(\psi_n) \geq \frac{|g_n|_{p+2}^2}{\alpha^{2/(p+2)}}N_{\delta,c,\alpha} + \frac{|h_n|_{p+2}^2}{(\lambda-\alpha)^{2/(p+2)}}N_{\delta,c,\lambda-\alpha} - \epsilon.
    \end{equation*}
    The items 1 and 2 in the same lemma imply then
    \begin{equation*}
        J_{\delta,c}(\psi_n) \geq \frac{(\alpha-\epsilon)^{2/(p+2)}}{\alpha^{2/(p+2)}}N_{\delta,c,\alpha} + \frac{(\lambda-\alpha-\epsilon)^{2/(p+2)}}{(\lambda-\alpha)^{2/(p+2)}}N_{\delta,c,\lambda-\alpha} - \epsilon.
    \end{equation*}
    Taking the limit $n\to\infty$ and then the limit $\epsilon\to0$ we get the desired inequality, which contradicts Lemma \ref{Lemma: Sub Aditivity to rule out dichotomy}. Therefore $\alpha\notin(0,\lambda)$, and since $\alpha>0$ by Theorem \ref{Theorem: Ruling out vanishing}, we conclude that $\alpha=\lambda$.
\end{proof}

At this point, we have proved that $\alpha=\lambda$, which corresponds to the compactness alternative in the concentration-compactness principle. The next theorem makes this assertion precise.

\begin{theorem}[\textbf{Compactness}]\label{Theorem: Compactness}
    Suppose $\alpha=\lambda$. Then, there exists a sequence of numbers $y_n$ with the property that for all $z\in (\lambda/2,\lambda)$, there exists some $w=w(z)\in\R$ and some $n_0=n_0(z)\in\N$ such that
    \begin{equation*}
        \int_{y_n-w}^{y_n+w}|\psi_n(x)|^{p+2}dx>z, \quad \forall n\geq n_0.
    \end{equation*}
\end{theorem}
\begin{proof}
    Fix an arbitrary $z\in(\lambda/2,\lambda)$. Then, there is some $\epsilon>0$ such that $z=\lambda-\epsilon$ and so we can take some $w_1=w_1(z)$ such that
    \begin{equation*}
        \lambda-\frac{\epsilon}{2}<\rho(w_1)<\lambda+\frac{\epsilon}{2}.
    \end{equation*}
    Also, there exists some $n_1=n_1(w_1,\epsilon)=n_1(z)$ such that for all $n\geq n_1$: $|\rho(w_1)-\rho_n(w_1)|<\epsilon/2$. Hence,
    \begin{equation*}
        \rho_n(w_1) = \sup_{y\in\R}\int_{y-w_1}^{y+w_1}|\psi_n|^{p+2}dx>\lambda-\epsilon = z, \quad \forall n\geq n_1.
    \end{equation*}
    In particular, for each $z$ there is some sequence $y_n(z)$ satisfying
    \begin{equation*}
        \int_{y_n(z)-w_1(z)}^{y_n(z)+w_1(z)}|\psi_n(x)|^{p+2}dx>z, \quad \forall n\geq n_1(z).
    \end{equation*}
    Moreover, using the same arguments for $z=\lambda/2$, we can find some $w_1=w_1(\lambda/2)\in\R$ and some $n_1=n_1(\lambda/2)\in\N$ that allow us to pick some sequence $y_n=y_n(\lambda/2)$ such that
    \begin{equation*}
        \int_{y_n(\lambda/2)-w_1(\lambda/2)}^{y_n(\lambda/2)+w_1(\lambda/2)}|\psi_n(x)|^{p+2}dx>\frac{\lambda}{2}, \quad \forall n\geq n_1(\lambda/2).
    \end{equation*}
    The key part of the proof is to note that the intervals $[y_n(\lambda/2)-w_1(\lambda/2),y_n(\lambda/2)+w_1(\lambda/2)]$ and $[y_n(z)-w_1(z),y_n(z)+w_1(z)]$ must overlap. In fact, otherwise we have
    \begin{equation*}
        \lambda = \int_\R|\psi_n|^{p+2} \geq \int_{y_n(z)-w_1(z)}^{y_n(z)+w_1(z)}|\psi_n|^{p+2} + \int_{y_n(\lambda/2)-w_1(\lambda/2)}^{y_n(\lambda/2)+w_1(\lambda/2)}|\psi_n|^{p+2} > z+\frac{\lambda}{2} > \lambda.
    \end{equation*}
    This overlapping implies that $|y_n(\lambda/2)-y_n(z)|\leq w_1(\lambda/2)+w_1(z)$ for any $z\in(\lambda/2,\lambda)$.

    Now, defining $w:= 2w_1(z)+w_1(\lambda/2)$ and $n_0 = \max\{n_1(z),n_1(\lambda/2)\}$, we obtain that for all $n\geq n_0$:
    \begin{equation*}
        y_n(\lambda/2)-y_n(z) \leq w-w_1(z) \quad \text{and} \quad y_n(z)-y_n(\lambda/2) \leq w-w_1(z),
    \end{equation*}
    and therefore
    \begin{equation*}
        [y_n(z)-w_1(z),y_n(z)+w_1(z)]\subset [y_n(\lambda/2)-w,y_n(\lambda/2)+w].
    \end{equation*}
    Then, for any $z\in(\lambda/2,\lambda)$, there exists some $w=w(z)\in\R$ and some $n_0=n_0(z)\in\N$ such that
    \begin{equation*}
        \int_{y_n(\lambda/2)-w}^{y_n(\lambda/2)+w}|\psi_n(x)|^{p+2}dx\geq \int_{y_n(z)-w_1(z)}^{y_n(z)+w_1(z)}|\psi_n(x)|^{p+2}dx > z, \quad \forall n\geq n_0.
    \end{equation*}
\end{proof}

Now that we know there is a tight sub-sequence of translations of the minimizing sequence $\psi_n$, we can use this to overcome the lack of compactness in the embedding $H^{1/2}(\R)\embed L^{p+2}(\R)$ in order to obtain a convergent sub-sequence of translations and the existence of solution for $N_{\delta,c,\lambda}$.

\begin{theorem}[\textbf{Existence of minimizers for} $\boldsymbol{N_{\delta,c,\lambda}}$]\label{Theorem: Existence of solutions for N delta c lambda}
    Let $c,\delta,\lambda>0$ be fixed and consider a minimizing sequence $\psi_n$ for $N_{\delta,c,\lambda}$. Then, there is a sub-sequence $\psi_{n_k}$ and a sequence of real numbers $y_{n_k}$ such that $\psi_{n_k}(\cdot+y_{n_k})$ converge in $H^{1/2}$ to some $\psi\in H^{1/2}$. The limit $\psi$ is a minimizer for $N_{\delta,c,\lambda}$.
\end{theorem}
\begin{proof}
    \emph{Step 1: a weak limit $\psi$ of translates of $\psi_n$.} For each $k\in\N$ (taken large enough, depending on $\lambda$, if necessary), let $z=\lambda-1/k>\lambda/2$. By Theorem \ref{Theorem: Compactness}, there exists a sequence $y_n$ with the property that for any $k$, we can find some $w_k\in\R$ and some $n_k\in\N$ such that for all $n\geq n_k$:
    \begin{equation}\label{Eq: 0 in proof existence of minimizers for N}
        \int_{y_n-w_k}^{y_n+w_k}|\psi_n(x)|^{p+2}dx > \lambda-\frac{1}{k}.
    \end{equation}
    Define $\tilde\psi_n(x):= \psi_n(x+y_n)$. Since $\|\tilde\psi_n\|_{1/2} = \|\psi_n\|_{1/2} \leq B$ (see Lemma \ref{Lemma: Any minimizing sequence is bounded}), Banach-Alaoglu's theorem guarantees the existence of some $\psi\in H^{1/2}$ and a sub-sequence, which we still denote as $\tilde\psi_n$, such that $\tilde\psi_n\rightharpoonup \psi$ in $H^{1/2}$. We show below that $\psi$ is a minimizer for $N_{\delta,c,\lambda}$.

    \emph{Step 2: $\psi$ is admissible and $J_{\delta,c}(\psi)\leq N_{\delta,c,\lambda}$.} Sobolev's embedding theorem implies that $\tilde\psi_n\rightharpoonup\psi$ in $L^{p+2}$, so that the weak lower semicontinuity of the $L^{p+2}$-norm gives
    \begin{equation}\label{Eq: 1 in proof existence of minimizers for N}
        \int_\R |\psi|^{p+2} \leq \liminf_{n} \int_\R |\tilde\psi_n|^{p+2} = \lambda.
    \end{equation}
    Moreover, since $\sqrt{J_{\delta,c}(\cdot)}$ is a norm in $H^{1/2}$ equivalent to the standard one (see Corollary \ref{Corollary: J delta and norm 1/2}), the functional $J_{\delta,c}$ is weak lower semicontinuous in $H^{1/2}$; thus
    \begin{equation}\label{Eq: 2 in proof existence of minimizers for N}
        J_{\delta,c}(\psi) \leq  \liminf_{n}J_{\delta,c}(\tilde\psi_{n}) = N_{\delta,c,\lambda}.
    \end{equation}
    By \eqref{Eq: 1 in proof existence of minimizers for N} and \eqref{Eq: 2 in proof existence of minimizers for N}, to show that
    $\psi$ is a minimizer for $N_{\delta,c,\lambda}$ we only need to prove that $H_{N}(\psi)\geq\lambda$. This is the key part of the proof, where we use the tightness obtained from the concentration-compactness principle, stated in equation \eqref{Eq: 0 in proof existence of minimizers for N}.

    \emph{Step 3: strong convergence on bounded intervals.} Fix $k\in\N$ and let $I_k:=(-w_k,w_k)$. The restrictions $\tilde\psi_n|_{I_k}$ are bounded by $B$ in $H^{1/2}(I_k)$, and the embedding $H^{1/2}(I_k)\embed L^{p+2}(I_k)$ is compact by Rellich-Kondrachov's theorem. Hence every sub-sequence of $\tilde\psi_n|_{I_k}$ has a further sub-sequence converging strongly in $L^{p+2}(I_k)$, and its limit must be $\psi|_{I_k}$, since $\tilde\psi_n\rightharpoonup\psi$ in $L^{p+2}(\R)$. Therefore, the whole sequence satisfies $\tilde\psi_n\to\psi$ in $L^{p+2}(I_k)$, for each $k\in\N$.

    \emph{Step 4: conclusion.} Since \eqref{Eq: 0 in proof existence of minimizers for N} reads $\int_{-w_k}^{w_k}|\tilde\psi_n|^{p+2}>\lambda-1/k$ for $n\geq n_k$, Step 3 gives
    \begin{equation*}
        \int_\R |\psi|^{p+2}\geq \int_{-w_k}^{w_k}|\psi|^{p+2} = \lim_{n\to\infty} \int_{-w_k}^{w_k}|\tilde\psi_{n}|^{p+2} \geq \lambda - \frac{1}{k}, \quad \forall k\in\N,
    \end{equation*}
    so that $H_N(\psi)\geq\lambda$, and $\psi$ is a minimizer for $N_{\delta,c,\lambda}$. Finally, $\tilde\psi_{n}\rightharpoonup\psi$ in $H^{1/2}$ together with $J_{\delta,c}(\tilde\psi_n)\to J_{\delta,c}(\psi)$ upgrades this weak convergence to $\tilde\psi_{n}\to\psi$ in $H^{1/2}$.
\end{proof}

After we have attained our goal of proving the existence of minimizers for the variational problem $N_{\delta,c,\lambda}$, we use symmetric decreasing rearrangements to get the same result for the problems $I_{\delta,c,\lambda}$ and $L_{\delta,c,\lambda}$. This is possible because of the following lemma.

\begin{lemma}\label{Lemma: J delta and rearrangement}
    Let $\psi\in H^{1/2}$ and $\psi^*:=|\psi|^*$ (the symmetric decreasing rearrangement, see Appendix \ref{Appendix: Symmetric decreasing rearrangement}). Then,
    \begin{equation*}
        J_{\delta,c}(\psi^*) \leq J_{\delta,c}(\psi).
    \end{equation*}
\end{lemma}
\begin{proof}
    \emph{Step 1: reduction to a symbol inequality.} We want to show
    \begin{equation*}
        \frac{1}{2}\int\psi^*M_\delta\psi^* + \frac{c}{2} \int(\psi^*)^2 \leq \frac{1}{2}\int\psi M_\delta\psi + \frac{c}{2} \int\psi^2,
    \end{equation*}
    so that, since the rearrangement preserves $L^p$-norms (see Proposition \ref{Prop: Rearrangement properties}), it suffices to show
    \begin{equation*}
        \int m_\delta(\xi)|\fourier{\psi^*}(\xi)|^2d\xi \leq \int m_\delta(\xi)|\fourier{\psi}(\xi)|^2d\xi.
    \end{equation*}
    Since $m_\delta(\xi)=2\pi\xi\coth(2\pi\delta\xi)-1/\delta$ differs from $2\pi\xi\coth(2\pi\delta\xi)$ only by an additive constant, and the rearrangement preserves the $L^2$-norm (so $\int|\fourier{\psi^*}|^2=\int|\fourier{\psi}|^2$), this constant cancels on both sides; it is therefore enough to prove the inequality with $m_\delta(\xi)$ replaced by $2\pi\xi\coth(2\pi\delta\xi)$.

    \emph{Step 2: a Riesz rearrangement estimate.} For any $\nu>0$, define $\kappa_\nu$ as the function satisfying
    \begin{equation*}
        \fourier{\kappa_\nu}(\xi) = \frac{1}{\nu+2\pi\xi\coth(2\pi\delta\xi)}.
    \end{equation*}
    Since $2\pi\xi\coth(2\pi\delta\xi)\geq 2\pi|\xi|$, $\fourier{\kappa_\nu}\in L^2$, so $\kappa_\nu\in L^2$ is well defined. Moreover, $\kappa_\nu$ is the convolution kernel of the inverse of the operator with symbol $\nu+2\pi\xi\coth(2\pi\delta\xi)$, and by Lemma 1 in \cite{Albert1995} it is given by a series $\kappa_\nu(x)=\sum_{j\geq0}a_je^{-\theta_j|x|/\delta}$, with $a_j>0$ and $\theta_j>0$ for all $j$. In particular, $\kappa_\nu$ is strictly positive, even and strictly decreasing in $|x|$, so that $\kappa_\nu^*=\kappa_\nu$.

    Since any $f\in L^2$ satisfies $\mu(\{x:f(x)>y\})<\infty$ for $y>0$, and $\psi,\kappa_\nu\in L^2$, the Riesz rearrangement lemma (see Proposition \ref{Prop: Rearrangement properties}), together with $\kappa_\nu^*=\kappa_\nu$, gives
    \begin{align}
        \left|(\psi*\kappa_\nu*\psi)(0)\right| &\leq (\psi^**\kappa_\nu*\psi^*)(0) \notag \\
        \left|\int\psi(x)\kappa_\nu*\psi(-x)dx\right| &\leq \int\psi^*(x)\kappa_\nu*\psi^*(-x)dx \label{Eq: eq 1 in proof Lemma for rearrangement 1}
    \end{align}
    where the symmetry of $\kappa_\nu$ and the translation invariance of the integral imply
    \begin{align*}
        \kappa_\nu*\psi(-x) &= \int \kappa_\nu(-x-y)\psi(y)dy
        = \int \kappa_\nu(x+y)\psi(y)dy \\
        &= \int \kappa_\nu(y)\psi(y-x)dy
        = \int \kappa_\nu(-y)\psi(y-x)dy \\
        &= \int \kappa_\nu(x-y)\psi(y)dy
        = \kappa_\nu*\psi(x)
    \end{align*}
    and similarly $\kappa_\nu*\psi^*(-x)=\kappa_\nu*\psi^*(x)$. From \eqref{Eq: eq 1 in proof Lemma for rearrangement 1}, using that the Fourier transform is a unitary operator on $L^2$,
    \begin{align*}
        \left|\int\psi(x)\kappa_\nu*\psi(x)dx\right| &\leq \int\psi^*(x)\kappa_\nu*\psi^*(x)dx \\
        \int \fourier{\kappa_\nu}(\xi)|\fourier{\psi}(\xi)|^2d\xi &\leq \int \fourier{\kappa_\nu}(\xi)|\fourier{\psi^*}(\xi)|^2d\xi.
    \end{align*}

    \emph{Step 3: passing to the limit $\nu\to\infty$.} For any $\nu>0$, Step 2 gives
    \begin{equation*}
        \int(\nu-\nu^2\fourier{\kappa_\nu}(\xi))|\fourier{\psi^*}(\xi)|^2d\xi \leq \int(\nu-\nu^2\fourier{\kappa_\nu}(\xi))|\fourier{\psi}(\xi)|^2d\xi.
    \end{equation*}
    Now,
    \begin{align*}
        \nu-\nu^2\fourier{\kappa_\nu}(\xi) &= \nu\left[1-\frac{\nu}{\nu+2\pi\xi\coth(2\pi\delta\xi)}\right]
        = \nu\left[\frac{2\pi\xi\coth(2\pi\delta\xi)}{\nu+2\pi\xi\coth(2\pi\delta\xi)}\right]\\
        &\leq 2\pi\xi\coth(2\pi\delta\xi)
        \lesssim_\delta 1 + |\xi|,
    \end{align*}
    so that, since $\psi,\psi^*\in H^{1/2}$ (see Lemma \ref{Lemma rearrangement in H1/2}), dominated convergence as $\nu\to\infty$ yields
    \begin{equation*}
        \int2\pi\xi\coth(2\pi\delta\xi))|\fourier{\psi^*}(\xi)|^2d\xi \leq \int2\pi\xi\coth(2\pi\delta\xi))|\fourier{\psi}(\xi)|^2d\xi.
    \end{equation*}

    \emph{Step 4: conclusion.} Subtracting $\frac{1}{\delta}\int|\fourier{\psi^*}(\xi)|^2d\xi=\frac{1}{\delta}\int|\fourier{\psi}(\xi)|^2d\xi$ (equal by preservation of the $L^2$-norm) from both sides of the inequality in Step 3 gives $\int m_\delta(\xi)|\fourier{\psi^*}(\xi)|^2d\xi\leq\int m_\delta(\xi)|\fourier{\psi}(\xi)|^2d\xi$, which is the desired inequality.
\end{proof}

Let us now see that a rearrangement of any function in the set of minimizers $\mathcal{N}_{\delta,c,\lambda}$ gives an element in $\mathcal{N}_{\delta,c,\lambda}\cap\mathcal{I}_{\delta,c,\lambda}\cap\mathcal{L}_{\delta,c,\lambda}$.

\begin{proposition}\label{Proposition: Existence of minimizer of I and L}
    Let $\psi\in\mathcal{N}_{\delta,c,\lambda}$. Then, $\psi^*\in\mathcal{I}_{\delta,c,\lambda}$, $\psi^*\in\mathcal{N}_{\delta,c,\lambda}$ and $\psi^*\in\mathcal{L}_{\delta,c,\lambda}$. Moreover, $N_{\delta,c,\lambda}=I_{\delta,c,\lambda}=L_{\delta,c,\lambda}$.
\end{proposition}
\begin{proof}
    We begin showing that $\psi^*\in\mathcal{N}_{\delta,c,\lambda}$. Since the rearrangement preserves the $L^{p+2}$-norm, it follows that
    \begin{equation*}
        \lambda = H_N(\psi) = \int|\psi|^{p+2} = \int |\psi^*|^{p+2} = H_N(\psi^*),
    \end{equation*}
    then $\psi^*$ is admissible for $N_{\delta,c,\lambda}$, and by Lemma \ref{Lemma: J delta and rearrangement}:
    \begin{equation*}
        N_{\delta,c,\lambda} \leq J_{\delta,c}(\psi^*) \leq J_{\delta,c}(\psi) = N_{\delta,c,\lambda},
    \end{equation*}
    thus $\psi^*\in\mathcal{N}_{\delta,c,\lambda}$.
    
    Now, again, since the rearrangement preserves the $L^{p+2}$-norm, it follows that
    \begin{equation*}
        \lambda=H_N(\psi) = \int|\psi|^{p+2} = \int |\psi^*|^{p+2} = \int (\psi^*)^{p+2} = H_I(\psi^*)
    \end{equation*}
    and
    \begin{equation*}
        \lambda=H_N(\psi) = \int|\psi|^{p+2} = \int |\psi^*|^{p+2} = \int \psi^*|\psi^*|^{p+1} = H_L(\psi^*),
    \end{equation*}
    then $\psi^*$ is admissible for $I_{\delta,c,\lambda}$ and $L_{\delta,c,\lambda}$. In particular,
    \begin{equation*}
        I_{\delta,c,\lambda} \leq J_{\delta,c}(\psi^*)=N_{\delta,c,\lambda} \quad \text{ and } \quad L_{\delta,c,\lambda} \leq J_{\delta,c}(\psi^*)=N_{\delta,c,\lambda}.
    \end{equation*}
    Now, suppose that $\psi^*\notin \mathcal{I}_{\delta,c,\lambda}$. Then, there exists some $\varphi\in H^{1/2}$ such that $H_I(\varphi)=\lambda$ and $J_{\delta,c}(\varphi)< J_{\delta,c}(\psi^*)$. Define the function
    \begin{equation*}
        f=\frac{\lambda^{1/(p+2)}}{|\varphi|_{p+2}}\varphi,
    \end{equation*}
    then
    \begin{equation*}
        H_N(f) = \frac{\lambda}{|\varphi|_{p+2}^{p+2}}\int|\varphi|^{p+2} =\lambda,
    \end{equation*}
    so $f$ is admissible for $N_{\delta,c,\lambda}$ and
    \begin{equation*}
        J_{\delta,c}(f) = \frac{\lambda^{2/(p+2)}}{|\varphi|_{p+2}^2}J_{\delta,c}(\varphi) \leq J_{\delta,c}(\varphi),
    \end{equation*}
    where the inequality follows since $\lambda = H_I(\varphi) \leq H_N(\varphi) = |\varphi|_{p+2}^{p+2}$. Therefore,
    \begin{equation*}
        J_{\delta,c}(\psi^*) = N_{\delta,c,\lambda} \leq J_{\delta,c}(f) \leq J_{\delta,c}(\varphi) < J_{\delta,c}(\psi^*),
    \end{equation*}
    a contradiction. Thus, $\psi^*\in\mathcal{I}_{\delta,c,\lambda}$. 
    
    Similarly, since $H_L(\varphi)\leq H_N(\varphi)$, the same argument shows that $\psi^*\in\mathcal{L}_{\delta,c,\lambda}$.

    Since $\psi^*\in\mathcal{N}_{\delta,c,\lambda}\cap\mathcal{I}_{\delta,c,\lambda}\cap\mathcal{L}_{\delta,c,\lambda}$, it follows that $N_{\delta,c,\lambda}=J_{\delta,c}(\psi^*)=I_{\delta,c,\lambda}=L_{\delta,c,\lambda}$.
\end{proof}

We are now in position to prove the main theorem of this section.

\begin{theorem}[\textbf{Existence of suitable solitary waves}]\label{Th: Existence of suitable solitary waves}
    For any velocity $c>0$, any depth parameter $\delta>0$, any admissible order of nonlinearity $p\geq1$ (as in the definition of $g$), and any of the three nonlinear terms $g\in\{g_I,g_N,g_L\}$ described above, there exists an even, strictly positive, and radially decreasing solitary wave $\phi_{\delta,c}$ of equation \eqref{Eq: Broader class of GILW}.
\end{theorem}
\begin{proof}
    Fix arbitrary values $\lambda,c,\delta>0$ and $p\geq1$. Theorem \ref{Theorem: Existence of solutions for N delta c lambda} and Proposition \ref{Proposition: Existence of minimizer of I and L} imply the existence of even, nonnegative, and radially decreasing minimizers for $N_{\delta,c,\lambda}, I_{\delta,c,\lambda}$ and $L_{\delta,c,\lambda}$; strict positivity is shown below.

    Let $\psi_{c,I}\in \mathcal{I}_{\delta,c,\lambda}, \psi_{c,N}\in \mathcal{N}_{\delta,c,\lambda}$ and $\psi_{c,L}\in \mathcal{L}_{\delta,c,\lambda}$ be even, nonnegative, and radially decreasing minimizers. By Lagrange multiplier theorem, it follows that these minimizers satisfy
    \begin{equation}\label{Eq: Minimizer equation}
        M_\delta\psi_{c,k} + c\psi_{c,k} = \theta_{\delta,c,\lambda,k} G_k'(\psi_{c,k}) = \theta_{\delta,c,\lambda,k}(p+2)(p+1) g_k(\psi_{c,k}) = \theta_{\delta,c,\lambda,k}(p+2)\psi_{c,k}^{p+1},
    \end{equation}
    where $k\in\{I,N,L\}$, $\theta_{\delta,c,\lambda,k}\in\R$ is the Lagrange multiplier and $g_k$, a priori, depends on whether $\psi_{c,k}$ is in $\mathcal{I}_{\delta,c,\lambda}$, $\mathcal{N}_{\delta,c,\lambda}$ or $\mathcal{L}_{\delta,c,\lambda}$; but in fact, since $\psi_{c,k}$ is a nonnegative real-valued function, in any case we get $g_k(\psi_{c,k})=\psi_{c,k}^{p+1}/(p+1)$.

    Note that since $\psi_{c,k}\in H^{1/2}$, Sobolev's theorem and a bootstrap argument show $\psi_{c,k}\in H^{\infty}$.

    Let us see that $\theta_{\delta,c,\lambda,k}>0$. Since $\psi_{c,k}$ is a minimizer, $J_{\delta,c}(\psi_{c,k})$ equals $N_{\delta,c,\lambda}$, $I_{\delta,c,\lambda}$ or $L_{\delta,c,\lambda}$ according to the case $k$, and by Proposition \ref{Proposition: Existence of minimizer of I and L} these three values coincide and are positive (Lemma \ref{Lemma: Any minimizing sequence is bounded}); hence $J_{\delta,c}(\psi_{c,k})>0$ in every case. Multiplying equation \eqref{Eq: Minimizer equation} by $\psi_{c,k}$ and integrating we get
    \begin{equation*}
        0<2J_{\delta,c}(\psi_{c,k}) = \theta_{\delta,c,\lambda,k}(p+2)\int \psi_{c,k}^{p+2} = \theta_{\delta,c,\lambda,k}(p+2)\lambda
    \end{equation*}
    thus $\theta_{\delta,c,\lambda,k}>0$.

    Let us also see that $\psi_{c,k}$ is strictly positive. By \eqref{Eq: Minimizer equation}, $\psi_{c,k}=\theta_{\delta,c,\lambda,k}(p+2)\,\kappa*\psi_{c,k}^{p+1}$, where $\kappa$ is the convolution kernel of the inverse of $M_\delta+c$, that is, of the operator with symbol $2\pi\xi\coth(2\pi\delta\xi)+(c-1/\delta)$. Since $c-1/\delta>-1/\delta$, Lemma 1 in \cite{Albert1995} shows that $\kappa(x)>0$ for all $x\in\R$. As $\psi_{c,k}^{p+1}\geq0$ is not identically zero and $\theta_{\delta,c,\lambda,k}>0$, it follows that $\psi_{c,k}(x)>0$ for all $x\in\R$.

    Note that in any case, $g_k(C\psi) = C^{p+1}g_k(\psi)$ for any positive constant $C$, so defining $\eta_{\delta,c,\lambda,k}^p:=\theta_{\delta,c,\lambda,k}(p+2)(p+1)$ and $\phi_{c,k}=\eta_{\delta,c,\lambda,k}\psi_{c,k}$, it follows from equation \eqref{Eq: Minimizer equation} that
    \begin{equation*}
        M_\delta\phi_{c,k} + c\phi_{c,k} = \frac{1}{p+1}\phi_{c,k}^{p+1},
    \end{equation*}
    thus $\phi_{c,k}$ is the desired solitary wave in the respective case $k=I,N$ or $L$.
\end{proof}

\section{Convergence of solitary waves}\label{Section: Convergence of solitary waves}

In this section, we prove that, in the deep water limit $\delta\to\infty$, the solitary waves of gILW converge to solitary waves of gBO both in the energy space and uniformly; while, in the shallow water limit $\delta\to0$, suitably rescaled solitary waves of gILW converge to solitary waves of gKdV in $L^q(\R)$ for every $q\in[2,\infty]$, in particular, uniformly. From now on, consider the values $\lambda>0$, $c>0$ and an admissible $p\geq1$ to be arbitrary and fixed.

For each $\delta>0$, let $\phi_{\delta,c}$ be a suitable solitary wave of gILW with velocity $c$ and depth parameter $\delta$, given by Theorem \ref{Th: Existence of suitable solitary waves}. Recall that it has been obtained as
\begin{equation*}
    \phi_{\delta,c} = \eta_{\delta,c,\lambda}\psi_{\delta,c,\lambda}, \qquad \eta_{\delta,c,\lambda}^p:= \theta_{\delta,c,\lambda}(p+2)(p+1),
\end{equation*}
where $\psi_{\delta,c,\lambda}$ is a suitable solution of the variational problem $I_{\delta,c,\lambda}$ and $\theta_{\delta,c,\lambda}>0$ is the associated Lagrange multiplier.

In the shallow water limit, gILW converges to gKdV only after a rescaling (see Section \ref{Section: Introduction}). Accordingly, we consider the \emph{rescaled solitary waves}
\begin{equation*}
    \phi_{\delta,c}^{(s)} := \left(\frac{3}{\delta}\right)^{1/p}\phi_{\delta,\delta c/3}.
\end{equation*}
Straightforward computations show that $\phi_{\delta,c}^{(s)}$ satisfies
\begin{equation}\label{Eq: rescaled solitary wave equation}
    M_\delta^{(s)}\phi_{\delta,c}^{(s)} + c\phi_{\delta,c}^{(s)} - \frac{1}{p+1}(\phi_{\delta,c}^{(s)})^{p+1} = 0, \qquad M_\delta^{(s)}:=\frac{3}{\delta}M_\delta,
\end{equation}
where the symbol $m_\delta^{(s)}:=(3/\delta)m_\delta$ of $M_\delta^{(s)}$ converges point-wise to $4\pi^2\xi^2$, the symbol of gKdV, as $\delta\to0$. Moreover, the rescaled functional $J_{\delta,c}^{(s)}:=(3/\delta)J_{\delta,\delta c/3}$ is the functional of Section \ref{Section: General Long Wave type models} particularized to $M=M_\delta^{(s)}$, so that the variational problem $I_{\delta,c,\lambda}^{(s)}:=(3/\delta)I_{\delta,\delta c/3,\lambda}$ has the same minimizers as $I_{\delta,\delta c/3,\lambda}$. Writing
\begin{equation*}
    \psi_{\delta,c,\lambda}^{(s)}:=\psi_{\delta,\delta c/3,\lambda}, \qquad \theta_{\delta,c,\lambda}^{(s)}:=\frac{3}{\delta}\theta_{\delta,\delta c/3,\lambda}, \qquad \eta_{\delta,c,\lambda}^{(s)}:=\left(\frac{3}{\delta}\right)^{1/p}\eta_{\delta,\delta c/3,\lambda},
\end{equation*}
we have $\phi_{\delta,c}^{(s)}=\eta_{\delta,c,\lambda}^{(s)}\psi_{\delta,c,\lambda}^{(s)}$, where $\psi_{\delta,c,\lambda}^{(s)}$ is a suitable solution of $I_{\delta,c,\lambda}^{(s)}$ with Lagrange multiplier $\theta_{\delta,c,\lambda}^{(s)}>0$, and $(\eta_{\delta,c,\lambda}^{(s)})^p=\theta_{\delta,c,\lambda}^{(s)}(p+2)(p+1)$. In this way, every statement below about the rescaled objects is a statement about the solitary waves of gILW. Objects without superscript always refer to gILW, while the superscript $(s)$ is reserved for the rescaled ones; in statements that apply to both families simultaneously, we write $\phi_{\delta,c}^{(j)}$, $j\in\{s,d\}$, with the convention that $\phi_{\delta,c}^{(d)}:=\phi_{\delta,c}$, and similarly for $\psi$, $\theta$, $\eta$, $M_\delta$, $m_\delta$, $J_{\delta,c}$ and the convolution kernel $\mathcal{K}_{\delta,c}$ introduced in Section \ref{Section: Uniform convergence of solitary waves}.

Now, we define the natural candidates to be the limit of our solitary waves, following the notation $\infty$ (deep water limit, i.e.\ gBO) and $0$ (shallow water limit, i.e.\ gKdV) already introduced in Section \ref{Section: Known results for gKdV and gBO}. Recall from Section \ref{Section: General Long Wave type models} that $J_{\infty,c}$ and $J_{0,c}$ denote the functional $J_c$ particularized to $M=\mathcal{H}\partial_x$ (gBO) and $M=-\partial_x^2$ (gKdV), respectively. Let $\phi_{\infty,c}$ and $\phi_{0,c}$ be even, strictly positive, and radially decreasing solitary waves of gBO and gKdV, respectively, with velocity $c$. Also, assume that $\phi_{\infty,c}$ and $\phi_{0,c}$ have the form
\begin{equation*}
    \phi_{\infty,c} = \eta_{\infty,c,\lambda}\psi_{\infty,c,\lambda} \quad \text{and} \quad \phi_{0,c} = \eta_{0,c,\lambda}\psi_{0,c,\lambda},
\end{equation*}
where $\psi_{\infty,c,\lambda}$ and $\psi_{0,c,\lambda}$ are solutions for the variational problems $I_{\infty,c,\lambda}$ and $I_{0,c,\lambda}$, respectively, and the values $\eta_{\infty,c,\lambda}$ and $\eta_{0,c,\lambda}$ are defined as
\begin{equation*}
    \eta_{\infty,c,\lambda}^p:= \theta_{\infty,c,\lambda}(p+2)(p+1)\quad \text{and}
    \quad \eta_{0,c,\lambda}^p:=\theta_{0,c,\lambda}(p+2)(p+1),
\end{equation*}
where $\theta_{\infty,c,\lambda},\theta_{0,c,\lambda}>0$ are Lagrange multipliers.

The existence of solitary waves of gBO and gKdV with this particular variational structure is justified as follows. For gBO, this is Theorem \ref{Theorem: Existence for gBO}, and, by Proposition \ref{Prop: gBO waves are FL ground states}, the resulting solitary wave $\phi_{\infty,c}$ is unique. For gKdV, the structure of the dispersive term $M=-\partial_x^2$ allows us to run the very same proof scheme used in Section \ref{Section: Existence of suitable solitary waves for GILW} for the existence of suitable solitary waves, which again yields solitary waves of the stated form.

The strategy we use to show the uniform convergence of solitary waves is based on uniform estimates of the waves $\phi_{\delta,c}^{(j)}$, $j\in\{s,d\}$. In order to get these estimates, we use the information of $\psi_{\delta,c,\lambda}^{(j)}$ as a minimizer for the associated variational problem. Of course, since the relation between the waves $\psi_{\delta,c}^{(j)}$ and $\phi_{\delta,c}^{(j)}$ is given through the multiplier $\theta_{\delta,c,\lambda}^{(j)}$, we need to uniformly estimate these multipliers. The following lemma gives the desired bounds.

\begin{lemma}\label{Lemma: thetas are bounded}
    There are some $C_1,C_2>0$, depending only on $c$ and $\lambda$, such that $0<\theta_{\delta,c,\lambda}<C_1$ and $0<\theta_{\delta,c,\lambda}^{(s)}<C_2$ for any $\delta>0$. Moreover, the map $\delta\mapsto\theta_{\delta,c,\lambda}^{(s)}$ is non-increasing on $(0,\infty)$.
\end{lemma}
\begin{proof}
    Functions $\psi_{\delta,c,\lambda}$ and $\psi_{\delta,c,\lambda}^{(s)}$ satisfy equations
    \begin{equation*}
        M_{\delta}\psi_{\delta,c,\lambda} + c\psi_{\delta,c,\lambda} = \theta_{\delta,c,\lambda}(p+2)(\psi_{\delta,c,\lambda})^{p+1} \quad \text{and}
    \end{equation*}
    \begin{equation*}
        M_{\delta}^{(s)}\psi_{\delta,c,\lambda}^{(s)} + c\psi_{\delta,c,\lambda}^{(s)} = \theta_{\delta,c,\lambda}^{(s)}(p+2)(\psi_{\delta,c,\lambda}^{(s)})^{p+1},
    \end{equation*}
    respectively. Then, multiplying the above equations by $\psi_{\delta,c,\lambda}$ and $\psi_{\delta,c,\lambda}^{(s)}$, and integrating over $\R$, we get
    \begin{equation*}
        0<2I_{\delta,c,\lambda}=2J_{\delta,c}(\psi_{\delta,c,\lambda}) = \theta_{\delta,c,\lambda}(p+2)\lambda \quad \text{and} \quad 0<2I_{\delta,c,\lambda}^{(s)}=2J_{\delta,c}^{(s)}(\psi_{\delta,c,\lambda}^{(s)}) = \theta_{\delta,c,\lambda}^{(s)}(p+2)\lambda.
    \end{equation*}
    Doing the same with the Euler-Lagrange equations of the limit problems $I_{\infty,c,\lambda}$ and $I_{0,c,\lambda}$, and their solutions $\psi_{\infty,c,\lambda}$ and $\psi_{0,c,\lambda}$, we get
    \begin{equation*}
        0<2I_{\infty,c,\lambda}=2J_{\infty,c}(\psi_{\infty,c,\lambda}) = \theta_{\infty,c,\lambda}(p+2)\lambda
    \end{equation*}
    and
    \begin{equation*}
        0<2I_{0,c,\lambda}=2J_{0,c}(\psi_{0,c,\lambda}) = \theta_{0,c,\lambda}(p+2)\lambda.
    \end{equation*}
    Note that, recalling the decomposition $m_\delta^{(s)}(\xi)=4\pi^2\xi^2-\xi^2h(\xi,\delta)$ from Lemma 8.2.1 in \cite{Abdelouhab1989}, for any $\psi\in H^1$ with $\int\psi^{p+2}=\lambda$:
    \begin{align*}
        2J_{\delta,c}^{(s)}(\psi) &= \int \{4\pi^2\xi^2-\xi^2h(\xi,\delta)\}|\fourier{\psi}(\xi)|^2 + c\int |\fourier{\psi}(\xi)|^2 \\
        &\leq \int 4\pi^2\xi^2|\fourier{\psi}(\xi)|^2 + c\int |\fourier{\psi}(\xi)|^2 \\
        &= 2J_{0,c}(\psi).
    \end{align*}
    Thus, taking infimum over the set $\{\psi\in H^1: |\psi|_{p+2}^{p+2}=\lambda\}$ we get
    \begin{equation*}
        I_{\delta,c,\lambda}^{(s)} \leq I_{0,c,\lambda}, \quad \forall\delta>0.
    \end{equation*}
    Similarly, we obtain that
    \begin{equation*}
        I_{\delta,c,\lambda} \leq I_{\infty,c,\lambda}, \quad \forall\delta>0.
    \end{equation*}
    Now, writing these inequalities in terms of the Lagrange multipliers $\theta$ gives
    \begin{equation*}
        \theta_{\delta,c,\lambda}^{(s)} \leq \theta_{0,c,\lambda} \quad \text{and} \quad \theta_{\delta,c,\lambda} \leq \theta_{\infty,c,\lambda},
    \end{equation*}
    so setting $C_1:= \theta_{\infty,c,\lambda}$ and $C_2:= \theta_{0,c,\lambda}$, we get the bounds.

    For the final part, note that each term in the series defining $h(\xi,\delta)$ is non-decreasing in $\delta$, so that, for every $\xi\in\R$, $m_\delta^{(s)}(\xi)=4\pi^2\xi^2-\xi^2h(\xi,\delta)$ is non-increasing in $\delta$. Hence $J_{\delta',c}^{(s)}(\psi)\geq J_{\delta,c}^{(s)}(\psi)$ for all $\psi$ whenever $0<\delta'\leq\delta$, and taking infimum we get $I_{\delta',c,\lambda}^{(s)}\geq I_{\delta,c,\lambda}^{(s)}$, that is, $\theta_{\delta',c,\lambda}^{(s)}\geq\theta_{\delta,c,\lambda}^{(s)}$.
\end{proof}

\subsection{Convergence in the energy space}

We begin with the deep water limit in the energy space $H^{1/2}$, which will be used below to obtain the uniform convergence in that limit; the underlying idea follows the approach of \cite{Abdallah2024}. Due to the singular character of the shallow water limit $\delta\to0$, we have not been able to establish convergence in the energy space for the rescaled solitary waves; for them, we only prove the uniform bounds of Proposition \ref{Prop: psi shallow are H 1/2 uniformly bounded} below, which are enough to obtain convergence in $L^q$ in Section \ref{Section: Uniform convergence of solitary waves}.

\begin{lemma}\label{Lemma: Sufficient condition for convergence in the energy space}
    Let $\psi_{\delta,c,\lambda}$ be suitable solutions for $I_{\delta,c,\lambda}$. Then, $I_{\delta,c,\lambda}\to I_{\infty,c,\lambda}$ and $J_{\infty,c}(\psi_{\delta,c,\lambda})\to I_{\infty,c,\lambda}$ as $\delta\to\infty$. In particular, $\{\psi_{\delta,c,\lambda}\}_\delta$ is a minimizing sequence for $I_{\infty,c,\lambda}$.
\end{lemma}
\begin{proof}
    Let $\Theta_\delta:=\mathcal{H}\partial_x-M_\delta$, which, by Lemma \ref{Lemma: m delta and |xi|}, is a Fourier multiplier with symbol $2\pi|\xi|-m_\delta(\xi)\in[0,2/\delta]$. By Lemma \ref{Lemma: thetas are bounded}, $I_{\delta,c,\lambda}\leq I_{\infty,c,\lambda}$ for every $\delta>0$; since $m_\delta\geq0$, this gives $\frac{c}{2}|\psi_{\delta,c,\lambda}|_2^2\leq J_{\delta,c}(\psi_{\delta,c,\lambda})=I_{\delta,c,\lambda}\leq I_{\infty,c,\lambda}$. Hence, since $\psi_{\delta,c,\lambda}$ is admissible for $I_{\infty,c,\lambda}$,
    \begin{equation*}
        I_{\infty,c,\lambda} \leq J_{\infty,c}(\psi_{\delta,c,\lambda}) = J_{\delta,c}(\psi_{\delta,c,\lambda}) + \frac{1}{2}\int\psi_{\delta,c,\lambda}\Theta_\delta\psi_{\delta,c,\lambda} \leq I_{\delta,c,\lambda} + \frac{1}{\delta}|\psi_{\delta,c,\lambda}|_2^2 \leq I_{\infty,c,\lambda} + \frac{2I_{\infty,c,\lambda}}{c\,\delta},
    \end{equation*}
    and the result follows by letting $\delta\to\infty$.
\end{proof}

\begin{theorem}\label{Th: Deep water convergence of solitary waves}
    Let $c>0$ and $g$ be one of the admissible nonlinearities for equation \eqref{Eq: Broader class of GILW}: $g_I,g_N$ or $g_L$. For any $\delta\gg1$, let $\phi_{\delta,c}$ be a suitable solitary wave for equation \eqref{Eq: Broader class of GILW} with velocity $c$ and depth parameter $\delta$. Then, $\phi_{\delta,c}\to\phi_{\infty,c}$ in $H^{1/2}$ as $\delta\to\infty$, where $\phi_{\infty,c}$ is the suitable solitary wave of gBO for the same velocity $c$ and non-linearity $g$.
\end{theorem}
\begin{proof}
    Since $\phi_{\delta,c}$ is obtained from Theorem \ref{Th: Existence of suitable solitary waves}, we have $\phi_{\delta,c}=\eta_{\delta,c,\lambda}\psi_{\delta,c,\lambda}$, where $\psi_{\delta,c,\lambda}$ is a suitable solution of $I_{\delta,c,\lambda}$, $N_{\delta,c,\lambda}$ or $L_{\delta,c,\lambda}$; since it is suitable, it is a solution of $I_{\delta,c,\lambda}$ independently of the nonlinearity $g$ (Proposition \ref{Proposition: Existence of minimizer of I and L}). Similarly, $\phi_{\infty,c}=\eta_{\infty,c,\lambda}\psi_{\infty,c,\lambda}$, where $\psi_{\infty,c,\lambda}$ is a suitable solution of $I_{\infty,c,\lambda}$ and $\eta_{\infty,c,\lambda}^p=\theta_{\infty,c,\lambda}(p+2)(p+1)$.

    \emph{Step 1: convergence of the multipliers.} Multiplying the Euler-Lagrange equations of $I_{\delta,c,\lambda}$ and $I_{\infty,c,\lambda}$ by their respective solutions and integrating, we get $2I_{\delta,c,\lambda}=\theta_{\delta,c,\lambda}(p+2)\lambda$ and $2I_{\infty,c,\lambda}=\theta_{\infty,c,\lambda}(p+2)\lambda$. In particular, $\eta_{\infty,c,\lambda}$ is determined by $I_{\infty,c,\lambda}$, so that, by the uniqueness of $\phi_{\infty,c}$ (Proposition \ref{Prop: gBO waves are FL ground states}), $\psi_{\infty,c,\lambda}$ is the unique suitable solution of $I_{\infty,c,\lambda}$. Moreover, Lemma \ref{Lemma: Sufficient condition for convergence in the energy space} gives $\theta_{\delta,c,\lambda}\to\theta_{\infty,c,\lambda}$, and therefore $\eta_{\delta,c,\lambda}\to\eta_{\infty,c,\lambda}$.

    \emph{Step 2: convergence of the minimizers.} It suffices to show that any sequence $\delta_n\to\infty$ has a sub-sequence along which $\psi_{\delta_n,c,\lambda}\to\psi_{\infty,c,\lambda}$ in $H^{1/2}$. By Lemma \ref{Lemma: Sufficient condition for convergence in the energy space}, $\{\psi_{\delta_n,c,\lambda}\}_n$ is a minimizing sequence for $I_{\infty,c,\lambda}$ and, since these functions are nonnegative, also for the corresponding problem $N_{\infty,c,\lambda}$, whose value coincides with $I_{\infty,c,\lambda}$ (as in Proposition \ref{Proposition: Existence of minimizer of I and L}). The concentration-compactness argument of Section \ref{Section: Existence of suitable solitary waves for GILW} applies verbatim with $M_\delta$ replaced by $\mathcal{H}\partial_x$ (see also Theorem 8.3 in \cite{Angulo2009}), and yields a sub-sequence, still denoted $\psi_{\delta_n,c,\lambda}$, points $y_n\in\R$ and a minimizer $\tilde\psi$ such that $\psi_{\delta_n,c,\lambda}(\cdot+y_n)\to\tilde\psi$ in $H^{1/2}$.

    The sequence $y_n$ is bounded. Indeed, choose $R>0$ such that $\int_{-R}^{R}|\tilde\psi|^{p+2}>\lambda/2$. Since $\psi_{\delta_n,c,\lambda}$ is even and non-increasing in $|x|$, we have $\int_{-R}^{R}\psi_{\delta_n,c,\lambda}^{p+2}\geq\int_{y_n-R}^{y_n+R}\psi_{\delta_n,c,\lambda}^{p+2}$, and the last integral converges to $\int_{-R}^{R}|\tilde\psi|^{p+2}>\lambda/2$. If $|y_n|\geq2R$ along a sub-sequence, the intervals $(-R,R)$ and $(y_n-R,y_n+R)$ would be disjoint, and we would get $\int\psi_{\delta_n,c,\lambda}^{p+2}>\lambda$ for $n$ large, a contradiction. Hence, passing to a further sub-sequence, $y_n\to y\in\R$, and therefore $\psi_{\delta_n,c,\lambda}\to\tilde\psi(\cdot-y)$ in $H^{1/2}$. The limit $\tilde\psi(\cdot-y)$ is a minimizer of $I_{\infty,c,\lambda}$ and, being the point-wise a.e. limit of a sub-sequence of suitable functions, it is even, nonnegative and non-increasing in $|x|$; by the uniqueness established in Step 1, $\tilde\psi(\cdot-y)=\psi_{\infty,c,\lambda}$.

    \emph{Step 3: conclusion.} By Steps 1 and 2, $\phi_{\delta,c}=\eta_{\delta,c,\lambda}\psi_{\delta,c,\lambda}\to\eta_{\infty,c,\lambda}\psi_{\infty,c,\lambda}=\phi_{\infty,c}$ in $H^{1/2}$.
\end{proof}

For the rescaled solitary waves we cannot argue as above, but Lemma \ref{Lemma: thetas are bounded} still provides uniform bounds for the minimizers $\psi_{\delta,c,\lambda}^{(s)}$, which will be used later to prove uniform convergence in the shallow water limit.

\begin{proposition}\label{Prop: psi shallow are H 1/2 uniformly bounded}
    Let $\psi_{\delta,c,\lambda}^{(s)}$ be the rescaled suitable minimizers. Then,
    \begin{enumerate}
        \item $|\psi_{\delta,c,\lambda}^{(s)}|_2\lesssim 1$ for every $\delta>0$.
        \item For the constant $\alpha_1$ given in Lemma \ref{Lemma: Lower estimate for m_s_delta}:
            \begin{equation*}
                \int_{\{|\xi|\leq \alpha_1/\delta\}}|\xi|^2|\fourier{\psi_{\delta,c,\lambda}^{(s)}}|^2 \lesssim 1, \quad \text{for every } \delta>0.
            \end{equation*}
        \item $\|\psi_{\delta,c,\lambda}^{(s)}\|_{1/2}\lesssim 1$ for every $\delta\in(0,1]$.
    \end{enumerate}
\end{proposition}
\begin{proof}
    By Lemma \ref{Lemma: thetas are bounded}, for every $\delta>0$ we have $J_{\delta,c}^{(s)}(\psi_{\delta,c,\lambda}^{(s)})=I_{\delta,c,\lambda}^{(s)}\leq I_{0,c,\lambda}$, that is,
    \begin{equation*}
        \frac{1}{2}\int m_{\delta}^{(s)}(\xi)|\fourier{\psi_{\delta,c,\lambda}^{(s)}}|^2 + \frac{c}{2}\int |\fourier{\psi_{\delta,c,\lambda}^{(s)}}|^2 \leq I_{0,c,\lambda}.
    \end{equation*}
    Now, using the lower estimate for the symbol $m_{\delta}^{(s)}=(3/\delta)m_\delta$ given in Lemma \ref{Lemma: Lower estimate for m_s_delta}, we obtain, for every $\delta>0$,
    \begin{align}\label{Eq: Estimate in corollary of Prop sufficient conditions for H1 convergence}
        \frac{5\pi}{2}\int_{\{|\xi|\leq \alpha_1/\delta\}}|\xi|^2|\fourier{\psi_{\delta,c,\lambda}^{(s)}}|^2 &+ \frac{5\pi\alpha_1^2}{2\delta^2}\int_{\{\alpha_1/\delta\leq |\xi|\leq \alpha_2/\delta\}}|\fourier{\psi_{\delta,c,\lambda}^{(s)}}|^2 \notag\\
        &+ \frac{5\pi}{2\delta}\int_{\{\alpha_2/\delta\leq|\xi|\}}|\xi||\fourier{\psi_{\delta,c,\lambda}^{(s)}}|^2 + \frac{c}{2}\int|\fourier{\psi_{\delta,c,\lambda}^{(s)}}|^2 \leq I_{0,c,\lambda}.
    \end{align}
    Since all the terms on the left-hand side are nonnegative, the first two items follow.

    To prove the last item we estimate the $L^2$-norm of $|\xi|^{1/2}|\fourier{\psi_{\delta,c,\lambda}^{(s)}}|$:
    \begin{align*}
        \int_\R |\xi||\fourier{\psi_{\delta,c,\lambda}^{(s)}}|^2 &\leq \int_{\{|\xi|\leq 1\}}|\fourier{\psi_{\delta,c,\lambda}^{(s)}}|^2 + \int_{\{1\leq|\xi|\leq \alpha_1/\delta\}}|\xi|^2|\fourier{\psi_{\delta,c,\lambda}^{(s)}}|^2 + \frac{\alpha_2}{\delta}\int_{\{\alpha_1/\delta\leq |\xi|\leq \alpha_2/\delta\}}|\fourier{\psi_{\delta,c,\lambda}^{(s)}}|^2 \\
        &\quad\quad\quad +  \int_{\{\alpha_2/\delta\leq |\xi|\}}|\xi||\fourier{\psi_{\delta,c,\lambda}^{(s)}}|^2.
    \end{align*}
    The first two terms on the right-hand side are bounded by items 1 and 2, while, by \eqref{Eq: Estimate in corollary of Prop sufficient conditions for H1 convergence}, the last two are bounded by $C\delta$ for some constant $C$ independent of $\delta$. Therefore, $\|\psi_{\delta,c,\lambda}^{(s)}\|_{1/2}\lesssim 1$ for every $\delta\in(0,1]$.
\end{proof}

\subsection{Uniform convergence of solitary waves}\label{Section: Uniform convergence of solitary waves}

We now turn to uniform convergence, valid in both the deep and shallow water limits. The strategy is to first obtain estimates for the solitary waves $\phi_{\delta,c}^{(j)}$, $j\in\{s,d\}$, that are uniform in $\delta$, and for the convolution kernel $\mathcal{K}_{\delta,c}^{(j)}$ through which $\phi_{\delta,c}^{(j)}$ is expressed, and then combine both to pass to the limit in the convolution representation of the solitary wave equation.

We begin with two naive estimates.

\begin{lemma}[\textbf{Naive uniform estimates}]\label{Lemma: Naive uniform estimates on phi}
    Let $\phi_{\delta,c}$ be the suitable solitary waves of gILW and $\phi_{\delta,c}^{(s)}$ their rescalings. For any $q\in [2,p+2]$ and $j\in\{s,d\}$, there is some constant $C=C(c,q,j)>0$ such that
    \begin{equation*}
        |\phi_{\delta,c}^{(j)}|_q \leq C, \quad \forall \delta>0.
    \end{equation*}
\end{lemma}
\begin{proof}
    Let $j\in\{s,d\}$. On the one hand, for any $\delta>0$ we have
    \begin{equation*}
        |\phi_{\delta,c}^{(j)}|_{p+2}^{p+2}=(\eta_{\delta,c,\lambda}^{(j)})^{p+2}\lambda \leq C\lambda.
    \end{equation*}
    On the other hand, multiplying the solitary wave equation $M_\delta^{(j)}\phi_{\delta,c}^{(j)}+c\phi_{\delta,c}^{(j)}=\frac{1}{p+1}(\phi_{\delta,c}^{(j)})^{p+1}$ by $\phi_{\delta,c}^{(j)}$ and integrating:
    \begin{align*}
        c|\phi_{\delta,c}^{(j)}|_2^2 &= \frac{(\eta_{\delta,c,\lambda}^{(j)})^{p+2}}{p+1}\lambda - \int m_{\delta}^{(j)}|\fourier{\phi_{\delta,c}^{(j)}}|^2 \leq \frac{C\lambda}{p+1}, \\
        |\phi_{\delta,c}^{(j)}|_2 &\leq \left(\frac{C\lambda}{(p+1)c}\right)^{1/2}.
    \end{align*}
    Finally, interpolation gives the desired result.
\end{proof}

\begin{lemma}[\textbf{Naive decayment estimate}]\label{Lemma: uniform point-wise upper bound for solitary waves}
    Let $\phi_{\delta,c}$ be the suitable solitary waves of gILW and $\phi_{\delta,c}^{(s)}$ their rescalings. Then,
    \begin{equation*}
        \phi_{\delta,c}^{(j)}(x) \leq C_{c,\lambda,j}\frac{1}{|x|^{1/(p+2)}}, \quad \forall\delta>0.
    \end{equation*}
\end{lemma}
\begin{proof}
    Since $\phi_{\delta,c}^{(j)}=\eta_{\delta,c,\lambda}^{(j)}\psi_{\delta,c,\lambda}^{(j)}$ and the family $\eta_{\delta,c,\lambda}^{(j)}$ is uniformly bounded with respect to $\delta$, we only need to show that
    \begin{equation*}
        \psi_{\delta,c,\lambda}^{(j)}(x) \lesssim_{c,\lambda,j} \frac{1}{|x|^{1/(p+2)}}, \quad \forall\delta>0.
    \end{equation*}
    First, we prove that for all $\delta>0$:
    \begin{equation*}
        \psi_{\delta,c,\lambda}^{(j)}(x)\leq \left(\frac{\lambda}{2R}\right)^{1/(p+2)}, \quad \forall |x|\geq R.
    \end{equation*}
    Suppose that there exist some $\tilde\delta>0$ such that, for some $\tilde x,\tilde R \in \R$ satisfying $|\tilde x|\geq \tilde R$ we have
    \begin{equation*}
        |\psi_{\tilde \delta,c,\lambda}^{(j)}(\tilde x)| > \left(\frac{\lambda}{2\tilde R}\right)^{1/(p+2)}.
    \end{equation*}
    Since $\psi_{\tilde\delta,c,\lambda}^{(j)}$ is even, strictly positive and radially decreasing, the above inequality holds for any $x\in B_{\tilde R}$; then
    \begin{equation*}
        \lambda = \int (\psi_{\tilde\delta,c,\lambda}^{(j)})^{p+2} > \int_{-\tilde R}^{+\tilde R} (\psi_{\tilde\delta,c,\lambda}^{(j)})^{p+2} > \lambda,
    \end{equation*}
    a contradiction.

    Now, for any $x\in\R$, $x\neq0$, and $\epsilon>0$, we have:
    \begin{equation*}
        \psi_{\delta,c,\lambda}^{(j)}(x+\epsilon) \leq \left(\frac{\lambda}{2}\right)^{1/(p+2)}\frac{1}{|x|^{1/(p+2)}}
    \end{equation*}
    and letting $\epsilon\to0$ we get the desired inequality.
\end{proof}

Writing the solitary wave equation $M_\delta^{(j)}\phi_{\delta,c}^{(j)}+c\phi_{\delta,c}^{(j)}=\frac{1}{p+1}(\phi_{\delta,c}^{(j)})^{p+1}$ in Fourier variables and solving for $\fourier{\phi_{\delta,c}^{(j)}}$, we obtain
\begin{equation}
    \fourier{\phi_{\delta,c}^{(j)}}(\xi) = \left(\frac{1}{p+1}\right)\frac{1}{c+m_\delta^{(j)}(\xi)}[(\phi_{\delta,c}^{(j)})^{p+1}]^{\wedge}(\xi)
\end{equation}
Thus, defining the convolution kernel $\mathcal{K}_{\delta,c}^{(j)}$ as the function satisfying
\begin{equation*}
    \fourier{\mathcal{K}_{\delta,c}^{(j)}}(\xi):=\frac{1}{c+m_\delta^{(j)}(\xi)},
\end{equation*}
the solitary wave equation can equivalently be written in convolution form as
\begin{equation*}
    \phi_{\delta,c}^{(j)} = \frac{1}{p+1}\mathcal{K}_{\delta,c}^{(j)}*(\phi_{\delta,c}^{(j)})^{p+1}.
\end{equation*}
The next two lemmas establish uniform (in $\delta$) estimates for $\mathcal{K}_{\delta,c}^{(j)}$, which will then be transferred to the solitary waves $\phi_{\delta,c}^{(j)}$ themselves via this convolution representation.

\begin{lemma}\label{Lemma: Uniform estimates of fourier transform of convolution kernel K_c_delta}
    Let $c>0$, $j\in\{s,d\}$ and $r\in(1,\infty]$. Then, there is some constant $C=C(j,r,c)$ such that
    \begin{equation*}
        |\fourier{\mathcal{K}_{\delta,c}^{(j)}}|_r \leq C,
    \end{equation*}
    for every $\delta\geq 4/c$ if $j=d$, and for every $\delta\in(0,1]$ if $j=s$.
\end{lemma}
\begin{proof}
    We begin with the case $j=d$. It follows from Lemma \ref{Lemma: m delta and |xi|} that, for $\delta\geq4/c$:
    \begin{equation*}
        \fourier{\mathcal{K}_{\delta,c}}(\xi) = \frac{1}{c+m_\delta(\xi)} \leq \frac{1}{c+2\pi|\xi|-\frac{2}{\delta}} \leq \frac{1}{\frac{c}{2}+2\pi|\xi|}.
    \end{equation*}
    Then, $|\fourier{\mathcal{K}_{\delta,c}}|_r \leq C$ for any $r\in(1,\infty]$ and $\delta\geq 4/c$.

    For the case $j=s$, the case $r=\infty$ follows from $m_\delta^{(s)}\geq0$, which gives $\fourier{\mathcal{K}_{\delta,c}^{(s)}}\leq1/c$. For $r\in(1,\infty)$, Lemma \ref{Lemma: Lower estimate for m_s_delta} gives
    \begin{align*}
        \int_\R \left|\frac{1}{c+m_{\delta}^{(s)}(\xi)}\right|^{r} &\leq \int_{\R}\frac{d\xi}{(c+5\pi\xi^2)^r} + \int_{\{\alpha_1/\delta\leq|\xi|\leq \alpha_2/\delta\}}\left(\frac{\delta^2}{\delta^2c+5\pi\alpha_1^2}\right)^r \\ 
        & \quad \quad + \int_{\{\alpha_2/\delta\leq|\xi|\}} \left(\frac{\delta}{\delta c+5\pi|\xi|}\right)^r \\
        &=: I+II+III.
    \end{align*}
    Since $r>1/2$, $I$ is a finite constant depending only on $c$ and $r$. Moreover,
    \begin{equation*}
        II \leq \frac{2\alpha_2}{\delta}\left(\frac{\delta^2}{5\pi\alpha_1^2}\right)^r = \frac{2\alpha_2\,\delta^{2r-1}}{(5\pi\alpha_1^2)^r} \quad\text{ and }\quad III \leq 2\left(\frac{\delta}{5\pi}\right)^r\int_{\alpha_2/\delta}^{\infty}\frac{d\xi}{\xi^r} = \frac{2\,\alpha_2^{1-r}\,\delta^{2r-1}}{(5\pi)^r(r-1)}.
    \end{equation*}
    Since $2r-1>0$, both $II$ and $III$ are bounded for $\delta\in(0,1]$, which proves the claim.
\end{proof}

\begin{lemma}\label{Lemma: Uniform estimates of convolution kernel K_c_delta}
    For any $q\in[2,\infty)$ and $j\in\{s,d\}$, there is some constant $C=C(c,q,j)>0$ such that:
    \begin{equation*}
        |\mathcal{K}_{\delta,c}^{(s)}|_q\leq C, \quad \forall \delta\ll1 \quad \text{ and } \quad |\mathcal{K}_{\delta,c}|_q\leq C, \quad \forall \delta\gg1.
    \end{equation*}
\end{lemma}
\begin{proof}
    Lemma \ref{Lemma: Uniform estimates of fourier transform of convolution kernel K_c_delta} implies that $\mathcal{K}_{\delta,c}^{(j)}$ are functions, so its $L^q$-norm, even if infinite, is well-defined. 
    
    Note that, since $\fourier{\mathcal{K}_{\delta,c}^{(j)}}$ is even, then:
    \begin{equation*}
        \mathcal{K}_{\delta,c}^{(j)}(x) = [\fourier{\mathcal{K}_{\delta,c}^{(j)}}(\xi)]^{\vee}(x) = [\fourier{\mathcal{K}_{\delta,c}^{(j)}}(-\xi)]^{\wedge}(x) = [\fourier{\mathcal{K}_{\delta,c}^{(j)}}(\xi)]^{\wedge}(x).
    \end{equation*}
    Thus, Hausdorff-Young's inequality (Theorem 2.3 in \cite{Linares2015}) implies, for any $q\in [2,\infty)$:
    \begin{equation*}
        |\mathcal{K}_{\delta,c}^{(j)}|_q \leq |\fourier{\mathcal{K}_{\delta,c}^{(j)}}|_{r},
    \end{equation*}
    where $1/q+1/r=1$, so $r\in(1,2]$ and therefore the result follows from Lemma \ref{Lemma: Uniform estimates of fourier transform of convolution kernel K_c_delta}.
\end{proof}

\begin{proposition}\label{Prop: Uniform estimates of solitary waves}
    Let $\phi_{\delta,c}$ be the suitable solitary waves of gILW and $\phi_{\delta,c}^{(s)}$ their rescalings. Then, for any $r\in[2,\infty]$ and $j\in\{s,d\}$, there is some $C=C(r,c,j)>0$ such that
    \begin{equation*}
        |\phi_{\delta,c}^{(j)}|_r \leq C,
    \end{equation*}
    for $\delta\ll1$ if $j=s$, and $\delta\gg1$ if $j=d$.
\end{proposition}
\begin{proof}
    By interpolation, we only need to prove the result for $r=2$ and $r=\infty$. We prove both cases using Young's inequality. 
    
    On the one hand, for the case $r=2$ we have
    \begin{equation*}
        |\phi_{\delta,c}^{(j)}|_2 = |\mathcal{K}_{\delta,c}^{(j)}*(\phi_{\delta,c}^{(j)})^{p+1}|_2 \leq |\mathcal{K}_{\delta,c}^{(j)}|_2|(\phi_{\delta,c}^{(j)})^{p+1}|_1 = |\mathcal{K}_{\delta,c}^{(j)}|_2|\phi_{\delta,c}^{(j)}|_{p+1}^{p+1} \leq C_{c,2,j},
    \end{equation*}
    where last inequality follows from Lemma \ref{Lemma: Uniform estimates of convolution kernel K_c_delta} for $\delta\ll1$ if $j=s$, and $\delta\gg1$ if $j=d$.
    
    On the other hand, for the case $r=\infty$:
    \begin{equation*}
        |\phi_{\delta,c}^{(j)}|_\infty = |\mathcal{K}_{\delta,c}*(\phi_{\delta,c}^{(j)})^{p+1}|_\infty \leq |\mathcal{K}_{\delta,c}|_{q'}|(\phi_{\delta,c}^{(j)})^{p+1}|_{q} = |\mathcal{K}_{\delta,c}|_{q'}|\phi_{\delta,c}^{(j)}|_{q(p+1)}^{p+1} \leq C_{c,q,p,j},
    \end{equation*}
    where $1/q+1/q' =1$ and we have chosen $1<q<(p+2)/(p+1)<2$ in order to get the $\delta$-uniform estimates: $|\phi_{\delta,c}^{(j)}|_{q(p+1)}\lesssim_{c,q,p,j} 1$ from Lemma \ref{Lemma: Naive uniform estimates on phi} and $|\mathcal{K}_{\delta,c}^{(j)}|_{q'}\lesssim 1$ from Lemma \ref{Lemma: Uniform estimates of convolution kernel K_c_delta} for $\delta\ll1$ or $\delta\gg1$, depending on the choice of $j$.
\end{proof}

We now turn to the proofs of the uniform convergence results.

The two uniform convergence results below share their last step, which we isolate in the following elementary lemma.

\begin{lemma}\label{Lemma: From pointwise to uniform convergence}
    Let $\{f_n\}_n$ and $f$ be even, nonnegative functions which are non-increasing in $|x|$, with $f$ continuous. Assume that $f_n\to f$ point-wise in $\R$ and that there is some $C>0$ such that $f_n(x),f(x)\leq C|x|^{-1/(p+2)}$ for all $x\neq0$ and all $n$. Then, $f_n\to f$ uniformly in $\R$.
\end{lemma}
\begin{proof}
    Let $\epsilon>0$ and choose $R>0$ such that $CR^{-1/(p+2)}<\epsilon$; then $|f_n(x)-f(x)|<\epsilon$ for all $|x|\geq R$ and all $n$. By evenness, it remains to consider $x\in[0,R]$. By the uniform continuity of $f$ on $[0,R]$, there is a partition $0=x_0<x_1<\cdots<x_N=R$ such that $f(x_{i})-f(x_{i+1})<\epsilon$ for all $i$, and, by the point-wise convergence at these finitely many points, there is $n_0$ such that $|f_n(x_i)-f(x_i)|<\epsilon$ for all $i$ and $n\geq n_0$. For $x\in[x_i,x_{i+1}]$ and $n\geq n_0$, the monotonicity of $f_n$ and $f$ gives
    \begin{gather*}
        f_n(x)\leq f_n(x_i) < f(x_i)+\epsilon < f(x)+2\epsilon, \\ f_n(x)\geq f_n(x_{i+1}) > f(x_{i+1})-\epsilon > f(x)-2\epsilon,
    \end{gather*}
    so that $|f_n(x)-f(x)|<2\epsilon$ for all $x\in[0,R]$ and $n\geq n_0$.
\end{proof}

\begin{theorem}\label{Th: Deep water uniform convergence of solitary waves}
    Let $c>0$ and $g$ be one of the admissible nonlinearities for equation \eqref{Eq: Broader class of GILW}. For any $\delta\gg1$, let $\phi_{\delta,c}$ be a suitable solitary wave for equation \eqref{Eq: Broader class of GILW} with velocity $c$ and depth parameter $\delta$. Then, $\phi_{\delta,c}\to\phi_{\infty,c}$ in $L^\infty$ as $\delta\to\infty$, where $\phi_{\infty,c}$ is a suitable solitary wave of gBO for the same velocity $c$ and non-linearity $g$.
\end{theorem}
\begin{proof}
    By Theorem \ref{Th: Deep water convergence of solitary waves}, there is some $\phi_{\infty,c}\in H^{1/2}$ such that $\phi_{\delta,c}\to \phi_{\infty,c}$ in $H^{1/2}$. We are going to show that this convergence is uniform.

    Note that it is enough to prove that any sub-sequence of $\phi_{\delta,c}$ has a sub-sequence which converge to $\phi_{\infty,c}$ in $L^\infty$. Therefore, we consider an arbitrary sub-sequence, which we still denote as $\phi_{\delta,c}$.
    
    Since the solitary waves $\phi_{\delta,c}$ are suitable and, by Proposition \ref{Prop: Uniform estimates of solitary waves}, uniformly bounded in $L^{\infty}$ for $\delta\gg1$, by Helly's selection theorem (Theorem \ref{Theorem: Helly selection theorem} in Appendix \ref{Appendix: Helly selection theorem}) we know there exists a point-wise convergent sub-sequence, which we still denote as $\phi_{\delta,c}$, that converge to some (radially) non-increasing function $\tilde\phi_{\infty,c}$; we show that $\tilde\phi_{\infty,c}=\phi_{\infty,c}$.

    In fact, since $\phi_{\delta,c}\to \phi_{\infty,c}$ in $H^{1/2}$, passing to a sub-sequence if necessary, we have that $\phi_{\delta,c}\to \phi_{\infty,c}$ point-wise a.e. in $\R$, and by uniqueness of the point-wise limit, $\tilde\phi_{\infty,c} = \phi_{\infty,c}$ a.e. in $\R$. Note that both, $\tilde\phi_{\infty,c}$ and $\phi_{\infty,c}$, are (radially) non-increasing, and $\phi_{\infty,c}$ is continuous, so necessarily $\tilde\phi_{\infty,c} = \phi_{\infty,c}$ in all $\R$. In particular, $\phi_{\delta,c}\to \phi_{\infty,c}$ point-wise in $\R$.

    Finally, since $\phi_{\infty,c}$ is continuous and, by Lemma \ref{Lemma: uniform point-wise upper bound for solitary waves}, the functions $\phi_{\delta,c}$ and $\phi_{\infty,c}$ satisfy a common decay bound, Lemma \ref{Lemma: From pointwise to uniform convergence} shows that the convergence is uniform in $\R$.
\end{proof}

\begin{theorem}\label{Th: Shallow water uniform convergence of solitary waves}
    Let $c>0$ and $g$ be one of the nonlinearities admissible for equation \eqref{Eq: Broader class of GILW}. For any $\delta\ll1$, let $\phi_{\delta,\delta c/3}$ be a suitable solitary wave for equation \eqref{Eq: Broader class of GILW} with velocity $\delta c/3$ and depth parameter $\delta$, and let $\phi_{\delta,c}^{(s)}=(3/\delta)^{1/p}\phi_{\delta,\delta c/3}$ be its rescaling. Then, $\phi_{\delta,c}^{(s)}\to\phi_{0,c}$ in $L^\infty$ as $\delta\to0$, where $\phi_{0,c}$ is the suitable solitary wave of gKdV for the same velocity $c$ and non-linearity $g$.
\end{theorem}
\begin{proof}
    As in the deep-water case, it is enough to prove that any sub-sequence of $\phi_{\delta,c}^{(s)}$ has a sub-sequence which converge to $\phi_{0,c}$ in $L^\infty$. Therefore, we consider an arbitrary sub-sequence, which we still denote as $\phi_{\delta,c}^{(s)}$.

    Again, since the solitary waves $\phi_{\delta,c}^{(s)}$ are suitable and, by Proposition \ref{Prop: Uniform estimates of solitary waves}, uniformly bounded in $L^{\infty}$ for $\delta\ll1$, by Helly's selection theorem (Theorem \ref{Theorem: Helly selection theorem} in Appendix \ref{Appendix: Helly selection theorem}) we know there exists a point-wise convergent sub-sequence, which we still denote as $\phi_{\delta,c}^{(s)}$, that converge to some (radially) non-increasing function $\tilde\phi_{0,c}$. We have to show that $\tilde\phi_{0,c}=\phi_{0,c}$, but, in contrast with the deep-water case, we do not have the convergence in the energy space; therefore, we work with the solitary wave equation in frequency.

    For each $\delta\ll1$, we have
    \begin{equation}\label{Eq: 1 in proof uniform convergence shallow water}
        m_{\delta}^{(s)}\fourier{\phi_{\delta,c}^{(s)}} + c\fourier{\phi_{\delta,c}^{(s)}} = \frac{1}{p+1}[(\phi_{\delta,c}^{(s)})^{p+1}]^{\wedge}.
    \end{equation}

    By Proposition \ref{Prop: Uniform estimates of solitary waves}, $\phi_{\delta,c}^{(s)}$ is uniformly bounded in $L^2$, so, passing to a sub-sequence if necessary, there exists some $\tilde{\tilde\phi}_{0,c}\in L^2$ such that $\phi_{\delta,c}^{(s)}\rightharpoonup \tilde{\tilde\phi}_{0,c}$ in $L^2$. This means that for any $\zeta\in C^{\infty}_0$:
    \begin{equation*}
        \int\tilde{\tilde\phi}_{0,c}\zeta = \lim_{\delta\to0}\int \phi_{\delta,c}^{(s)}\zeta = \int \tilde\phi_{0,c}\zeta,
    \end{equation*}
    where last equality follows by dominated convergence, since $\phi_{\delta,c}^{(s)}$ is $L^\infty$-uniformly bounded and converge point-wise to $\tilde\phi_{0,c}$. In particular, $\tilde{\tilde\phi}_{0,c}= \tilde\phi_{0,c}$ a.e. in $\R$.

    Similarly, since $\phi_{\delta,c}^{(s)}$ is also uniformly bounded in $L^{2(p+1)}$, the same argument of the paragraph above implies that there exists a sub-sequence such that $(\phi_{\delta,c}^{(s)})^{p+1}\rightharpoonup (\tilde\phi_{0,c})^{p+1}$ in $L^2$.

    These two $L^2$-weak convergence imply that for any $\zeta\in C^{\infty}_0$:
    \begin{equation*}
        \int \fourier{\phi_{\delta,c}^{(s)}} \zeta = \int \phi_{\delta,c}^{(s)} \fourier{\zeta} \to \int \tilde\phi_{0,c} \fourier{\zeta} = \int \fourier{\tilde\phi_{0,c}} \zeta,
    \end{equation*}
    and
    \begin{equation*}
        \int [(\phi_{\delta,c}^{(s)})^{p+1}]^{\wedge} \zeta \to \int [(\tilde\phi_{0,c})^{p+1}]^{\wedge} \zeta.
    \end{equation*}

    We will also need $\fourier{\phi_{\delta,c}^{(s)}}$ to be uniformly bounded in $L^\infty$. Solving for $\fourier{\phi_{\delta,c}^{(s)}}$ in \eqref{Eq: 1 in proof uniform convergence shallow water} gives $\fourier{\phi_{\delta,c}^{(s)}} = \frac{1}{p+1}\fourier{\mathcal{K}_{\delta,c}^{(s)}}\cdot[(\phi_{\delta,c}^{(s)})^{p+1}]^{\wedge}$, so that
    \begin{equation*}
        |\fourier{\phi_{\delta,c}^{(s)}}|_\infty \leq \frac{1}{p+1}|\fourier{\mathcal{K}_{\delta,c}^{(s)}}|_\infty\left|(\phi_{\delta,c}^{(s)})^{p+1}\right|_1 \leq \frac{1}{(p+1)c}|\phi_{\delta,c}^{(s)}|_{p+1}^{p+1} \lesssim 1,
    \end{equation*}
    where we used $\fourier{\mathcal{K}_{\delta,c}^{(s)}}(\xi)=1/(c+m_\delta^{(s)}(\xi))\leq1/c$, since $m_\delta^{(s)}\geq0$, and the uniform bound $|\phi_{\delta,c}^{(s)}|_{p+1}\lesssim1$ given by Lemma \ref{Lemma: Naive uniform estimates on phi}.

    Also, recalling the decomposition $m_\delta^{(s)}(\xi)=4\pi^2\xi^2-\xi^2h(\xi,\delta)$ from Lemma 8.2.1 in \cite{Abdelouhab1989}, for each $\zeta\in C^{\infty}_0$:
    \begin{equation*}
        \int m_{\delta}^{(s)}\fourier{\phi_{\delta,c}^{(s)}}\zeta = \int_{\text{Supp}(\zeta)}4\pi^2\xi^2\fourier{\phi_{\delta,c}^{(s)}}(\xi)\zeta(\xi)d\xi - \int_{\text{Supp}(\zeta)}\xi^2h(\xi,\delta)\fourier{\phi_{\delta,c}^{(s)}}(\xi)\zeta(\xi)d\xi,
    \end{equation*}
    where we take the limit $\delta\to0$ and use the compact support of $\zeta$, together with the uniform bound $|\fourier{\phi_{\delta,c}^{(s)}}|_\infty\lesssim1$ established above, to apply dominated convergence in the last term of the RHS; then we have
    \begin{equation*}
        \lim_{\delta\to0}\int m_{\delta}^{(s)}\fourier{\phi_{\delta,c}^{(s)}}\zeta = \lim_{\delta\to0}\int 4\pi^2\xi^2\fourier{\phi_{\delta,c}^{(s)}}(\xi)\zeta(\xi)d\xi = \int 4\pi^2\xi^2\fourier{\tilde\phi_{0,c}}(\xi)\zeta(\xi)d\xi,
    \end{equation*}
    where the last equality follows from the weak $L^2$-convergence $\fourier{\phi_{\delta,c}^{(s)}}\rightharpoonup\fourier{\tilde\phi_{0,c}}$ established above, tested against $4\pi^2\xi^2\zeta(\xi)\in C^\infty_0$.

    Now, multiplying equation \eqref{Eq: 1 in proof uniform convergence shallow water} by $\zeta\in C^{\infty}_0$ an integrating:
    \begin{equation*}
        \int m_{\delta}^{(s)}\fourier{\phi_{\delta,c}^{(s)}}\zeta + c\int\fourier{\phi_{\delta,c}^{(s)}}\zeta = \frac{1}{p+1}\int[(\phi_{\delta,c}^{(s)})^{p+1}]^{\wedge}\zeta,
    \end{equation*}
    so, taking the limit $\delta\to0$ we get
    \begin{equation*}
        \int \left(4\pi^2\xi^2\fourier{\tilde\phi_{0,c}} + c\fourier{\tilde\phi_{0,c}} - \frac{1}{p+1}[(\tilde\phi_{0,c})^{p+1}]^{\wedge}\right)\zeta=0, \quad \forall \zeta\in C^{\infty}_0,
    \end{equation*}
    therefore, the function $\tilde\phi_{0,c}\in L^2$ satisfies the equation
    \begin{equation}\label{Eq: 2 in proof uniform convergence shallow water}
        4\pi^2\xi^2\fourier{\tilde\phi_{0,c}} + c\fourier{\tilde\phi_{0,c}} = \frac{1}{p+1}[(\tilde\phi_{0,c})^{p+1}]^{\wedge}, \quad \text{ for a.e. } \xi\in\R.
    \end{equation}
    Note that $\fourier{\tilde\phi_{0,c}}$ and $[(\tilde\phi_{0,c})^{p+1}]^{\wedge}$ are in $L^2$, thus $4\pi^2\xi^2\fourier{\tilde\phi_{0,c}} \in L^2$ and it has an inverse Fourier transform in $L^2$, which is given by $-\partial_x^2\tilde\phi_{0,c}$. Now, taking inverse Fourier transform to equation \eqref{Eq: 2 in proof uniform convergence shallow water}:
    \begin{equation}\label{Eq: 3 in proof uniform convergence shallow water}
        -\partial_x^2\tilde\phi_{0,c} + c\tilde\phi_{0,c} = \frac{1}{p+1}(\tilde\phi_{0,c})^{p+1},
    \end{equation}
    in the $L^2$ sense. This means that $\tilde\phi_{0,c}$ is a suitable solitary wave for gKdV. Since there is only one non-trivial suitable solitary wave of gKdV, named $\phi_{0,c}$, we are only left to show that $\tilde\phi_{0,c}$ is non-trivial.

    The non-triviality of $\tilde\phi_{0,c}$ is obtained by the same technical tool we have used to rule-out vanishing when applying the concentration-compactness principle in Section \ref{Section: Existence of suitable solitary waves for GILW}. Recall that $\phi_{\delta,c}^{(s)}=\eta_{\delta,c,\lambda}^{(s)}\psi_{\delta,c,\lambda}^{(s)}$ with $|\psi_{\delta,c,\lambda}^{(s)}|_{p+2}^{p+2}=\lambda$ independent of $\delta$, while, by Lemma \ref{Lemma: thetas are bounded}, the map $\delta\mapsto\eta_{\delta,c,\lambda}^{(s)}>0$ is non-increasing on $(0,\infty)$. Hence, fixing any $\delta_0>0$, we get $\eta_{\delta,c,\lambda}^{(s)}\geq\eta_{\delta_0,c,\lambda}^{(s)}>0$ for all $\delta\leq\delta_0$, and therefore $|\phi_{\delta,c}^{(s)}|_{p+2}=\eta_{\delta,c,\lambda}^{(s)}\lambda^{1/(p+2)}\geq\eta_{\delta_0,c,\lambda}^{(s)}\lambda^{1/(p+2)}>0$ uniformly for $\delta\leq\delta_0$. Combining this with Proposition \ref{Prop: psi shallow are H 1/2 uniformly bounded}, it follows that, for $\delta\ll1$, there exists some $\epsilon>0$ and $B>0$ such that $\epsilon\leq|\phi_{\delta,c}^{(s)}|_{p+2}$ and $\|\phi_{\delta,c}^{(s)}\|_{1/2}\leq B$; therefore, by Lemma \ref{Lemma: to rule out vanishing}, there exists some $\eta>0$, independent of $\delta$, such that
    \begin{equation*}
        \sup_{y\in\R}\int_{y-2}^{y+2}|\phi_{\delta,c}^{(s)}(x)|^{p+2}dx\geq\eta, \quad \forall\delta\ll1.
    \end{equation*}
    Since the solitary waves $\phi_{\delta,c}^{(s)}$ are suitable, the $\sup$ is always attained at $y=0$, so we have the uniform inequality
    \begin{equation*}
        \int_{-2}^{2}|\phi_{\delta,c}^{(s)}(x)|^{p+2}dx\geq\eta,\quad\forall\delta\ll1.
    \end{equation*}
    Finally, taking the limit $\delta\to0$ in the inequality above and applying dominated convergence, it follows that $\tilde\phi_{0,c}$ is non-trivial, and therefore $\tilde\phi_{0,c}=\phi_{0,c}$. In particular, $\phi_{\delta,c}^{(s)}\to \phi_{0,c}$ point-wise in $\R$.

    Finally, as in the deep water case, Lemmas \ref{Lemma: uniform point-wise upper bound for solitary waves} and \ref{Lemma: From pointwise to uniform convergence} show that the convergence is uniform in $\R$.
\end{proof}

\begin{corollary}\label{Cor: Lq convergence shallow water}\label{Cor: L2 convergence shallow water}
    Let $c>0$ and $g$ be one of the admissible nonlinearities for equation \eqref{Eq: Broader class of GILW}. For any $\delta\ll1$, let $\phi_{\delta,c}^{(s)}=(3/\delta)^{1/p}\phi_{\delta,\delta c/3}$ be the rescaling of a suitable solitary wave $\phi_{\delta,\delta c/3}$ of equation \eqref{Eq: Broader class of GILW} with velocity $\delta c/3$ and depth parameter $\delta$. Then, for every $q\in[2,\infty]$, $\phi_{\delta,c}^{(s)}\to\phi_{0,c}$ in $L^q$ as $\delta\to0$, where $\phi_{0,c}$ is the suitable solitary wave of gKdV for the same velocity $c$ and non-linearity $g$.
\end{corollary}
\begin{proof}
    The case $q=\infty$ is Theorem \ref{Th: Shallow water uniform convergence of solitary waves}. For $q\in(2,\infty)$, note that $\phi_{\delta,c}^{(s)}-\phi_{0,c}$ is bounded in $L^2$ for $\delta\ll1$, by Proposition \ref{Prop: Uniform estimates of solitary waves} and since $\phi_{0,c}\in L^2$; the interpolation inequality $|f|_q\leq|f|_2^{2/q}|f|_\infty^{1-2/q}$ and the case $q=\infty$ then give the result.

    The case $q=2$ is not reachable by this interpolation argument, so we argue directly on the solitary wave equations.

    Writing \eqref{Eq: rescaled solitary wave equation}, and the corresponding equation for gKdV, in Fourier variables and solving for the transform of the solitary wave, we obtain
    \begin{equation*}
        \fourier{\phi_{\delta,c}^{(s)}}(\xi) = \frac{1}{p+1}\frac{[(\phi_{\delta,c}^{(s)})^{p+1}]^{\wedge}(\xi)}{c+m_\delta^{(s)}(\xi)} \quad\text{ and }\quad \fourier{\phi_{0,c}}(\xi) = \frac{1}{p+1}\frac{[\phi_{0,c}^{p+1}]^{\wedge}(\xi)}{c+4\pi^2\xi^2},
    \end{equation*}
    where we used that the symbol of $-\partial_x^2$ is $4\pi^2\xi^2$. Subtracting, we get the decomposition
    \begin{equation*}
        \fourier{\phi_{\delta,c}^{(s)}}-\fourier{\phi_{0,c}} = \underbrace{\frac{1}{p+1}\cdot\frac{[(\phi_{\delta,c}^{(s)})^{p+1}]^{\wedge}-[\phi_{0,c}^{p+1}]^{\wedge}}{c+m_\delta^{(s)}}}_{=:A_\delta} + \underbrace{\frac{1}{p+1}\left(\frac{1}{c+m_\delta^{(s)}}-\frac{1}{c+4\pi^2\xi^2}\right)[\phi_{0,c}^{p+1}]^{\wedge}}_{=:B_\delta},
    \end{equation*}
    so that, by Plancherel's theorem, $|\phi_{\delta,c}^{(s)}-\phi_{0,c}|_2 \leq |A_\delta|_2 + |B_\delta|_2$. We estimate both terms separately.

    \emph{The term $A_\delta$.} By Lemma \ref{Lemma: Lower estimate for m_s_delta} we have $m_\delta^{(s)}\geq0$, so that $(c+m_\delta^{(s)})^{-1}\leq c^{-1}$, and Plancherel's theorem gives
    \begin{equation*}
        |A_\delta|_2 \leq \frac{1}{(p+1)c}\left|(\phi_{\delta,c}^{(s)})^{p+1}-\phi_{0,c}^{p+1}\right|_2.
    \end{equation*}
    Now, by the mean value theorem, $|a^{p+1}-b^{p+1}|\leq (p+1)(a^p+b^p)|a-b|$ for any $a,b\geq0$, hence, by the Cauchy-Schwarz inequality,
    \begin{equation*}
        \left|(\phi_{\delta,c}^{(s)})^{p+1}-\phi_{0,c}^{p+1}\right|_2 \leq (p+1)\left(\left|\phi_{\delta,c}^{(s)}\right|_{4p}^{p}+\left|\phi_{0,c}\right|_{4p}^{p}\right)\left|\phi_{\delta,c}^{(s)}-\phi_{0,c}\right|_{4}.
    \end{equation*}
    The point of this splitting is that the exponent $4$ appearing in the last factor is strictly larger than $2$, so that no circularity arises. Indeed, $|\phi_{\delta,c}^{(s)}|_{4p}\lesssim_{c,p}1$ for $\delta\ll1$ by Proposition \ref{Prop: Uniform estimates of solitary waves} (note $4p\geq4\geq2$), $|\phi_{0,c}|_{4p}<\infty$, and $|\phi_{\delta,c}^{(s)}-\phi_{0,c}|_4\to0$ by the case $q=4$ above. Therefore $|A_\delta|_2\to0$ as $\delta\to0$.

    \emph{The term $B_\delta$.} Recall from Lemma 8.2.1 in \cite{Abdelouhab1989} that $m_\delta^{(s)}(\xi) = 4\pi^2\xi^2 - \xi^2h(\xi,\delta)$; dropping the term $4\delta^2\xi^2$ from the denominator in the definition of $h(\xi,\delta)$ gives
    \begin{equation*}
        0\leq \xi^2h(\xi,\delta) = \sum_{n=1}^\infty\frac{96\delta^2\xi^4}{n^2(n^2+4\delta^2\xi^2)} \leq 96\delta^2\xi^4\frac{\pi^4}{90}.
    \end{equation*}
    Hence, for each fixed $\xi\in\R$, $m_\delta^{(s)}(\xi)\to4\pi^2\xi^2$ as $\delta\to0$, and consequently the integrand of $|B_\delta|_2^2$ converges to $0$ point-wise. Moreover, since both $(c+m_\delta^{(s)})^{-1}$ and $(c+4\pi^2\xi^2)^{-1}$ are bounded by $c^{-1}$, this integrand is dominated by $\frac{4}{(p+1)^2c^2}\left|[\phi_{0,c}^{p+1}]^{\wedge}\right|^2$, which is integrable because $\phi_{0,c}^{p+1}\in L^2$. The dominated convergence theorem then gives $|B_\delta|_2\to0$ as $\delta\to0$, which finishes the proof.
\end{proof}

\section{Orbital stability and instability of solitary waves of gILW}\label{Section: Orbital instability of solitary waves in shallow water}

In this section we prove the main results of the paper: the orbital stability and instability of the solitary waves of gILW when the product $\delta c$ of the depth parameter and the velocity is small (slow waves, Section \ref{Section: Shallow water regime}) or large (fast waves, Section \ref{Section: Deep water regime}). Both stability and instability are obtained simultaneously from a single criterion, Theorem \ref{Theorem: Stability instability criterion Albert}, which requires three kinds of information: the structural conditions (C1)-(C3) on the symbol $m_\delta$, three spectral conditions on the solitary wave operator, and the sign of the scalar quantity $I=\langle\psi,\phi\rangle$, where $\psi\in L^2$ solves $\mathcal{L}\psi=\phi$.

The symbol and spectral conditions are established in Section \ref{Section: Spectral conditions}, together with their invariance under the natural scaling of gILW. The sign of $I$ is then obtained by transferring, as $\delta\to0$, the sign of the corresponding quantity $I_0$ for the gKdV operator $\mathcal{L}_{0,c}$, and, as $\delta\to\infty$, the sign of the quantity $I_\infty$ for the gBO operator $\mathcal{L}_{\infty,c}$; their values reproduce the classical thresholds $p=4$ for gKdV \cite{Bona1987} and $p=2$ for gBO. In the slow regime, the rescaled solitary waves $\phi_{\delta,c}^{(s)}$ of Section \ref{Section: Convergence of solitary waves} are used only as a technical device: they are the ones that converge to the solitary waves of gKdV, and the information obtained for them is transferred directly to the solitary waves of gILW.

\subsection{Spectral conditions, symbol conditions and scaling}\label{Section: Spectral conditions}

The perturbation theory tools for closed operators introduced in Section \ref{Section: Perturbation theory tools} (see, in particular, Propositions \ref{Prop: Basic properties of gap metric} and \ref{Prop: Perturbation theorem}) are used throughout. We consider the solitary wave operators of gILW and of the rescaled solitary waves of Section \ref{Section: Convergence of solitary waves},
\begin{equation*}
    \mathcal{L}_{\delta,c}:= M_{\delta} + c - \phi_{\delta,c}^p \quad\text{ and }\quad \mathcal{L}_{\delta,c}^{(s)}:= M_{\delta}^{(s)} + c - (\phi_{\delta,c}^{(s)})^p,
\end{equation*}
and, according to the convention of Section \ref{Section: Convergence of solitary waves}, we also write $\mathcal{L}_{\delta,c}^{(d)}:=\mathcal{L}_{\delta,c}$. Straightforward computations show that $\mathcal{L}_{\delta,c}^{(s)}=(3/\delta)\mathcal{L}_{\delta,\delta c/3}$; in particular, both operators have proportional spectra and the same eigenfunctions, so that each of the spectral conditions below holds for $\mathcal{L}_{\delta,c}^{(s)}$ if and only if it holds for $\mathcal{L}_{\delta,\delta c/3}$.

\begin{theorem}\label{Theorem: Convergence of solitary wave operators}
    The following limits are satisfied:
    \begin{enumerate}
        \item $d(\mathcal{L}_{\delta,c},\mathcal{L}_{\infty,c})\to0$ as $\delta\to\infty$, where $\mathcal{L}_{\infty,c} = \mathcal{H}\partial_x + c - \phi_{\infty,c}^{p}$,
        \item $d(\mathcal{L}_{\delta,c}^{(s)},\mathcal{L}_{0,c})\to0$ as $\delta\to 0$, where $\mathcal{L}_{0,c} = -\partial^2_x + c - \phi_{0,c}^{p}$.
    \end{enumerate}
\end{theorem}
\begin{proof}
    For the first limit, note that
    \begin{align*}
        d(\mathcal{L}_{\delta,c},\mathcal{L}_{\infty,c}) &= d(M_{\delta}+c-(\phi_{\delta,c})^p, \mathcal{H}\partial_x + c - \phi_{\infty,c}^{p}) \\
        &\leq 2(1+\|c-(\phi_{\delta,c})^p\|_{2,2}^2) d(M_{\delta}, \mathcal{H}\partial_x - \phi_{\infty,c}^{p}+(\phi_{\delta,c})^p)\\
        &\lesssim d(M_{\delta}, \mathcal{H}\partial_x) + d(\mathcal{H}\partial_x, \mathcal{H}\partial_x - \phi_{\infty,c}^{p}+(\phi_{\delta,c})^p) \\
        &\leq d(M_{\delta}, M_{\delta}+\mathcal{H}\partial_x-M_{\delta}) + \|\phi_{\infty,c}^{p}-(\phi_{\delta,c})^p\|_{2,2}\\
        &\leq \|\mathcal{H}\partial_x-M_{\delta}\|_{2,2}+\|\phi_{\infty,c}^{p}-(\phi_{\delta,c})^p\|_{2,2}
    \end{align*}

    then, since $\phi_{\delta,c}\to\phi_{\infty,c}$ uniformly and 
    \begin{equation*}
        \left|[\mathcal{H}\partial_x-M_{\delta}]^{\wedge}(\xi)\right| = \left|2\pi|\xi|- 2\pi\xi\coth(2\pi\delta\xi) + \frac{1}{\delta} \right| \lesssim \frac{1}{\delta}, 
    \end{equation*}
    we conclude the desired limit.

    For the second limit, similarly we get
    \begin{align*}
        d(\mathcal{L}_{\delta,c}^{(s)},\mathcal{L}_{0,c}) &= d(M_{\delta}^{(s)}+c-(\phi_{\delta,c}^{(s)})^p, -\partial_x^2 + c - \phi_{0,c}^{p}) \\
        &\leq 2(1+\|c-(\phi_{\delta,c}^{(s)})^p\|_{2,2}^2) d(M_{\delta}^{(s)}, -\partial_x^2 - \phi_{0,c}^{p}+(\phi_{\delta,c}^{(s)})^p)\\
        &\lesssim d(M_{\delta}^{(s)}, -\partial_x^2) + d(-\partial_x^2, -\partial_x^2 - \phi_{0,c}^{p}+(\phi_{\delta,c}^{(s)})^p) \\
        &\leq d(M_{\delta}^{(s)}+I-I, -\partial_x^2+I-I) + \|\phi_{0,c}^{p}-(\phi_{\delta,c}^{(s)})^p\|_{2,2}\\
        &\lesssim d(M_{\delta}^{(s)}+I, -\partial_x^2+I) + \|\phi_{0,c}^{p}-(\phi_{\delta,c}^{(s)})^p\|_{2,2}\\
        &= d((M_{\delta}^{(s)}+I)^{-1}, (-\partial_x^2+I)^{-1}) + \|\phi_{0,c}^{p}-(\phi_{\delta,c}^{(s)})^p\|_{2,2} \\
        &\leq \|(M_{\delta}^{(s)}+I)^{-1} - (-\partial_x^2+I)^{-1}\|_{2,2} + \|\phi_{0,c}^{p}-(\phi_{\delta,c}^{(s)})^p\|_{2,2},
    \end{align*}

    where last estimates follows from boundedness of $(M_{\delta}^{(s)}+I)^{-1}$ and $(-\partial_x^2+I)^{-1}$. 

    Since $\phi_{\delta,c}^{(s)}\to\phi_{0,c}$ uniformly, the second term in the RHS of estimate above goes to zero; while the first term can be estimated as
    \begin{align*}
        \|(M_{\delta}^{(s)}+I)^{-1} - (-\partial_x^2+I)^{-1}\|_{2,2} &\leq \left\|\frac{1}{m_{\delta}^{(s)}+1}-\frac{1}{4\pi^2\xi^2+1}\right\|_{\infty}\\
        &= \sup_{\xi\in\R} \left|\frac{4\pi^2\xi^2-m_{\delta}^{(s)}(\xi)}{(4\pi^2\xi^2+1)(m_{\delta}^{(s)}(\xi)+1)}\right| \\
        &= \sup_{\xi\in\R} \left|\frac{4\pi^2\xi^2-\frac{3}{\delta^2}(2\pi\delta\xi\coth(2\pi \delta\xi)-1)}{(4\pi^2\xi^2+1)(\frac{3}{\delta^2}(2\pi\delta\xi\coth(2\pi \delta\xi)-1)+1)}\right| \\
        [y:=2\pi\delta\xi, \beta(y)=(y\coth(y)-1)/y^2]&= \sup_{\xi\in\R} \left|\frac{4\pi^2\xi^2-\frac{3}{\delta^2}\beta(y)y^2}{(4\pi^2\xi^2+1)(\frac{3}{\delta^2}\beta(y)y^2+1)}\right| \\
        &= \sup_{\xi\in\R} \left|\frac{4\pi^2\xi^2-12\pi^2\xi^2\beta(y)}{(4\pi^2\xi^2+1)(12\pi^2\xi^2\beta(y)+1)}\right| \\
        &\lesssim \sup_{\xi\in\R} \left|\frac{\xi^2-3\xi^2\beta(y)}{\xi^2+1 + \xi^4\beta(y)+\xi^2\beta(y)}\right| \\
        &= \sup_{\xi\in\R} \left|\frac{1-3\beta(y)}{1+\frac{1}{\xi^2} + \xi^2\beta(y)+\beta(y)}\right| \\
        &\leq \sup_{\xi\in\R} \left|\frac{1-3\beta(y)}{1+ \xi^2\beta(y)}\right| \\
        &= \sup_{y\in\R} \left|\frac{1-3\beta(y)}{1+ \frac{y^2}{4\pi^2\delta^2}\beta(y)}\right|.
    \end{align*}

    Note that $\beta$ has the following properties:
    \begin{enumerate}
        \item $\lim_{|y|\to0}(3\beta(y)-1)=0$,
        \item $y^2\beta(y)$ is an increasing function of $y$ and
        \item there is some $\mu>0$ such that $0\leq\beta(y)<\mu$ for any $y\in\R$.
    \end{enumerate}

    Suppose that any $\epsilon>0$ is given, pick some $\eta>0$ such that $|3\beta(y)-1|<\epsilon$ for $|y|\leq\eta$. Now, for any $\delta$, on the one hand
    \begin{equation*}
        \sup_{|y|\leq\eta} \left|\frac{1-3\beta(y)}{1+ \frac{y^2}{4\pi^2\delta^2}\beta(y)}\right| <\epsilon,
    \end{equation*}
    while on the other hand, picking some $\delta_0\ll1$ such that $\frac{3\mu+1}{1+\frac{y^2}{4\pi^2\delta_0^2}\beta(y)}<\epsilon$ for any $|y|\geq\eta$,
    \begin{equation*}
        \sup_{|y|\geq\eta} \left|\frac{1-3\beta(y)}{1+ \frac{y^2}{4\pi^2\delta^2}\beta(y)}\right|<\epsilon, \quad\forall\delta\leq\delta_0.
    \end{equation*}

    These estimates imply the desired limit $d(\mathcal{L}_{\delta,c}^{(s)},\mathcal{L}_{0,c})\to 0$ as $\delta\to0$.
\end{proof}

\begin{lemma}[\textbf{Spectral conditions in shallow water}]\label{Lemma: Spectral conditions for shallow-gILW operator}
    Let $c>0$ and $p\in\N$ be fixed, and let $\mu_0$, $\gamma_0>0$ be as in Theorem \ref{Theorem: Spectrum of gKdV solitary wave operator}. Then there exists some $\epsilon=\epsilon(c,p)>0$ such that, for any $0<\delta<\epsilon$, the solitary wave operator $\mathcal{L}_{\delta,c}^{(s)}$ satisfies:
    \begin{enumerate}
        \item the eigenvalue $0$ of $\mathcal{L}_{\delta,c}^{(s)}$ is simple, and $\ker\mathcal{L}_{\delta,c}^{(s)}=\text{span}\{\partial_x\phi_{\delta,c}^{(s)}\}$;
        \item $\mathcal{L}_{\delta,c}^{(s)}$ has exactly one negative eigenvalue $\mu_\delta$, it is simple, and $\mu_\delta\to\mu_0$ as $\delta\to0$;
        \item $\gamma_\delta := \text{dist}\left(0,\text{spec}(\mathcal{L}_{\delta,c}^{(s)})\backslash\{0\}\right) \geq \dfrac{\gamma_0}{2}$.
    \end{enumerate}
    In particular, $\mathcal{L}_{\delta,c}^{(s)}$ satisfies conditions 1, 2, and 3 in Theorem \ref{Theorem: Stability instability criterion Albert}.
\end{lemma}
\begin{proof}
    Throughout the proof we use that $d(\mathcal{L}_{\delta,c}^{(s)},\mathcal{L}_{0,c})\to0$ as $\delta\to0$, by Theorem \ref{Theorem: Convergence of solitary wave operators}, together with the description of $\text{spec}(\mathcal{L}_{0,c})$ given by Theorem \ref{Theorem: Spectrum of gKdV solitary wave operator}. Note first that, differentiating the solitary wave equation \eqref{Eq: rescaled solitary wave equation} with respect to $x$,
    \begin{equation}\label{Eq: derivative of phi lies in the kernel}
        \mathcal{L}_{\delta,c}^{(s)}\left[\partial_x\phi_{\delta,c}^{(s)}\right] = M_\delta^{(s)}\partial_x\phi_{\delta,c}^{(s)} + c\,\partial_x\phi_{\delta,c}^{(s)} - (\phi_{\delta,c}^{(s)})^p\,\partial_x\phi_{\delta,c}^{(s)} = 0,
    \end{equation}
    and that $\partial_x\phi_{\delta,c}^{(s)}\in L^2\backslash\{0\}$, since $\phi_{\delta,c}^{(s)}\in H^\infty$ is not constant. Thus $0$ is always an eigenvalue of $\mathcal{L}_{\delta,c}^{(s)}$.

    \emph{Items 1 and 3.} Consider the open disk $\Omega:=\{z\in\C:|z|<\gamma_0/2\}$, whose boundary is a smooth contour. By the definition of $\gamma_0$ we have $\text{spec}(\mathcal{L}_{0,c})\cap\overline{\Omega}=\{0\}$, a single eigenvalue of multiplicity one by Theorem \ref{Theorem: Spectrum of gKdV solitary wave operator}. Hence Proposition \ref{Prop: Perturbation theorem} provides some $\epsilon_1=\epsilon_1(c,p)>0$ such that, for any $0<\delta<\epsilon_1$, the set $\text{spec}(\mathcal{L}_{\delta,c}^{(s)})\cap\overline{\Omega}$ consists of finitely many eigenvalues of finite multiplicity whose multiplicities add up to one. Since $0$ belongs to this set by \eqref{Eq: derivative of phi lies in the kernel}, it is the only point of $\text{spec}(\mathcal{L}_{\delta,c}^{(s)})$ in $\overline{\Omega}$ and its multiplicity is one. This proves items 1 and 3 simultaneously.

    \emph{Item 2.} Since $m_\delta^{(s)}\geq0$ by Lemma \ref{Lemma: Lower estimate for m_s_delta}, and since $|\phi_{\delta,c}^{(s)}|_\infty\leq C_\infty$ uniformly for $\delta\ll1$ by Proposition \ref{Prop: Uniform estimates of solitary waves} with $r=\infty$, we get, for any $u$ in the domain of $\mathcal{L}_{\delta,c}^{(s)}$,
    \begin{equation*}
        \left\langle\mathcal{L}_{\delta,c}^{(s)}u,u\right\rangle = \left\langle M_\delta^{(s)}u,u\right\rangle + c|u|_2^2 - \int_\R(\phi_{\delta,c}^{(s)})^pu^2 \geq \left(c-C_\infty^p\right)|u|_2^2,
    \end{equation*}
    and analogously for $\mathcal{L}_{0,c}$. Hence, by Theorem \ref{Theorem: Basic information about spectrum}, there is some $a>0$, depending only on $c$ and $p$, such that
    \begin{equation}\label{Eq: uniform lower bound for the spectra}
        \text{spec}(\mathcal{L}_{\delta,c}^{(s)})\subset(-a,\infty) \quad\text{ and }\quad \text{spec}(\mathcal{L}_{0,c})\subset(-a,\infty),
    \end{equation}
    for any $\delta$ small enough. This uniform bound is what prevents negative eigenvalues from escaping to $-\infty$ as $\delta\to0$, and is therefore essential in what follows.

    Let now $\Omega'\subset\C$ be the open disk having the segment $(-a,-\gamma_0/3)$ as a diameter, so that its boundary is a smooth contour meeting $\R$ exactly at the two points $-a$ and $-\gamma_0/3$. Neither of them belongs to $\text{spec}(\mathcal{L}_{0,c})$: the first one by \eqref{Eq: uniform lower bound for the spectra}, and the second one because $\text{spec}(\mathcal{L}_{0,c})\cap(-\gamma_0,\gamma_0)=\{0\}$ by the definition of $\gamma_0$. Since the spectrum of $\mathcal{L}_{0,c}$ is real, this gives $\text{spec}(\mathcal{L}_{0,c})\cap\partial\Omega'=\emptyset$; moreover, by Theorem \ref{Theorem: Spectrum of gKdV solitary wave operator} and since $|\mu_0|\geq\gamma_0>\gamma_0/3$,
    \begin{equation*}
        \text{spec}(\mathcal{L}_{0,c})\cap\Omega' = \text{spec}(\mathcal{L}_{0,c})\cap\left(-a,-\frac{\gamma_0}{3}\right) = \{\mu_0\},
    \end{equation*}
    a single eigenvalue of multiplicity one. Proposition \ref{Prop: Perturbation theorem} then provides some $\epsilon_2>0$ such that, for any $0<\delta<\epsilon_2$, the multiplicities of the eigenvalues in $\text{spec}(\mathcal{L}_{\delta,c}^{(s)})\cap\overline{\Omega'}$ add up to one.

    Set $\epsilon:=\min\{\epsilon_1,\epsilon_2\}$ and let $0<\delta<\epsilon$. If $\lambda$ is any negative point of $\text{spec}(\mathcal{L}_{\delta,c}^{(s)})$, then $\lambda\leq-\gamma_0/2<-\gamma_0/3$ by item 3, and $\lambda>-a$ by \eqref{Eq: uniform lower bound for the spectra}; hence $\lambda\in\Omega'$. By the previous paragraph there is exactly one such $\lambda$, which we denote by $\mu_\delta$, and it is a simple eigenvalue. Finally, applying the last assertion of Proposition \ref{Prop: Perturbation theorem} to an open disk centered at $\mu_0$, contained in $\Omega'$, and small enough that it meets $\text{spec}(\mathcal{L}_{0,c})$ only at $\mu_0$, we obtain a simple eigenvalue of $\mathcal{L}_{\delta,c}^{(s)}$ in it converging to $\mu_0$; being negative, it must coincide with $\mu_\delta$, so $\mu_\delta\to\mu_0$. This proves item 2.

    The final assertion of the lemma is now immediate, since conditions 1, 2, and 3 in Theorem \ref{Theorem: Stability instability criterion Albert} amount to items 1 and 2 above.
\end{proof}

The same perturbative argument, with gBO playing the role of gKdV, gives the corresponding spectral information for the solitary wave operators $\mathcal{L}_{\delta,c}$ of gILW in the deep water regime, now without any rescaling. Here, $\phi_{\infty,c}$ is the limit given by Theorem \ref{Th: Deep water uniform convergence of solitary waves}, which, by Proposition \ref{Prop: gBO waves are FL ground states} and Remark \ref{Remark: consequences of FL for the deep water limit}, is the unique suitable solitary wave of gBO with velocity $c$.

\begin{lemma}[\textbf{Spectral conditions in deep water}]\label{Lemma: Spectral conditions for gILW operator in deep water}
    Let $c>0$ and $p\in\N$ be fixed, and let $\mu_\infty$, $\gamma_\infty>0$ be as in Proposition \ref{Prop: gBO waves are FL ground states}. Then there exists some $\Lambda=\Lambda(c,p)>0$ such that, for any $\delta>\Lambda$, the solitary wave operator $\mathcal{L}_{\delta,c}$ satisfies:
    \begin{enumerate}
        \item the eigenvalue $0$ of $\mathcal{L}_{\delta,c}$ is simple, and $\ker\mathcal{L}_{\delta,c}=\text{span}\{\partial_x\phi_{\delta,c}\}$;
        \item $\mathcal{L}_{\delta,c}$ has exactly one negative eigenvalue $\mu_\delta^{(d)}$, it is simple, and $\mu_\delta^{(d)}\to\mu_\infty$ as $\delta\to\infty$;
        \item $\gamma_\delta^{(d)} := \text{dist}\left(0,\text{spec}(\mathcal{L}_{\delta,c})\backslash\{0\}\right) \geq \dfrac{\gamma_\infty}{2}$.
    \end{enumerate}
    In particular, $\mathcal{L}_{\delta,c}$ satisfies conditions 1, 2, and 3 in Theorem \ref{Theorem: Stability instability criterion Albert}.
\end{lemma}
\begin{proof}
    The proof is the same as that of Lemma \ref{Lemma: Spectral conditions for shallow-gILW operator}, with the limit $\delta\to0$ replaced by $\delta\to\infty$, the operators $\mathcal{L}_{\delta,c}^{(s)}$ and $\mathcal{L}_{0,c}$ replaced by $\mathcal{L}_{\delta,c}$ and $\mathcal{L}_{\infty,c}$, respectively, and Theorem \ref{Theorem: Spectrum of gKdV solitary wave operator} replaced by item 2 of Proposition \ref{Prop: gBO waves are FL ground states}; the convergence $d(\mathcal{L}_{\delta,c},\mathcal{L}_{\infty,c})\to0$ as $\delta\to\infty$, used throughout that proof, is item 1 of Theorem \ref{Theorem: Convergence of solitary wave operators}. The only other ingredient of that proof which depends on the regime is the uniform lower bound \eqref{Eq: uniform lower bound for the spectra}, which also holds in deep water: since $y\coth y\geq1$ for every $y\in\R$, we have $m_\delta\geq0$, and $|\phi_{\delta,c}|_\infty\leq C_\infty$ uniformly for $\delta\gg1$ by Proposition \ref{Prop: Uniform estimates of solitary waves} with $j=d$ and $r=\infty$; hence
    \begin{equation*}
        \left\langle\mathcal{L}_{\delta,c}u,u\right\rangle \geq \left(c-C_\infty^p\right)|u|_2^2,
    \end{equation*}
    for every $u$ in the domain of $\mathcal{L}_{\delta,c}$ and every $\delta\gg1$, and analogously for $\mathcal{L}_{\infty,c}$, since $\mathcal{H}\partial_x\geq0$ and $\phi_{\infty,c}\in L^\infty$.
\end{proof}

\begin{lemma}\label{Lemma: symbols m s d delta satisfies C1 C2 C3}
    For every $\delta>0$, the symbol $m_\delta$ of gILW satisfies conditions (C1), (C2) and (C3).
\end{lemma}
\begin{proof}
    Since $m^{(s)}_\delta=(3/\delta)m_\delta$ and conditions (C1), (C2) and (C3) are invariant under multiplication of the symbol by a positive constant, it is enough to prove them for $m^{(s)}_\delta$.
    
    \emph{Step 1: condition C1.} It follows from its definition that $m^{(s)}_\delta \in C^\infty(\R\backslash\{0\})$. The first derivative of the symbol is given by
    \begin{equation*}
        (m^{(s)}_\delta)'(\xi) = \frac{6\pi}{\delta}\left[\coth{(2\pi\delta\xi) - 2\pi\delta\xi\text{csch}^2{(2\pi\delta\xi)}}\right],
    \end{equation*}
    and one can show, by L'Hôpital's rule, that it satisfies the limit
    \begin{equation*}
        \lim_{\xi\to0}(m^{(s)}_\delta)'(\xi) = 0,
    \end{equation*}
    so $m^{(s)}_\delta \in C^1(\R)$, and therefore satisfies condition C1.

    \emph{Step 2: condition C2.} Let $l(\xi):= \xi m^{(s)}_\delta(\xi)$. Then,
    \begin{equation*}
        l'(\xi) = m^{(s)}_\delta(\xi) + \xi (m^{(s)}_\delta)'(\xi),
    \end{equation*}
    and
    \begin{align*}
        l''(\xi) &= 2(m^{(s)}_\delta)'(\xi) + \xi (m^{(s)}_\delta)''(\xi) \\
        &=  \frac{12\pi}{\delta}\left[\coth{(2\pi\delta\xi) - 2\pi\delta\xi\text{csch}^2{(2\pi\delta\xi)}}\right] \\
        &\quad + \frac{6\pi\xi}{\delta} [-2\pi\delta\text{csch}^2{(2\pi\delta\xi)}-2\pi\delta\text{csch}^2{(2\pi\delta\xi)} + 8\pi^2\delta^2\xi\text{csch}^2{(2\pi\delta\xi)}\coth{(2\pi\delta\xi)}] \\
        &= \frac{12\pi}{\delta}\left[\coth{(2\pi\delta\xi) - 2\pi\delta\xi\text{csch}^2{(2\pi\delta\xi)}}\right] \\
        &\quad + 24\pi^2\xi[-\text{csch}^2{(2\pi\delta\xi)} + 2\pi\delta\xi\text{csch}^2{(2\pi\delta\xi)}\coth{(2\pi\delta\xi)}].\\
        &= \frac{12\pi}{\delta}\coth{(2\pi\delta\xi)} - 48\pi^2\xi\text{csch}^2{(2\pi\delta\xi)} + 48\pi^3\delta\xi^2\text{csch}^2{(2\pi\delta\xi)}\coth(2\pi\delta\xi).
    \end{align*}
    Note that using the L'Hôpital's rule we can prove that $\xi\text{csch}(2\pi\delta\xi)$ is bounded away from $0$, and that it is enough for $l''$ to be bounded away from $0$ too.

    Let $x:=2\pi\delta\xi$, then $l''(\xi)$ can be written in terms of $x$ as
    \begin{align*}
        l''(\xi) &= \frac{12\pi}{\delta}\coth{x} - \frac{24\pi x}{\delta}\text{csch}^2{x} + \frac{12\pi x^2}{\delta}\text{csch}^2{x}\coth{x} \\
        &= \frac{12\pi}{\delta}\frac{\cosh{x}}{\sinh{x}} - \frac{24\pi x}{\delta}\frac{1}{\sinh^2{x}} + \frac{12\pi x^2}{\delta}\frac{\cosh{x}}{\sinh^3{x}} \\
        &= \frac{12\pi}{\delta}\left[\frac{\sinh^2{x}\cosh{x}-2x\sinh{x} + x^2\cosh{x}}{\sinh^3{x}}\right].
    \end{align*}
    Since $(\sinh{x} - x)^2 \geq 0$ and equality holds only for $x=0$, we have
    \begin{equation*}
        \sinh^2(x) + x^2 > 2x\sinh{x}, \quad \forall x\neq0,
    \end{equation*}
    then, using $\cosh{x}>1$ for any $x\in\R\backslash\{0\}$, 
    \begin{equation*}
        \sinh^2{x}\cosh{x}-2x\sinh{x} + x^2\cosh{x} > \cosh{x}(\sinh^2{x}+x^2) - 2x\sinh{x} > 0.
    \end{equation*}
    Therefore, the only zero of $l''$ is $\xi=0$. 
    
    To show that the zero is non-degenerate we analyze the behavior of $l''$ around $\xi=0$. The behavior of $\sinh$ and $\cosh$ around $0$ is given by
    \begin{equation*}
        \sinh(x) = x + \frac{x^3}{6} + O(x^5),
    \end{equation*}
    and
    \begin{equation*}
        \cosh(x) = 1 + \frac{x^2}{2} + O(x^4).
    \end{equation*}
    Then, letting 
    \begin{equation*}
        g(x):=\frac{\sinh^2{x}\cosh{x}-2x\sinh{x} + x^2\cosh{x}}{\sinh^3{x}},
    \end{equation*}
    we have
    \begin{align*}
        g(x) &= \frac{x^4 + O(x^6)}{x^3+ O(x^5)} = x + O(x^3) \\
        g'(x) &= 1 + O(x^2) \\
        g'(0) &= 1.
    \end{align*}
    Now, by the chain rule,
    \begin{equation*}
        l'''(\xi)\Big|_{\xi=0} = \frac{12\pi}{\delta}g'(x)2\pi\delta\Big|_{x=0} = 24\pi^2.
    \end{equation*}
    Thus, $\xi=0$ is a non-degenerate zero of $l''$. This shows that $m_{\delta}^{(s)}$ satisfies condition C2.

    \emph{Step 3: condition C3.} Since
    \begin{equation*}
        l''(\xi) = \frac{12\pi}{\delta}g(x) = 24\pi^2\xi + O(\xi^3), \quad \forall |\xi| \ll 1,
    \end{equation*}
    the symbol satisfies condition C3.
\end{proof}

The hypotheses of Theorem \ref{Theorem: Stability instability criterion Albert}, as well as the sign of $I$, are invariant under the scaling $D_\theta$ of gILW given by Lemma \ref{Lemma: Scaling symmetries}.

\begin{lemma}[\textbf{Invariance under the scaling}]\label{Lemma: Scaling invariance of the criterion}
    Let $p\in\N$ and $\delta_1,c_1>0$, and let $\phi_{\delta_1,c_1}$ be a suitable solitary wave of gILW with depth parameter $\delta_1$ and velocity $c_1$, which satisfies condition (i) in Definition \ref{Def: Orbital stability} and whose solitary wave operator satisfies conditions 1, 2 and 3 in Theorem \ref{Theorem: Stability instability criterion Albert}. Let $\delta_2,c_2>0$ be such that $\delta_2c_2=\delta_1c_1$, set $\theta:=\delta_1/\delta_2=c_2/c_1$, and let $\phi_{\delta_2,c_2}:=D_\theta\phi_{\delta_1,c_1}$. Then:
    \begin{enumerate}
        \item $\phi_{\delta_2,c_2}$ is a suitable solitary wave of gILW with depth parameter $\delta_2$ and velocity $c_2$, which satisfies condition (i) in Definition \ref{Def: Orbital stability}, and whose solitary wave operator satisfies conditions 1, 2 and 3 in Theorem \ref{Theorem: Stability instability criterion Albert};
        \item if $\psi_1\in L^2$ solves $\mathcal{L}_{\delta_1,c_1}\psi_1=\phi_{\delta_1,c_1}$, then $\psi_2:=\theta^{1/p-1}\tau_\theta\psi_1$ solves $\mathcal{L}_{\delta_2,c_2}\psi_2=\phi_{\delta_2,c_2}$, where $\tau_\theta f(x)=f(\theta x)$, and
        \begin{equation*}
            \left\langle\psi_2,\phi_{\delta_2,c_2}\right\rangle = \theta^{\frac{2}{p}-2}\left\langle\psi_1,\phi_{\delta_1,c_1}\right\rangle.
        \end{equation*}
    \end{enumerate}
    In particular, $\phi_{\delta_2,c_2}$ is orbitally stable if $\langle\psi_1,\phi_{\delta_1,c_1}\rangle<0$, and unstable if $\langle\psi_1,\phi_{\delta_1,c_1}\rangle>0$.
\end{lemma}
\begin{proof}
    By items 3 and 5 of Lemma \ref{Lemma: Scaling symmetries}, $\phi_{\delta_2,c_2}$ is a suitable solitary wave of gILW with depth parameter $\delta_2$ and velocity $c_2$, and, by item 6 of the same lemma, it satisfies condition (i). By item 4,
    \begin{equation*}
        \mathcal{L}_{\delta_2,c_2}=\theta\left[\tau_\theta\mathcal{L}_{\delta_1,c_1}\tau_\theta^{-1}\right],
    \end{equation*}
    so that the spectrum of $\mathcal{L}_{\delta_2,c_2}$ is that of $\mathcal{L}_{\delta_1,c_1}$ multiplied by $\theta>0$, with corresponding eigenfunctions related by $\tau_\theta$ and the same multiplicities; hence conditions 1, 2 and 3 in Theorem \ref{Theorem: Stability instability criterion Albert} hold for $\mathcal{L}_{\delta_2,c_2}$. This proves item 1. For item 2,
    \begin{equation*}
        \mathcal{L}_{\delta_2,c_2}\psi_2 = \theta\,\tau_\theta\mathcal{L}_{\delta_1,c_1}\tau_\theta^{-1}\left[\theta^{1/p-1}\tau_\theta\psi_1\right] = \theta^{1/p}\tau_\theta\phi_{\delta_1,c_1} = \phi_{\delta_2,c_2},
    \end{equation*}
    and the identity for $\langle\psi_2,\phi_{\delta_2,c_2}\rangle$ follows from the change of variables $y=\theta x$. The last assertion follows from items 1 and 2, Lemma \ref{Lemma: symbols m s d delta satisfies C1 C2 C3} and Theorem \ref{Theorem: Stability instability criterion Albert}, since $\theta^{2/p-2}>0$.
\end{proof}

\subsection{The shallow water regime}\label{Section: Shallow water regime}

In the slow regime, the sign of $I$ is studied through the rescaled operators $\mathcal{L}_{\delta,c}^{(s)}=(3/\delta)\mathcal{L}_{\delta,\delta c/3}$, which converge to $\mathcal{L}_{0,c}$ as $\delta\to0$. Since the equation $\mathcal{L}_{\delta,c}^{(s)}\psi=\phi_{\delta,c}^{(s)}$ does not determine $\psi$ uniquely, our first task is to single out a canonical solution $\psi_{\delta,c}^{(s)}$, uniform in $\delta$, on which $I$ can be computed and its limit as $\delta\to0$ can be studied.

\begin{lemma}\label{Lemma: L2 bound for psi in shallow water}
    Let $c>0$ and $p\in\mathbb{N}$ be fixed. For any $\delta>0$, let $\phi_{\delta,c}^{(s)}$ be the rescaled suitable solitary wave of Section \ref{Section: Convergence of solitary waves}. Then, there exists some $\epsilon=\epsilon(c,p)>0$ such that, for any $0<\delta<\epsilon$, the equation
    \begin{equation}\label{Eq: Equation for psi in shallow water}
        \mathcal{L}_{\delta,c}^{(s)}\psi = \phi_{\delta,c}^{(s)}
    \end{equation}
    admits a unique solution $\psi_{\delta,c}^{(s)}\in L^2$ orthogonal to $\partial_x\phi_{\delta,c}^{(s)}$, the set of all its $L^2$-solutions being $\psi_{\delta,c}^{(s)}+\text{span}\{\partial_x\phi_{\delta,c}^{(s)}\}$. Moreover,
    \begin{equation*}
        |\psi_{\delta,c}^{(s)}|_{2} \lesssim_{c,p} 1,
    \end{equation*}
    and the quantity $\langle\psi,\phi_{\delta,c}^{(s)}\rangle$ takes the same value for every $L^2$-solution $\psi$ of \eqref{Eq: Equation for psi in shallow water}.
\end{lemma}
\begin{proof}
    Note first that the normalization $\psi\perp\partial_x\phi_{\delta,c}^{(s)}$ cannot be dispensed with: by Lemma \ref{Lemma: Spectral conditions for shallow-gILW operator}, $\partial_x\phi_{\delta,c}^{(s)}$ spans the kernel of $\mathcal{L}_{\delta,c}^{(s)}$, so that together with any solution $\psi$ of \eqref{Eq: Equation for psi in shallow water} the whole family $\psi+t\,\partial_x\phi_{\delta,c}^{(s)}$, $t\in\R$, also solves it, and the $L^2$-norms of the members of this family are unbounded. The last assertion of the lemma shows, however, that this normalization is harmless for our purposes, since the number $I=\langle\psi,\phi_{\delta,c}^{(s)}\rangle$ appearing in Theorem \ref{Theorem: Stability instability criterion Albert} does not depend on which solution is used to compute it.

    Throughout the proof, $\epsilon=\epsilon(c,p)>0$ denotes the constant given by Lemma \ref{Lemma: Spectral conditions for shallow-gILW operator}, and $\gamma_0>0$ is as in Theorem \ref{Theorem: Spectrum of gKdV solitary wave operator}. We shall freely use items 1 and 3 of that lemma, namely that, for any $0<\delta<\epsilon$,
    \begin{equation}\label{Eq: kernel and gap of L delta c s}
        \ker\mathcal{L}_{\delta,c}^{(s)} = \text{span}\{\partial_x\phi_{\delta,c}^{(s)}\} \quad\text{ and }\quad \gamma_\delta:=\text{dist}\left(0,\text{spec}(\mathcal{L}_{\delta,c}^{(s)})\backslash\{0\}\right)\geq\frac{\gamma_0}{2}.
    \end{equation}

    \emph{Step 1: solvability of \eqref{Eq: Equation for psi in shallow water} and the choice of $\psi_{\delta,c}^{(s)}$.} Fix $0<\delta<\epsilon$ and let $N_\delta:=(\ker\mathcal{L}_{\delta,c}^{(s)})^{\perp}$. Since $\phi_{\delta,c}^{(s)}\in H^\infty$, we have
    \begin{equation*}
        \left\langle \phi_{\delta,c}^{(s)},\partial_x\phi_{\delta,c}^{(s)}\right\rangle = \int_\R \phi_{\delta,c}^{(s)}\partial_x\phi_{\delta,c}^{(s)} = \frac{1}{2}\int_\R\partial_x\left[(\phi_{\delta,c}^{(s)})^2\right] = 0,
    \end{equation*}
    so that $\phi_{\delta,c}^{(s)}\in N_\delta$ by \eqref{Eq: kernel and gap of L delta c s}. Here, the integral converges absolutely because $\phi_{\delta,c}^{(s)},\partial_x\phi_{\delta,c}^{(s)}\in L^2$, so that $\phi_{\delta,c}^{(s)}\partial_x\phi_{\delta,c}^{(s)}\in L^1$ by the Cauchy-Schwarz inequality; and it vanishes because, by the fundamental theorem of calculus,
    \begin{equation*}
        \int_\R\partial_x\left[(\phi_{\delta,c}^{(s)})^2\right] = \lim_{R\to\infty}\left[\phi_{\delta,c}^{(s)}(R)^2 - \phi_{\delta,c}^{(s)}(-R)^2\right] = 0,
    \end{equation*}
    since $\phi_{\delta,c}^{(s)}(x)\to0$ as $|x|\to\infty$. Moreover, $\mathcal{L}_{\delta,c}^{(s)}$ is self-adjoint by Theorem \ref{Theorem: Basic information about spectrum}, and $0$ is an isolated point of its spectrum whose associated eigenspace is exactly $\ker\mathcal{L}_{\delta,c}^{(s)}$. Hence, by the spectral theorem, the closed subspace $N_\delta=(\ker\mathcal{L}_{\delta,c}^{(s)})^{\perp}$ reduces $\mathcal{L}_{\delta,c}^{(s)}$ — that is, both $N_\delta$ and $\ker\mathcal{L}_{\delta,c}^{(s)}$ are invariant under $\mathcal{L}_{\delta,c}^{(s)}$ — and the restriction $\mathcal{L}_{\delta,c}^{(s)}|_{N_\delta}$ is a self-adjoint operator on $N_\delta$ with
    \begin{equation*}
        \text{spec}\left(\mathcal{L}_{\delta,c}^{(s)}|_{N_\delta}\right) = \text{spec}(\mathcal{L}_{\delta,c}^{(s)})\backslash\{0\} \subset \{\lambda\in\R: |\lambda|\geq\gamma_\delta\}.
    \end{equation*}
    In particular, $\mathcal{L}_{\delta,c}^{(s)}|_{N_\delta}$ is invertible and its inverse is a bounded, self-adjoint operator on $N_\delta$ with
    \begin{equation}\label{Eq: bound for the inverse on N delta}
        \left\|\left(\mathcal{L}_{\delta,c}^{(s)}|_{N_\delta}\right)^{-1}\right\|_{2,2} \leq \frac{1}{\gamma_\delta}.
    \end{equation}
    Estimate \eqref{Eq: bound for the inverse on N delta} is a consequence of the fact that, for a bounded self-adjoint operator $B$, the operator norm and the spectral radius coincide, $\|B\|_{2,2} = \sup\{|\lambda|:\lambda\in\text{spec}(B)\}$. Applying this to $B=(\mathcal{L}_{\delta,c}^{(s)}|_{N_\delta})^{-1}$, whose spectrum is $\{1/\lambda:\lambda\in\text{spec}(\mathcal{L}_{\delta,c}^{(s)}|_{N_\delta})\}\subset\{z:|z|\leq1/\gamma_\delta\}$ by the spectral mapping theorem, we obtain \eqref{Eq: bound for the inverse on N delta}.
    We may therefore define
    \begin{equation*}
        \psi_{\delta,c}^{(s)} := \left(\mathcal{L}_{\delta,c}^{(s)}|_{N_\delta}\right)^{-1}\phi_{\delta,c}^{(s)} \in N_\delta,
    \end{equation*}
    which is the unique solution of \eqref{Eq: Equation for psi in shallow water} lying in $N_\delta$, that is, the unique solution orthogonal to $\partial_x\phi_{\delta,c}^{(s)}$. Since two solutions of \eqref{Eq: Equation for psi in shallow water} differ by an element of $\ker\mathcal{L}_{\delta,c}^{(s)}$, it follows from \eqref{Eq: kernel and gap of L delta c s} that the set of all $L^2$-solutions of \eqref{Eq: Equation for psi in shallow water} is exactly $\psi_{\delta,c}^{(s)}+\text{span}\{\partial_x\phi_{\delta,c}^{(s)}\}$.

    \emph{Step 2: the uniform bound.} Combining \eqref{Eq: bound for the inverse on N delta} with the lower bound for $\gamma_\delta$ in \eqref{Eq: kernel and gap of L delta c s}, we get, for any $0<\delta<\epsilon$,
    \begin{equation*}
        |\psi_{\delta,c}^{(s)}|_2 \leq \frac{1}{\gamma_\delta}\left|\phi_{\delta,c}^{(s)}\right|_2 \leq \frac{2}{\gamma_0}\left|\phi_{\delta,c}^{(s)}\right|_2 \lesssim_{c,p} 1,
    \end{equation*}
    where the last estimate is Proposition \ref{Prop: Uniform estimates of solitary waves} with $r=2$ and $j=s$, valid for $\delta\ll1$, and where we used that $\gamma_0$ depends only on $c$ and $p$. Shrinking $\epsilon$ if necessary, the above holds for any $0<\delta<\epsilon$.

    \emph{Step 3: independence of $\langle\psi,\phi_{\delta,c}^{(s)}\rangle$ from the choice of solution.} By Step 1, any $L^2$-solution of \eqref{Eq: Equation for psi in shallow water} is of the form $\psi=\psi_{\delta,c}^{(s)}+t\,\partial_x\phi_{\delta,c}^{(s)}$, $t\in\R$, and $\langle\partial_x\phi_{\delta,c}^{(s)},\phi_{\delta,c}^{(s)}\rangle=0$.
\end{proof}

Having fixed the normalized solution $\psi_{\delta,c}^{(s)}$, we can now speak unambiguously of $I_\delta:=\langle\psi_{\delta,c}^{(s)},\phi_{\delta,c}^{(s)}\rangle$. The next lemma shows that $I_\delta$ converges, as $\delta\to0$, to the analogous quantity $I_0$ for the limit gKdV operator; its computation is left for the proof of Theorem \ref{Th: Stability and instability for shallow-gILW} below.

\begin{lemma}\label{Lemma: Convergence of the Albert functional}
    Let $c>0$ and $p\in\N$ be fixed, and let $\epsilon=\epsilon(c,p)>0$ and $\psi_{\delta,c}^{(s)}$ be as in Lemma \ref{Lemma: L2 bound for psi in shallow water}. Let also $\psi_{0,c}\in L^2$ be the unique solution of $\mathcal{L}_{0,c}\psi=\phi_{0,c}$ orthogonal to $\partial_x\phi_{0,c}$, and define
    \begin{equation*}
        I_\delta := \left\langle\psi_{\delta,c}^{(s)},\phi_{\delta,c}^{(s)}\right\rangle, \quad 0<\delta<\epsilon, \qquad\text{ and }\qquad I_0 := \left\langle\psi_{0,c},\phi_{0,c}\right\rangle.
    \end{equation*}
    Then $I_\delta\to I_0$ as $\delta\to0$.
\end{lemma}
\begin{proof}
    Note first that Step 1 in the proof of Lemma \ref{Lemma: L2 bound for psi in shallow water} applies verbatim to $\mathcal{L}_{0,c}$ instead of $\mathcal{L}_{\delta,c}^{(s)}$, since Theorem \ref{Theorem: Spectrum of gKdV solitary wave operator} provides for $\mathcal{L}_{0,c}$ exactly the spectral information that Lemma \ref{Lemma: Spectral conditions for shallow-gILW operator} provides for $\mathcal{L}_{\delta,c}^{(s)}$. Hence $\psi_{0,c}$ is well defined, and the set of all $L^2$-solutions of $\mathcal{L}_{0,c}\psi=\phi_{0,c}$ is $\psi_{0,c}+\text{span}\{\partial_x\phi_{0,c}\}$.

    The proof is by a weak compactness argument. Let us emphasize from the outset that the conclusion concerns the scalar quantities $I_\delta$ and not the functions $\psi_{\delta,c}^{(s)}$ themselves: as will be seen in Step 3, the weak limits of $\{\psi_{\delta,c}^{(s)}\}$ are determined only up to an element of $\ker\mathcal{L}_{0,c}$, and it is the last assertion of Lemma \ref{Lemma: L2 bound for psi in shallow water} which makes this ambiguity harmless.

    Let $\{\delta_n\}\subset(0,\epsilon)$ be an arbitrary sequence with $\delta_n\to0$, and abbreviate $\psi_n:=\psi_{\delta_n,c}^{(s)}$ and $\phi_n:=\phi_{\delta_n,c}^{(s)}$.

    \emph{Step 1: extraction of a weakly convergent subsequence.} By Lemma \ref{Lemma: L2 bound for psi in shallow water} we have $|\psi_n|_2\lesssim_{c,p}1$, so that $\{\psi_n\}$ is a bounded sequence in the Hilbert space $L^2$. Hence there exists a subsequence, which we still denote by $\{\psi_n\}$, and some $\psi_*\in L^2$, such that
    \begin{equation*}
        \psi_n \rightharpoonup \psi_* \quad\text{ weakly in } L^2.
    \end{equation*}

    \emph{Step 2: the weak limit solves the limit equation.} We claim that $\mathcal{L}_{0,c}\psi_*=\phi_{0,c}$. Let $v\in\mathcal{S}(\R)$ be arbitrary. Since $m_\delta^{(s)}$ grows at most quadratically, $v$ belongs to the domain of $\mathcal{L}_{\delta_n,c}^{(s)}$, which is self-adjoint by Theorem \ref{Theorem: Basic information about spectrum}; therefore
    \begin{equation}\label{Eq: weak formulation for psi n}
        \left\langle\psi_n,\mathcal{L}_{\delta_n,c}^{(s)}v\right\rangle = \left\langle\mathcal{L}_{\delta_n,c}^{(s)}\psi_n,v\right\rangle = \left\langle\phi_n,v\right\rangle.
    \end{equation}
    We pass to the limit on both sides of \eqref{Eq: weak formulation for psi n}. On the right-hand side, $\langle\phi_n,v\rangle\to\langle\phi_{0,c},v\rangle$ by Corollary \ref{Cor: L2 convergence shallow water}. On the left-hand side we split
    \begin{equation*}
        \left\langle\psi_n,\mathcal{L}_{\delta_n,c}^{(s)}v\right\rangle = \left\langle\psi_n,\mathcal{L}_{0,c}v\right\rangle + \left\langle\psi_n,\left(\mathcal{L}_{\delta_n,c}^{(s)}-\mathcal{L}_{0,c}\right)v\right\rangle.
    \end{equation*}
    The first term converges to $\langle\psi_*,\mathcal{L}_{0,c}v\rangle$, because $\mathcal{L}_{0,c}v\in L^2$ is a fixed function and $\psi_n\rightharpoonup\psi_*$. The second term is bounded in modulus by $|\psi_n|_2\left|(\mathcal{L}_{\delta_n,c}^{(s)}-\mathcal{L}_{0,c})v\right|_2$, so, the first factor being bounded, it suffices to check that
    \begin{equation}\label{Eq: strong convergence of the operators on Schwartz functions}
        \left|\left(\mathcal{L}_{\delta,c}^{(s)}-\mathcal{L}_{0,c}\right)v\right|_2 \leq \left|\left(M_\delta^{(s)}+\partial_x^2\right)v\right|_2 + \left|\left((\phi_{\delta,c}^{(s)})^p-\phi_{0,c}^p\right)v\right|_2 \longrightarrow 0, \quad \text{ as }\delta\to0.
    \end{equation}
    For the first summand in \eqref{Eq: strong convergence of the operators on Schwartz functions}, recall from Lemma 8.2.1 in \cite{Abdelouhab1989} that $m_\delta^{(s)}(\xi)=4\pi^2\xi^2-\xi^2h(\xi,\delta)$, with $0\leq\xi^2h(\xi,\delta)\leq96\delta^2\xi^4\pi^4/90$ as shown in the proof of Corollary \ref{Cor: L2 convergence shallow water}; this gives simultaneously the point-wise convergence $m_\delta^{(s)}(\xi)\to4\pi^2\xi^2$ and the uniform bound $0\leq m_\delta^{(s)}(\xi)\leq4\pi^2\xi^2$. Hence, by Plancherel's theorem,
    \begin{equation*}
        \left|\left(M_\delta^{(s)}+\partial_x^2\right)v\right|_2^2 = \int_\R\left|m_\delta^{(s)}(\xi)-4\pi^2\xi^2\right|^2\left|\fourier{v}(\xi)\right|^2d\xi \longrightarrow 0,
    \end{equation*}
    by the dominated convergence theorem, the dominating function being $16\pi^4\xi^4|\fourier{v}(\xi)|^2$, which is integrable because $v\in\mathcal{S}(\R)$. For the second summand in \eqref{Eq: strong convergence of the operators on Schwartz functions}, the mean value theorem gives $|a^p-b^p|\leq p\max\{a,b\}^{p-1}|a-b|$ for $a,b\geq0$, so that
    \begin{equation*}
        \left|\left((\phi_{\delta,c}^{(s)})^p-\phi_{0,c}^p\right)v\right|_2 \leq p\left(\left|\phi_{\delta,c}^{(s)}\right|_\infty+\left|\phi_{0,c}\right|_\infty\right)^{p-1}\left|\phi_{\delta,c}^{(s)}-\phi_{0,c}\right|_\infty|v|_2 \longrightarrow 0,
    \end{equation*}
    by Theorem \ref{Th: Shallow water uniform convergence of solitary waves} together with the $\delta$-uniform bound of Proposition \ref{Prop: Uniform estimates of solitary waves} with $r=\infty$. This proves \eqref{Eq: strong convergence of the operators on Schwartz functions}.

    Passing to the limit in \eqref{Eq: weak formulation for psi n} we therefore obtain
    \begin{equation*}
        \left\langle\psi_*,\mathcal{L}_{0,c}v\right\rangle = \left\langle\phi_{0,c},v\right\rangle, \quad \forall v\in\mathcal{S}(\R),
    \end{equation*}
    that is, $\mathcal{L}_{0,c}\psi_*=\phi_{0,c}$ in the sense of distributions. In particular $-\partial_x^2\psi_*=\phi_{0,c}-c\psi_*+\phi_{0,c}^p\psi_*\in L^2$, so that $\psi_*\in H^2$, which is the domain of $\mathcal{L}_{0,c}$, and the identity holds in $L^2$. This proves the claim.

    \emph{Step 3: the pairing with $\phi_{0,c}$ does not depend on the subsequence.} By the first paragraph of this proof and Step 2, there is some $t\in\R$, a priori depending on the extracted subsequence, such that $\psi_*=\psi_{0,c}+t\,\partial_x\phi_{0,c}$. Nevertheless, by Step 3 in the proof of Lemma \ref{Lemma: L2 bound for psi in shallow water},
    \begin{equation}\label{Eq: pairing of the weak limit}
        \left\langle\psi_*,\phi_{0,c}\right\rangle = \left\langle\psi_{0,c},\phi_{0,c}\right\rangle + t\left\langle\partial_x\phi_{0,c},\phi_{0,c}\right\rangle = \left\langle\psi_{0,c},\phi_{0,c}\right\rangle = I_0,
    \end{equation}
    whatever the value of $t$. Let us point out that we cannot assert that $t=0$: although $\psi_n$ is orthogonal to $\partial_x\phi_n$ for every $n$, passing to the limit in this orthogonality relation would require the strong convergence of $\partial_x\phi_{\delta,c}^{(s)}$ in $L^2$, which is not available to us. This is precisely why the argument is carried out at the level of the scalars $I_\delta$ and not at the level of the functions $\psi_{\delta,c}^{(s)}$.

    \emph{Step 4: convergence of the pairing along the subsequence.} Write
    \begin{equation*}
        I_{\delta_n} = \left\langle\psi_n,\phi_n\right\rangle = \left\langle\psi_n,\phi_n-\phi_{0,c}\right\rangle + \left\langle\psi_n,\phi_{0,c}\right\rangle.
    \end{equation*}
    The first term is bounded in modulus, by the Cauchy-Schwarz inequality, by $|\psi_n|_2|\phi_n-\phi_{0,c}|_2$, which tends to $0$ because $\{|\psi_n|_2\}$ is bounded by Lemma \ref{Lemma: L2 bound for psi in shallow water} and $|\phi_n-\phi_{0,c}|_2\to0$ by Corollary \ref{Cor: L2 convergence shallow water}. The second term converges to $\langle\psi_*,\phi_{0,c}\rangle$, since $\psi_n\rightharpoonup\psi_*$ and $\phi_{0,c}$ is a fixed element of $L^2$. Hence, by \eqref{Eq: pairing of the weak limit},
    \begin{equation*}
        I_{\delta_n} \longrightarrow I_0 \quad\text{ along the extracted subsequence}.
    \end{equation*}

    \emph{Step 5: conclusion.} Steps 1--4 show that every subsequence of the real sequence $\{I_{\delta_n}\}_n$ admits a further subsequence converging to $I_0$, the limit being the same number $I_0$ for all such extractions. This forces $I_{\delta_n}\to I_0$: otherwise there would exist $\eta>0$ and a subsequence with $|I_{\delta_{n_k}}-I_0|\geq\eta$ for every $k$, and this subsequence would admit no further subsequence converging to $I_0$, a contradiction. Since the sequence $\{\delta_n\}\subset(0,\epsilon)$ with $\delta_n\to0$ was arbitrary, we conclude that $I_\delta\to I_0$ as $\delta\to0$.
\end{proof}

All the pieces are now in place: the symbol and spectral conditions of Theorem \ref{Theorem: Stability instability criterion Albert} are verified, and the sign of $I_\delta$ is determined, for $\delta$ small, by the sign of its limit $I_0$. We conclude with the computation of $I_0$ and the main result for slow solitary waves.

\begin{theorem}[\textbf{Slow solitary waves}]\label{Th: Stability and instability for shallow-gILW}
    Let $c>0$ and $p\in\mathbb{N}$, $p\neq4$, be fixed. There exists some $\epsilon=\epsilon(c,p)>0$ such that, for any $0<\delta<\epsilon$, if $\phi_{\delta,\delta c/3}$ is a suitable solitary wave of gILW, given by Theorem \ref{Th: Existence of suitable solitary waves}, with depth parameter $\delta$ and velocity $\delta c/3$, which satisfies condition (i) in Definition \ref{Def: Orbital stability}, then:
    \begin{enumerate}
        \item if $1\leq p < 4$, $\phi_{\delta,\delta c/3}$ is orbitally stable;
        \item if $p>4$, $\phi_{\delta,\delta c/3}$ is unstable.
    \end{enumerate}
\end{theorem}
\begin{proof}
    Let $\phi_{\delta,c}^{(s)}=(3/\delta)^{1/p}\phi_{\delta,\delta c/3}$ be the corresponding rescaled solitary wave. We apply Theorem \ref{Theorem: Stability instability criterion Albert} to $\phi_{\delta,\delta c/3}$. Its symbol $m_\delta$ satisfies conditions (C1), (C2) and (C3) by Lemma \ref{Lemma: symbols m s d delta satisfies C1 C2 C3}, and condition (i) holds by hypothesis. Since $\mathcal{L}_{\delta,\delta c/3}=(\delta/3)\mathcal{L}_{\delta,c}^{(s)}$, the three spectral conditions hold for $\mathcal{L}_{\delta,\delta c/3}$, for any $0<\delta<\epsilon_1$, because they hold for $\mathcal{L}_{\delta,c}^{(s)}$ by Lemma \ref{Lemma: Spectral conditions for shallow-gILW operator}, where $\epsilon_1=\epsilon_1(c,p)>0$ is the constant provided by that lemma. Moreover, if $\psi_{\delta,c}^{(s)}$ is the solution of $\mathcal{L}_{\delta,c}^{(s)}\psi=\phi_{\delta,c}^{(s)}$ given by Lemma \ref{Lemma: L2 bound for psi in shallow water}, then $\psi:=(\delta/3)^{1/p-1}\psi_{\delta,c}^{(s)}$ solves $\mathcal{L}_{\delta,\delta c/3}\psi=\phi_{\delta,\delta c/3}$, and
    \begin{equation*}
        \left\langle\psi,\phi_{\delta,\delta c/3}\right\rangle = \left(\frac{\delta}{3}\right)^{\frac{2}{p}-1} I_\delta, \qquad I_\delta = \left\langle\psi_{\delta,c}^{(s)},\phi_{\delta,c}^{(s)}\right\rangle.
    \end{equation*}
    Since, by Lemma \ref{Lemma: L2 bound for psi in shallow water}, the quantity $I$ in Theorem \ref{Theorem: Stability instability criterion Albert} does not depend on the solution used to compute it, it only remains to determine the sign of $I_\delta$, which we do in the two steps below.

    \emph{Step 1: computation of $I_0$.} By the explicit formula of Theorem \ref{Theorem: Existence and uniqueness for gKdV}, the solitary waves of gKdV satisfy the scaling relation $\phi_{0,c}(x)=c^{1/p}\phi_{0,1}(\sqrt c\,x)$, so that the map $c\mapsto\phi_{0,c}$ is differentiable with values in $H^2$, the domain of $\mathcal{L}_{0,c}$. Differentiating the gKdV solitary wave equation
    \begin{equation*}
        -\partial_x^2\phi_{0,c}+c\phi_{0,c}-\frac{1}{p+1}\phi_{0,c}^{p+1}=0
    \end{equation*}
    with respect to $c$, we get $\mathcal{L}_{0,c}[\partial_c\phi_{0,c}]=-\phi_{0,c}$. Since any two $L^2$-solutions of $\mathcal{L}_{0,c}\psi=\phi_{0,c}$ differ by a multiple of $\partial_x\phi_{0,c}$, which is orthogonal to $\phi_{0,c}$, we obtain (see also Lemma 4.1 in \cite{Albert1992})
    \begin{equation*}
        I_0 = -\left\langle\partial_c\phi_{0,c},\phi_{0,c}\right\rangle = -\frac{1}{2}\frac{d}{dc}\left|\phi_{0,c}\right|_2^2 = -\frac{1}{2}\frac{d}{dc}\left(c^{\frac{2}{p}-\frac{1}{2}}\left|\phi_{0,1}\right|_2^2\right),
    \end{equation*}
    that is,
    \begin{equation}\label{Eq: sign of I zero}
        I_0 = \frac{p-4}{4p}\,c^{\frac{2}{p}-\frac{3}{2}}\left|\phi_{0,1}\right|_2^2,
    \end{equation}
    so that $I_0<0$ if $p<4$, $I_0=0$ if $p=4$, and $I_0>0$ if $p>4$. This reproduces the threshold $p=2\mu$, where $\mu=2$ is the order of the symbol of gKdV, given by the Corollary of Theorem 3.1 in \cite{Bona1987}. In particular, $p=4$ is the only case in which the criterion of Theorem \ref{Theorem: Stability instability criterion Albert} gives no information.

    \emph{Step 2: conclusion.} By \eqref{Eq: sign of I zero} we have $I_0\neq0$, since $p\neq4$, and by Lemma \ref{Lemma: Convergence of the Albert functional} we have $I_\delta\to I_0$ as $\delta\to0$. Hence there exists $\epsilon_2=\epsilon_2(c,p)>0$ such that $I_\delta$ and $I_0$ have the same sign for any $0<\delta<\epsilon_2$. Setting $\epsilon:=\min\{\epsilon_1,\epsilon_2\}$, which depends only on $c$ and $p$, we conclude from Theorem \ref{Theorem: Stability instability criterion Albert} that, for any $0<\delta<\epsilon$: if $1\leq p<4$, then $I_\delta<0$ and $\phi_{\delta,\delta c/3}$ is orbitally stable; while if $p>4$, then $I_\delta>0$ and $\phi_{\delta,\delta c/3}$ is unstable.
\end{proof}

By Lemma \ref{Lemma: Scaling invariance of the criterion}, the conclusions of Theorem \ref{Th: Stability and instability for shallow-gILW} are not confined to the pairs $(\delta,\delta c/3)$: they hold for every pair of parameters with small product.

\begin{corollary}[\textbf{Slow solitary waves, general form}]\label{Cor: Stability and instability for gILW}
    Let $p\in\N$, $p\neq4$, let $\epsilon_0:=\epsilon(1,p)>0$ be given by Theorem \ref{Th: Stability and instability for shallow-gILW} with $c=1$, and set $\kappa_p:=\epsilon_0^2/3$. Let $\delta,c>0$ be such that $\delta c<\kappa_p$, set $\delta_1:=(3\delta c)^{1/2}\in(0,\epsilon_0)$, and let $\phi_{\delta_1,\delta_1/3}$ be a suitable solitary wave of gILW as in Theorem \ref{Th: Stability and instability for shallow-gILW}, with $c=1$, satisfying condition (i) in Definition \ref{Def: Orbital stability}. Then $\phi_{\delta,c}:=D_{\delta_1/\delta}\,\phi_{\delta_1,\delta_1/3}$ is a suitable solitary wave of gILW with depth parameter $\delta$ and velocity $c$, which satisfies condition (i) and:
    \begin{enumerate}
        \item if $1\leq p<4$, $\phi_{\delta,c}$ is orbitally stable;
        \item if $p>4$, $\phi_{\delta,c}$ is unstable.
    \end{enumerate}
\end{corollary}
\begin{proof}
    Since $\delta_1\cdot\frac{\delta_1}{3}=\delta c$, this follows from Theorem \ref{Th: Stability and instability for shallow-gILW} and Lemma \ref{Lemma: Scaling invariance of the criterion}, applied with $(\delta_1,c_1)=(\delta_1,\delta_1/3)$ and $(\delta_2,c_2)=(\delta,c)$.
\end{proof}

Corollary \ref{Cor: Stability and instability for gILW} shows that what has to be small is the product $\delta c$, and not the depth parameter alone: arbitrarily large depths are allowed, provided the wave is slow enough. Let us also point out that the case $p=4$ remains open, since the criterion of Theorem \ref{Theorem: Stability instability criterion Albert} gives no information when $I=0$; see \eqref{Eq: sign of I zero}.

\subsection{The deep water regime}\label{Section: Deep water regime}

We now turn to the deep water regime $\delta\to\infty$. In contrast with the shallow water regime, no rescaling of gILW is needed here: by Theorem \ref{Th: Deep water uniform convergence of solitary waves}, the suitable solitary waves $\phi_{\delta,c}$ of gILW converge uniformly, as $\delta\to\infty$, to the suitable solitary wave $\phi_{\infty,c}$ of gBO, which is unique and whose solitary wave operator $\mathcal{L}_{\infty,c}$ is described near $0$, for every admissible $p\geq1$, by Proposition \ref{Prop: gBO waves are FL ground states}, thanks to the work of Frank and Lenzmann \cite{Frank2013}. The argument is then the same as in the shallow water regime, with gBO playing the role of gKdV, and we only indicate the necessary changes in the proofs. The essential difference lies in the sign of the limit quantity $I_\infty$, which, as we show below, is determined by the threshold $p=2$ of gBO, instead of the threshold $p=4$ of gKdV.

Throughout this subsection, $c>0$ and $p\in\N$ are fixed and, as in Section \ref{Section: Convergence of solitary waves}, the suitable solitary waves $\phi_{\delta,c}$ of gILW are those obtained, through Theorem \ref{Th: Existence of suitable solitary waves}, from minimizers of the variational problems $I_{\delta,c,\lambda}$; by Theorems \ref{Th: Deep water convergence of solitary waves} and \ref{Th: Deep water uniform convergence of solitary waves}, they converge to $\phi_{\infty,c}=A\,\tau_cQ$ both in $H^{1/2}$ and uniformly, where $A=((p+1)c)^{1/p}$ and $Q$ are as in Proposition \ref{Prop: gBO waves are FL ground states}.

\begin{lemma}\label{Lemma: L2 bound for psi in deep water}
    There exists some $\Lambda=\Lambda(c,p)>0$ such that, for any $\delta>\Lambda$, the equation
    \begin{equation}\label{Eq: Equation for psi in deep water}
        \mathcal{L}_{\delta,c}\psi = \phi_{\delta,c}
    \end{equation}
    admits a unique solution $\psi_{\delta,c}\in L^2$ orthogonal to $\partial_x\phi_{\delta,c}$, the set of all its $L^2$-solutions being $\psi_{\delta,c}+\text{span}\{\partial_x\phi_{\delta,c}\}$. Moreover,
    \begin{equation*}
        |\psi_{\delta,c}|_{2} \lesssim_{c,p} 1,
    \end{equation*}
    and the quantity $\langle\psi,\phi_{\delta,c}\rangle$ takes the same value for every $L^2$-solution $\psi$ of \eqref{Eq: Equation for psi in deep water}.
\end{lemma}
\begin{proof}
    The proof is the same as that of Lemma \ref{Lemma: L2 bound for psi in shallow water}, with $\Lambda$ given by Lemma \ref{Lemma: Spectral conditions for gILW operator in deep water} instead of $\epsilon$ given by Lemma \ref{Lemma: Spectral conditions for shallow-gILW operator}, and with $\gamma_0$ replaced by $\gamma_\infty$. Indeed, by items 1 and 3 of Lemma \ref{Lemma: Spectral conditions for gILW operator in deep water}, for any $\delta>\Lambda$,
    \begin{equation*}
        \ker\mathcal{L}_{\delta,c} = \text{span}\{\partial_x\phi_{\delta,c}\} \quad\text{ and }\quad \gamma_\delta^{(d)}\geq\frac{\gamma_\infty}{2},
    \end{equation*}
    which is the counterpart of \eqref{Eq: kernel and gap of L delta c s}. With this information, Steps 1 and 3 of that proof apply verbatim, while in Step 2 we use the uniform bound $|\phi_{\delta,c}|_2\lesssim_c1$ for $\delta\gg1$, given by Proposition \ref{Prop: Uniform estimates of solitary waves} with $r=2$ and $j=d$, to get
    \begin{equation*}
        |\psi_{\delta,c}|_2 \leq \frac{2}{\gamma_\infty}\left|\phi_{\delta,c}\right|_2 \lesssim_{c,p} 1,
    \end{equation*}
    after enlarging $\Lambda$ if necessary.
\end{proof}

Having fixed the normalized solution $\psi_{\delta,c}$, we can define $I_\delta^{(d)}:=\langle\psi_{\delta,c},\phi_{\delta,c}\rangle$ for $\delta>\Lambda$. As in the shallow water regime, $I_\delta^{(d)}$ converges to the analogous quantity for the limit operator.

\begin{lemma}\label{Lemma: Convergence of the Albert functional in deep water}
    Let $\Lambda=\Lambda(c,p)>0$ and $\psi_{\delta,c}$ be as in Lemma \ref{Lemma: L2 bound for psi in deep water}. Let also $\psi_{\infty,c}\in L^2$ be the unique solution of $\mathcal{L}_{\infty,c}\psi=\phi_{\infty,c}$ orthogonal to $\partial_x\phi_{\infty,c}$, and define
    \begin{equation*}
        I_\delta^{(d)} := \left\langle\psi_{\delta,c},\phi_{\delta,c}\right\rangle, \quad \delta>\Lambda, \qquad\text{ and }\qquad I_\infty := \left\langle\psi_{\infty,c},\phi_{\infty,c}\right\rangle.
    \end{equation*}
    Then $I_\delta^{(d)}\to I_\infty$ as $\delta\to\infty$.
\end{lemma}
\begin{proof}
    First, $\psi_{\infty,c}$ is well defined, and the set of all $L^2$-solutions of $\mathcal{L}_{\infty,c}\psi=\phi_{\infty,c}$ is $\psi_{\infty,c}+\text{span}\{\partial_x\phi_{\infty,c}\}$: this follows from Step 1 in the proof of Lemma \ref{Lemma: L2 bound for psi in shallow water} applied to $\mathcal{L}_{\infty,c}$, since item 2 of Proposition \ref{Prop: gBO waves are FL ground states} provides for $\mathcal{L}_{\infty,c}$ the spectral information that Lemma \ref{Lemma: Spectral conditions for shallow-gILW operator} provides for $\mathcal{L}_{\delta,c}^{(s)}$.

    The rest of the proof is the same as that of Lemma \ref{Lemma: Convergence of the Albert functional}, with sequences $\delta_n\to\infty$ instead of $\delta_n\to0$, and with $\mathcal{L}_{\delta,c}^{(s)}$, $\mathcal{L}_{0,c}$, $\phi_{\delta,c}^{(s)}$, $\phi_{0,c}$, $\psi_{\delta,c}^{(s)}$, $\psi_{0,c}$ and Lemma \ref{Lemma: L2 bound for psi in shallow water} replaced by $\mathcal{L}_{\delta,c}$, $\mathcal{L}_{\infty,c}$, $\phi_{\delta,c}$, $\phi_{\infty,c}$, $\psi_{\delta,c}$, $\psi_{\infty,c}$ and Lemma \ref{Lemma: L2 bound for psi in deep water}, respectively. The convergence $\phi_{\delta,c}\to\phi_{\infty,c}$ in $L^2$, used in Steps 2 and 4, follows from Theorem \ref{Th: Deep water convergence of solitary waves}, and the uniform convergence used in Step 2 is Theorem \ref{Th: Deep water uniform convergence of solitary waves}, together with Proposition \ref{Prop: Uniform estimates of solitary waves} with $j=d$. Only two points in Step 2 require a comment. First, the analogue of the first summand in \eqref{Eq: strong convergence of the operators on Schwartz functions} is now controlled by Lemma \ref{Lemma: m delta and |xi|}, which gives $0\leq2\pi|\xi|-m_\delta(\xi)\leq2/\delta$ for every $\xi\in\R$, and therefore
    \begin{equation*}
        \left|\left(M_\delta-\mathcal{H}\partial_x\right)v\right|_2\leq\frac{2}{\delta}|v|_2\longrightarrow0, \quad \text{ as }\delta\to\infty,
    \end{equation*}
    for every $v\in\mathcal{S}(\R)$. Second, the limit identity $\mathcal{L}_{\infty,c}\psi_*=\phi_{\infty,c}$, obtained in the sense of distributions, gives $\mathcal{H}\partial_x\psi_*=\phi_{\infty,c}-c\psi_*+\phi_{\infty,c}^p\psi_*\in L^2$, so that $\psi_*\in H^1$, which is the domain of $\mathcal{L}_{\infty,c}$, and the identity holds in $L^2$.
\end{proof}

We can now compute $I_\infty$ and state the main result of this subsection.

\begin{theorem}\label{Th: Stability and instability for gILW in deep water}
    Let $c>0$ and $p\in\N$, $p\neq2$, be fixed. For any $\delta\gg1$, consider $\phi_{\delta,c}$ to be a suitable solitary wave of gILW as above, satisfying condition (i) in Definition \ref{Def: Orbital stability}. Then, there exists some $\Lambda=\Lambda(c,p)>0$ such that:
    \begin{enumerate}
        \item If $p=1$, then, for any $\delta>\Lambda$, $\phi_{\delta,c}$ is orbitally stable.
        \item If $p\geq3$, then, for any $\delta>\Lambda$, $\phi_{\delta,c}$ is unstable.
    \end{enumerate}
\end{theorem}
\begin{proof}
    We apply the stability and instability criterion of Theorem \ref{Theorem: Stability instability criterion Albert}. Its hypotheses on the symbol, namely conditions (C1), (C2) and (C3), are satisfied by $m_\delta$ due to Lemma \ref{Lemma: symbols m s d delta satisfies C1 C2 C3}, and $\phi_{\delta,c}$ satisfies condition (i) in Definition \ref{Def: Orbital stability} by hypothesis. The three spectral conditions on $\mathcal{L}_{\delta,c}$ are given by Lemma \ref{Lemma: Spectral conditions for gILW operator in deep water}, for any $\delta>\Lambda_1$, where $\Lambda_1=\Lambda_1(c,p)>0$ is the constant provided by Lemma \ref{Lemma: L2 bound for psi in deep water}, which is not smaller than the one provided by Lemma \ref{Lemma: Spectral conditions for gILW operator in deep water}. It therefore only remains to determine the sign of $I_\delta^{(d)}$, which, by Lemma \ref{Lemma: L2 bound for psi in deep water}, does not depend on the solution of \eqref{Eq: Equation for psi in deep water} used to compute it.

    \emph{Step 1: computation of $I_\infty$.} By item 1 of Proposition \ref{Prop: gBO waves are FL ground states}, $\phi_{\infty,c}=A\,\tau_cQ$ with $A=((p+1)c)^{1/p}$, that is, $\phi_{\infty,c}(x)=c^{1/p}\phi_{\infty,1}(cx)$, and therefore $\partial_c\phi_{\infty,c}=\frac{A}{c}\,\tau_cR$, where $R=\frac{1}{p}Q+xQ'$. Since $R\in H^2$ by item 2 of Theorem \ref{Theorem: Frank-Lenzmann}, the map $c\mapsto\phi_{\infty,c}$ is differentiable with values in $H^1$, the domain of $\mathcal{L}_{\infty,c}$, and differentiating the gBO solitary wave equation with respect to $c$ gives $\mathcal{L}_{\infty,c}[\partial_c\phi_{\infty,c}]=-\phi_{\infty,c}$. Arguing as in the shallow water regime, we obtain
    \begin{equation*}
        I_\infty = -\frac{1}{2}\frac{d}{dc}\left|\phi_{\infty,c}\right|_2^2 = -\frac{1}{2}\frac{d}{dc}\left(c^{\frac{2}{p}-1}\left|\phi_{\infty,1}\right|_2^2\right),
    \end{equation*}
    that is,
    \begin{equation}\label{Eq: sign of I infinity}
        I_\infty = \frac{p-2}{2p}\,c^{\frac{2}{p}-2}\left|\phi_{\infty,1}\right|_2^2,
    \end{equation}
    so that $I_\infty<0$ if $p=1$, $I_\infty=0$ if $p=2$, and $I_\infty>0$ if $p\geq3$. This reproduces the threshold $p=2\mu$, where $\mu=1$ is the order of the symbol of gBO, given by the Corollary of Theorem 3.1 in \cite{Bona1987}. The value $p=2$, which is the $L^2$-critical exponent for gBO, plays here the role played by $p=4$ in the shallow water regime: it is the only case in which the criterion of Theorem \ref{Theorem: Stability instability criterion Albert} gives no information in the deep water limit.

    \emph{Step 2: conclusion.} Assume $p\neq2$. By \eqref{Eq: sign of I infinity} we have $I_\infty\neq0$, and by Lemma \ref{Lemma: Convergence of the Albert functional in deep water} we have $I_\delta^{(d)}\to I_\infty$ as $\delta\to\infty$. Hence there exists $\Lambda_2=\Lambda_2(c,p)>0$ such that $I_\delta^{(d)}$ and $I_\infty$ have the same sign for any $\delta>\Lambda_2$. Setting $\Lambda:=\max\{\Lambda_1,\Lambda_2\}$, which depends only on $c$ and $p$, we conclude from Theorem \ref{Theorem: Stability instability criterion Albert} that, for any $\delta>\Lambda$: if $p=1$ then $I_\delta^{(d)}<0$ and $\phi_{\delta,c}$ is orbitally stable; while if $p\geq3$ then $I_\delta^{(d)}>0$ and $\phi_{\delta,c}$ is unstable.
\end{proof}

As in the shallow water regime, the conclusions of Theorem \ref{Th: Stability and instability for gILW in deep water} propagate, by Lemma \ref{Lemma: Scaling invariance of the criterion}, to every pair of parameters with the same product $\delta c$.

\begin{corollary}[\textbf{Invariance under the deep water scaling}]\label{Cor: Scaling for stability and instability in deep water}
    Let $c_1>0$ and $p\in\N$, $p\neq2$, let $\Lambda=\Lambda(c_1,p)>0$ be as in Theorem \ref{Th: Stability and instability for gILW in deep water}, and let $\delta_1>\Lambda$ and $\phi_{\delta_1,c_1}$ be as in that theorem. Let $\delta_2,c_2>0$ satisfy $\delta_2c_2=\delta_1c_1$, set $\theta:=\delta_1/\delta_2=c_2/c_1$, and let $\phi_{\delta_2,c_2}:=D_\theta\phi_{\delta_1,c_1}$, which, by items 3 and 5 of Lemma \ref{Lemma: Scaling symmetries}, is a suitable solitary wave of gILW with depth parameter $\delta_2$ and velocity $c_2$. Then $\phi_{\delta_2,c_2}$ satisfies condition (i) in Definition \ref{Def: Orbital stability} and:
    \begin{enumerate}
        \item if $p=1$, then $\phi_{\delta_2,c_2}$ is orbitally stable;
        \item if $p\geq3$, then $\phi_{\delta_2,c_2}$ is unstable.
    \end{enumerate}
\end{corollary}
\begin{proof}
    This follows from Theorem \ref{Th: Stability and instability for gILW in deep water}, whose proof shows that $\mathcal{L}_{\delta_1,c_1}$ satisfies conditions 1, 2 and 3 in Theorem \ref{Theorem: Stability instability criterion Albert} and that $I_{\delta_1}^{(d)}$ is negative if $p=1$ and positive if $p\geq3$, together with Lemma \ref{Lemma: Scaling invariance of the criterion}.
\end{proof}

\begin{remark}\label{Remark: deep water in terms of the product}
    Taking $c_1=1$ and writing $\Lambda_0:=\Lambda(1,p)$, Corollary \ref{Cor: Scaling for stability and instability in deep water} shows that, under the same hypotheses, the conclusions of Theorem \ref{Th: Stability and instability for gILW in deep water} hold for every pair $(\delta,c)$ with $\delta c>\Lambda_0$: what has to be large for the solitary waves of gILW in the deep water regime is the product of the depth parameter and the velocity, and not the depth parameter alone. Together with Corollary \ref{Cor: Stability and instability for gILW}, which concerns pairs with $\delta c$ small, we obtain the following picture for the suitable solitary waves of gILW satisfying condition (i) in Definition \ref{Def: Orbital stability}:
    \begin{center}
    \begin{tabular}{lccccc}
        \hline
        & $p=1$ & $p=2$ & $p=3$ & $p=4$ & $p\geq5$ \\
        \hline
        $\delta c$ small & stable & stable & stable & ? & unstable \\
        $\delta c$ large & stable & ? & unstable & unstable & unstable \\
        \hline
    \end{tabular}
    \end{center}
    In particular, for $p=3$ the solitary waves of gILW are orbitally stable when $\delta c$ is small and unstable when $\delta c$ is large, so that their stability properties change between the two regimes, at values of $\delta c$ which are not reached by the perturbative arguments of this paper. For $p=1$, orbital stability is known for every $\delta,c>0$ \cite{Albert1991, Albert1992}.
\end{remark}

\begin{appendices}

\section{Helly's Selection Theorem}\label{Appendix: Helly selection theorem}
Here we give a proof of an elementary but useful result in real analysis: Helly's selection theorem. This theorem, first introduced in \cite{Helly1912}, is useful in Section \ref{Section: Existence of suitable solitary waves for GILW} to guarantee the existence of the concentration function as a limit, and is a key part of the argument used to show the convergence of solitary waves in Section \ref{Section: Convergence of solitary waves}. We begin showing a previous result.

\begin{lemma}\label{Lemma: Apendix Real analysis 1}
    Let $f_n:[a,b]\to\R$ be a sequence of monotone functions such that $f_n\to h$ point-wise in $[a,b]$. If $h$ is continuous, then $f_n\to h$ uniformly in $[a,b]$.
\end{lemma}
\begin{proof}
    We make the proof for non-decreasing functions but the arguments work for any type of monotonicity.
    
    Let $\epsilon>0$ be given. Since $h$ is uniformly continuous in $[a,b]$, then there exists some $\delta>0$ such that if $|x-y|<\delta$ then $|h(x)-h(y)|<\epsilon/2$.

    Now, consider the partition of the interval $[a,b]$ given by
    \begin{equation*}
        a=x_0<x_1<\cdots<x_n=b,
    \end{equation*}
    where $|x_i-x_{i+1}|<\delta$ for any $i$.

    Since $f_n$ converge point-wise, there exists some $N\in\N$ such that for any $n\geq N$, $|f_{n}(x_i)-h(x_i)|<\epsilon/2$ for any $i$. 

    For any $x\in[a,b]$, there are some $x_i,x_{i+1}$ such that $x_i\leq x \leq x_{i+1}$. We have
    \begin{equation*}
        f_n(x_i) \leq f_n(x) \leq f_n(x_{i+1})
    \end{equation*}
    and 
    \begin{equation*}
        h(x_i) \leq h(x) \leq h(x_{i+1}),
    \end{equation*}
    where we have used that the limit of monotone functions is also monotone.

    Now
    \begin{align*}
        f_n(x)-h(x) &\leq f_n(x_{i+1})-h(x_i) \\
        &\leq |f_n(x_{i+1})-h(x_{i+1})| + |h(x_{i+1})-h(x_i)| \\
        &< \epsilon,
    \end{align*}
    and
    \begin{align*}
        h(x)-f_n(x) &\leq h(x_{i+1})-f_n(x_i) \\
        &\leq |h(x_{i+1})-h(x_{i})| + |h(x_{i})-f_n(x_i)| \\
        &< \epsilon.
    \end{align*}
    Therefore, $|h(x)-f_n(x)|<\epsilon$ for any $n\geq N$.
\end{proof}

Now we prove Helly's selection theorem.

\begin{theorem}(\textbf{Helly's selection})\label{Theorem: Helly selection theorem}
    Let $f_n:[0,\infty)\to [0,\infty)$ be a sequence of monotone functions such that $|f_n|_{\infty}\leq B$ for all $n$. Then, there exists some monotone function $h$, with the same type of monotonicity, and a sub-sequence $f_{n_j}$ such that $f_{n_j}\to h$ point-wise in $[0,\infty)$. Moreover, if $h$ is continuous, then $f_{n_j}\to h$ uniformly in compact sets.
\end{theorem}
\begin{proof}
    We make the proof for non-increasing functions, but the arguments work for any kind of monotone functions.

    Note that if we prove the convergence $f_{n_j}\to h$ point-wise in $[0,\infty)$, then the last part of the theorem about uniform convergence in compact sets is a direct application of Lemma \ref{Lemma: Apendix Real analysis 1}.
    
    Let $A=\{a_k\}_{k\in\N}$ be an enumeration of $\mathbb{Q}\cap[0,\infty)$. Using a Cantor diagonal argument and the Bolzano-Weierstrass theorem, we can get a sub-sequence $f_{n_j}$ such that for all $a_k\in A$: $f_{n_j}(a_k)$ is convergent. 

    Define the function $h_0: A\to[0,B]$ as $h_0(a_k):=\lim_j f_{n_j}(a_k)$. Note that $h_0$ is non-increasing. Then, we extend this function to $(0,\infty)$ as 
    \begin{equation*}
        h(x) := \inf\{h_0(a_k): a_k\in A, a_k\leq x\},
    \end{equation*}
    and set $h(0):= h_0(0)$. It follows from the above definition that $h$ is non-increasing and that $f_{n_j}(a_k)\to h(a_k)$.

    Lets see that for any $x_0\in[0,\infty)$ where $h$ is continuous, we have $f_{n_j}(x_0)\to h(x_0)$. Take $a,a'\in A$ such that $a<x_0<a'$. Then, for any $j$,
    \begin{equation*}
        f_{n_j}(a) \geq f_{n_j}(x_0) \geq f_{n_j}(a'),
    \end{equation*}
    from where we get
    \begin{equation*}
        h(a)\geq\limsup_j f_{n_j}(x_0)\geq\liminf_j f_{n_j}(x_0) \geq h(a').
    \end{equation*}
    Since $h$ is continuous in $x_0$, we can take the limits $a\to x_0^{-}$ and $a\to x_0^{+}$ to finally obtain $f_{n_j}(x_0)\to h(x_0)$.

    We are only left to prove the point-wise convergence $f_{n_j}(x)\to h(x)$ in points $x$ where $h$ is discontinuous. Since $h$ is monotone, there is at most an enumerable set of discontinuities $D$, therefore, we can apply again a Cantor diagonal argument to extract a new sub-sequence, which we still denote as $f_{n_j}$, such that $f_{n_j}(x)$ is convergent for any $x\in[0,\infty)$. 
    
    Finally, changing the values of $h$ (if necessary) at the points $x\in D$ we get a new function $h$ which is the point-wise limit of $f_{n_j}$, and since the point-wise limit of monotone functions is monotone, we conclude the proof.
\end{proof}

\section{Symmetric decreasing rearrangement}\label{Appendix: Symmetric decreasing rearrangement}
Here we define the basic notions and stablish some results on symmetric decreasing rearrangements. These results are used in Section \ref{Section: Existence of suitable solitary waves for GILW} to prove the existence of suitable solitary waves for gILW. For more details on this topic see \cite{Lieb2001}.

Let
\begin{equation*}
    S:=\left\{f:\R\to[0,\infty]\hspace{0.1cm} \Big| \text{ } f \text{ is even, } f(x)\leq f(y) \text{ if } |x|\geq|y| \right\}
\end{equation*}
be the set of measurable functions which are even and radially decreasing.

For any non-negative measurable function $f$ on $\R$, define
\begin{equation*}
    \Gamma_y(f):=\{x\hspace{0.1cm}| f(x)\geq y\} \quad \text{and} \quad \Theta_y(f):=\mu(\Gamma_y(f)),
\end{equation*}
where $\mu$ is the Lebesgue measure. If $\Theta_a(f)<\infty$ for some $a<\infty$, then, any $f^*\in S$ and satisfies
\begin{equation*}
    \Theta_y(f^*)=\Theta_y(f), \quad \forall y>0,
\end{equation*}
is called a symmetric decreasing rearrangement of $f$, or rearrangement, for short. 

If $f:\R\to\mathbb{C}$ is a measurable function, we define $f^*:=|f|^*$ whenever $\Theta_a(|f|)<\infty$ for some $a<\infty$. The following are some useful properties of rearrangements.

\begin{proposition}\label{Prop: Rearrangement properties}
    Let $f:\R\to\mathbb{C}$ be a measurable function such that $\Theta_a(|f|)<\infty$ for some $a<\infty$, then,
    \begin{enumerate}
        \item $|f^*|_p=|f|_p$, for any $p\in[1,\infty]$.
        \item $(f^*)^r=(f^r)^*$ point-wise for all $r>0$.
        \item If $f,g\in L^2$, then
        \begin{equation*}
            \Big|\int_{\R}f(x)g(x)dx\Big| \leq \int_\R f^*(x)g^*(x)dx.
        \end{equation*}
        \item \textbf{(Riesz rearrangement lemma)} Let $f_1,\cdots,f_n$ be measurable functions on $\R$ such that they satisfy $\mu(\{x| f_i(x)\geq y\})<\infty$ for all $y>0$ and all $1\leq i \leq n$. Then,
        \begin{equation*}
            |(f_1*\cdots*f_n)(0)| \leq (f^*_1*\cdots *f^*_n)(0)
        \end{equation*}
        in the sense that if the right-hand side is finite, then the left-hand side exists and the inequality holds.
    \end{enumerate}
\end{proposition}

\begin{lemma}\label{Lemma rearrangement in H1/2}
    Let $\psi\in H^{1/2}$, then $|\psi|,\psi^*\in H^{1/2}$. Moreover,
    \begin{equation*}
        \norm{D^{1/2}|\psi|}^2 \leq \norm{D^{1/2}\psi}^2 \quad \text{ and } \quad \norm{D^{1/2}\psi^*}^2 \leq \norm{D^{1/2}\psi}^2.
    \end{equation*}
\end{lemma}
\begin{proof}
    See Lemma $8.15$ in \cite{Angulo2009}.
\end{proof}

\end{appendices}

\section*{Acknowledgments}
The author thanks the Instituto de Matemática Pura e Aplicada (IMPA) for its support and for providing excellent working conditions and facilities during the development of this work. The author also acknowledges the financial support of the Conselho Nacional de Desenvolvimento Científico e Tecnológico (CNPq), Brazil.
\bibliographystyle{amsplain}
\bibliography{referencias}

\end{document}